\documentclass[a4paper,10pt]{article}

\usepackage{bm}
\usepackage{a4wide}
\usepackage{amsmath,amssymb,amsthm}
\usepackage{subfig}
\usepackage{graphicx}
\usepackage[numbers]{natbib}

\usepackage{commath}

\usepackage{hyperref}
\usepackage{url}

\usepackage{titlesec}
\titleformat*{\section}{\bf\large}
\titleformat*{\subsection}{\bf\normalsize}
\titleformat*{\subsubsection}{\itshape}
\titleformat*{\paragraph}{\itshape}

\newtheorem{remark}{Remark}[section]

\begin{document}
\title{Optimal control with stochastic PDE constraints and uncertain controls}

\author{Eveline Rosseel$^\star$ \and Garth N. Wells$^\dagger$}

\date{\normalsize
      $^\star$Department of Computer Science,
      Katholieke Universiteit Leuven, Belgium \\
      \url{eveline.rosseel@cs.kuleuven.be}
      \\[1em]
      $^\dagger$Department of Engineering,
      University of Cambridge,
      United Kingdom \\
      \url{gnw20@cam.ac.uk}}
\maketitle
\begin{abstract}
\noindent The optimal control of problems that are constrained by
partial differential equations with uncertainties and with uncertain
controls is addressed. The Lagrangian that defines the problem
is postulated in terms of stochastic functions, with the control
function possibly decomposed into an unknown deterministic component
and a known zero-mean stochastic component. The extra freedom provided
by the stochastic dimension in defining cost functionals is explored,
demonstrating the scope for controlling statistical aspects of the system
response. One-shot stochastic finite element methods are used to find
approximate solutions to control problems. It is shown that applying the
stochastic collocation finite element method to the formulated problem leads to
a coupling between stochastic collocation points when a deterministic
optimal control is considered or when moments are included in the cost
functional, thereby forgoing the primary advantage of the collocation
method over the stochastic Galerkin method for the considered problem.
The application of the presented methods is demonstrated through a number
of numerical examples. The presented framework is sufficiently general
to also consider a class of inverse problems, and numerical examples of
this type are also presented.
\end{abstract}
\section{Introduction}

In many applications, forces or boundary conditions are to be determined
such that the response of a physical or engineering system is optimal in
some sense. These problems can often be formulated as the minimisation
of an objective functional subject to a set of constraint equations
in the form of partial differential equations (PDEs). For problems
that involve uncertainty, incorporating stochastic information into a
control formulation can lead to a quantification of the statistics of
the system response.  Moreover, there is scope to control a system not
only for an optimal mean response, but also to include statistics of the
response in a cost functional.  We note that stochastic PDE-constrained
optimisation problems are closely related to stochastic inverse
problems, where the control variable corresponds to the parameter to be
identified~\citep{Zabaras2008, Jin2008, Sankaran2009}.

We examine in this work the numerical solution of optimal control
problems constrained by stochastic PDEs and with uncertain controls.
Stochastic finite element-based solvers have been studied extensively
for a large range of stochastic PDEs~\citep{Xiu2009}. However, few
results and examples on solving optimisation problems constrained by
stochastic PDEs are available. Existence of an optimal solution to
stochastic optimal control problems constrained by stochastic elliptic
PDEs was studied by \citet{Hou2010} for deterministic control functions.
\citet{Borzi2009} and \citet{Borzi2010a} studied multigrid solvers for
stochastic collocation solutions of parabolic and elliptic optimal control
problems with random coefficients and stochastic control functions.
We contend that problems with an unknown stochastic `control function'
constitute stochastic inverse problems and are different from control
problems where the focus is on computing a deterministic component of the
control function which forms the control `signal'.  If the stochastic
properties of the control are computed, \emph{ad hoc} procedures are
required to extract a deterministic function, which will in general
not be the optimal control.  We choose to split the control function
into an unknown deterministic component, which is to be computed, and a
known zero-mean stochastic part that represents the uncertainty in the
controller response.

Control cost functionals will be formulated in terms of norms that
include both spatial and stochastic dimensions. The inclusion of the
stochastic dimension provides additional freedom in the definition
of cost functionals. We formulate a one-shot approach to control,
and solve the resulting equations via a stochastic collocation or a
Galerkin finite element method. The one-shot approach is in contrast
with methods that require iteration~\citep{Gunzburger2003}.  Overviews of
the stochastic collocation and Galerkin finite element methods can
be found in \citep{Babuska2010} and \citep{Babuska2004}, respectively.
The stochastic collocation method is often preferred over the Galerkin
approach as it converts a stochastic problem into a collection of
decoupled deterministic problems.  However, we will show that the
so-called non-intrusivity property of the collocation method that is often
exploited, including for stochastic inverse problems~\citep{Zabaras2008,
Borzi2009}, does not hold for a large class of stochastic PDE-constrained
optimisation problems. The stochastic Galerkin method, on the other
hand, can be applied straightforwardly to stochastic optimal control
problems.  The efficient solution of stochastic finite element problems,
and in particular for the stochastic Galerkin method, can hinge on the
development and application of effective preconditioners.  This aspect is
addressed for stochastic control problems, with two preconditioners that
take the specific structure of the Galerkin one-shot systems into account
presented.  Extensive numerical examples support our investigations,
and the computer code to reproduce all numerical results is available
under the GNU Lesser Public License (LGPL) as part of the supporting
material~\citep{supporting}.

The remainder of this paper is organised as
follows. Section~\ref{sec:model} presents a control problem
involving an elliptic stochastic PDE constraint and formulates
this problem as a coupled system of stochastic PDEs. The framework
developed in Section~\ref{sec:model} is formulated such that it is
sufficiently general to also address a class of inverse problems.
The stochastic finite element discretisation is presented in
Section~\ref{sec:sfem}. Section~\ref{sec:regul} considers additional
regularising terms and their impact in the context of stochastic finite
element solvers.  Iterative solvers and preconditioners for the one-shot
Galerkin system are discussed in Section~\ref{sec:iterative}, which is
followed in Section~\ref{sec:num1} by numerical examples of stochastic
optimal control problems. Numerical examples illustrating the solution
of stochastic inverse problems are given in Section~\ref{sec:num2},
and conclusions are drawn in Section~\ref{sec:conclu}.
\section{A control problem with stochastic PDE constraints}
\label{sec:model}

We consider optimal control problems constrained by partial differential
equations with stochastic coefficients. In general, these can be
formulated as:
\begin{equation}
  \min_{z\in Z, \: u\in U} \mathcal{J}(z,u) \quad \text{subject to} \quad
  c(z,u) = 0,
\label{eq:control_constr}
\end{equation}
where $\mathcal{J} : Z \times U \rightarrow \mathbb{R}$ is a \emph{cost}
(or \emph{objective}) functional, $c$ is a constraint, $z$ the
\emph{state} variable, $u$ the \emph{control} variable, and $Z$ and $U$
are suitably defined function spaces.

In this work, the constraint equations in~\eqref{eq:control_constr} will
be given by one or more PDEs with stochastic coefficients. The state
variable will be a stochastic function and the control can be either
deterministic or stochastic. Although the cost functional $\mathcal{J}$
contains stochastic functions, it will be defined such that its outcome
is deterministic.  Before presenting concrete examples, it is useful to
provide some definitions.
\subsection{Definitions}

We will consider problems posed on a polygonal spatial domain $D
\subset \mathbb{R}^{d}$ with boundary~$\partial D$ and where $1 \le d
\le 3$.  The boundary is decomposed into $\partial D_{D}$ and $\partial
D_{N}$ such that $\partial D_{D} \cup \partial D_{N} = \partial D$
and $\partial D_{D} \cap \partial D_{N} = \emptyset$. The outward unit
normal vector to $\partial D$ is denoted by~$n$.  Consider also a complete
probability space $(\Omega,\mathcal{F},\mathcal{P})$, where $\Omega$ is
a sample space, $\mathcal{F}$ is a $\sigma$-algebra and $\mathcal{P}$
is a probability measure. We define the tensor product Hilbert space
$L^2(D)\otimes L^2(\Omega)$ of second-order random fields as
\begin{equation}
 L^2(D)\otimes L^2(\Omega)
    = \cbr{z : D \otimes \Omega \rightarrow \mathbb{R} ,\:
    \int_\Omega\int_D |z|^2  \dif x \dif \mathcal{P} < \infty  }.
\end{equation}
This space is equipped with the norm
\begin{equation}
 \norm{z}_{L^{2}(D)\otimes L^{2}(\Omega)}
    = \del{\int_\Omega \int_D |z|^2 \dif x \dif \mathcal{P}}^{\frac{1}{2}}.
\end{equation}
Analogously, the tensor product spaces $H^1_0(D)\otimes L^2(\Omega)$
and $H^1(D)\otimes L^2(\Omega)$ can be defined~\citep{Babuska2004}.
\subsection{Model problem}
\subsubsection{Constraint equation}

As a model constraint, we consider a stochastic steady-state diffusion
equation:
\begin{align}
\label{eq:constr1}
 -\nabla_{x} \cdot \del{\kappa\del{x, \omega} \nabla_{x} z\del{x,\omega}}
      &= u\del{x, \omega} \qquad &x \in D, \ \omega \in \Omega,
\\
 \label{eq:constr2}
  z\del{x,\omega} &= z_D\del{x,\omega}
        \qquad &x \in \partial D_{D}, \ \omega \in \Omega,
\\
  \label{eq:constr3}
 \kappa\del{x,\omega} \nabla_{x} z\del{x,\omega}\cdot n &= g\del{x,\omega}
        \qquad &x \in \partial D_{N}, \ \omega \in \Omega,
\end{align}
where the function $z : D \times \Omega \rightarrow \mathbb{R}$
is to be found, the diffusion coefficient $\kappa : D \times \Omega
\rightarrow \mathbb{R}$ is prescribed, $u : D \times \Omega \rightarrow
\mathbb{R}$ is a source term, $z_{D} : \partial D_{D} \times \Omega
\rightarrow \mathbb{R}$ is a Dirichlet condition and $g : \partial
D_{N} \times \Omega \rightarrow \mathbb{R}$ is a Neumann condition.
The operator~$\nabla_{x}$ involves only derivatives with respect to the
spatial variable~$x$.  If $\kappa$, $u$, $z_{D}$ and $g$ are second-order
random fields and assuming boundedness and positivity of the diffusion
coefficient $\kappa$, existence and uniqueness of a weak solution $z$
can be proved~\citep{Babuska2004}.

In practice, for control problems $u$ and/or $g$ will not be prescribed,
but computed as part of a constrained optimisation problem.
\subsubsection{Cost functionals}

We consider two cost functionals which, while outwardly similar,
account differently for the stochastic nature of the problem.  The first
functional has the form:
\begin{multline}
\label{eq:cost}
  \mathcal{J}_{1}\del{z, u, g}
    := \frac{\alpha}{2} \norm{z - \Hat{z}}^2_{L^2(D)\otimes L^2(\Omega)}
        + \frac{\beta}{2} \norm{\mathrm{std}(z)}^2_{L^2(D)}
        \\
        + \frac{\gamma}{2}\norm{u}^2_{L^2(D)\otimes L^2(\Omega)}
        + \frac{\delta}{2} \norm{g}^2_{L^2(\partial D_N)\otimes L^2(\Omega)},
\end{multline}
where $\Hat{z}: D\times \Omega\rightarrow \mathbb{R}$ is a prescribed
target function and $\alpha$, $\beta$, $\gamma$ and $\delta$ are
positive constants.  Typically for control problems $\Hat{z}$ will be
deterministic, but we permit here a more general case.  The first term in
\eqref{eq:cost} is a measure of the distance, between the state variable
$z$ and the prescribed function $\Hat{z}$, in terms of the expectation of
$\del{z - \Hat{z}}^{2}$.  The second term measures the standard deviation
of~$z$, which is denoted by~$\mathrm{std}(z)$ and defined by:
\begin{equation}
\mathrm{std}(z):= \del{\int_\Omega \del{z
 - \int_\Omega z \dif \mathcal{P}}^2\dif \mathcal{P}}^{\frac{1}{2}}.
\end{equation}
By increasing the value of $\beta$ with respect to $\alpha$, a greater
relative contribution of the variance of $z$ to $\mathcal{J}_{1}$ is
implied. The final two terms in~\eqref{eq:cost} are regularisation terms
for the distributive control via~$u$ and the boundary control via~$g$.

The second considered functional has the form:
\begin{multline}
\label{eq:cost_alt}
  \mathcal{J}_{2}\del{z, u, g}
    := \frac{\alpha}{2} \norm{  \bar{z} - \hat{z}}^2_{ L^2(D)}
        + \frac{\beta}{2} \norm{\mathrm{std}(z)}^2_{L^2(D)}
        + \frac{\gamma}{2}\norm{u}^2_{L^2(D)\otimes L^2(\Omega)}
        + \frac{\delta}{2} \norm{g}^2_{L^2(\partial D_N)\otimes L^2(\Omega)},
\end{multline}
where $\Bar{z}$ is the expectation of~$z$ and $\hat{z} :D \rightarrow
\mathbb{R}$ is a prescribed target function. This functional is subtly,
but significantly, different from that in~\eqref{eq:cost}. The first
term in $\mathcal{J}_{2}$ measures the $L^{2}$-distance between the
expectation of $z$ and the target $\Hat{z}$. This includes no measure of
the variance of the actual response and is unaffected by the presence
of uncertainty in~$z$.

Setting $\alpha = 1$, $\beta = 0$ and $\delta = 0$, the functional
$\mathcal{J}_{1}$ in~\eqref{eq:cost} coincides with a functional presented in
\citet{Borzi2010a}, but in their work a method is constructed around a
modified version of~$\mathcal{J}_{1}$.  For the above parameters and for
deterministic $u$, $\mathcal{J}_{1}$ coincides with the cost functional
considered by \citet{Hou2010}.

\subsubsection{Structure of the control functions}

A feature of our work is that for control problems, as opposed
to inverse problems, we will consider control variables that are
decomposed additively into unknown deterministic (to be computed) and
known stochastic components.  We will consider $u$ to have the form
\begin{equation}
\label{eq:u}
 u\del{x,\omega} = \Bar{u}\del{x} + u'\del{x, \omega},
\end{equation}
where $\Bar{u} : D \rightarrow \mathbb{R}$ is deterministic and is
the mean of $u$ and $u': D \times \Omega \rightarrow \mathbb{R}$ is a
zero-mean stochastic part. The goal will be to compute $\bar{u}$, which
constitutes the `signal' sent to a control device. The actual controller
response is $u$, with $u'$ modelling the uncertainty in the controller
response for a given instruction.  The boundary control function will
be decomposed analogously,
\begin{equation}
\label{eq:g}
 g\del{x,\omega} = \Bar{g}\del{x} + g'\del{x, \omega},
\end{equation}
where $\Bar{g} : \partial D_{N} \rightarrow \mathbb{R}$ and $g':
\partial D_{N} \times \Omega \rightarrow \mathbb{R}$ is a zero-mean
stochastic part.
\subsection{Representation of stochastic fields}
\subsubsection{Finite-dimensional noise assumption}

In representing random fields, we employ the finite-dimensional
noise assumption \citep[Section 2.4]{Babuska2004}, which states that
the random fields $\kappa$, $z_D$, $g$ and $u$ can be approximated
using a prescribed finite number of random variables~$\xi = \cbr{
\xi_{i} }_{i=1}^{L}$, where $L \in \mathbb{N}$ and $\xi_i : \Omega
\rightarrow \Gamma_i \subseteq \mathbb{R}$.  We assume that each random
variable is independent and is characterised by a probability density
function~$\rho_i : \Gamma_i \rightarrow [0,1]$.  Defining the support
$\Gamma = \prod_{i=1}^L \Gamma_i \subset \mathbb{R}^L$, for a given $y =
(y_1, \hdots, y_L) \in \Gamma$ the joint probability density function
of $\xi$ is given by $\rho = \prod_{i=1}^L \rho_i(y_i)$.  The preceding
assumptions enable a parametrisation of the problem in $y$ in place of
the random events~$\omega$~\citep{Frauenfelder2005}.

As an example, consider a finite-term expansion of the stochastic
coefficient~$\kappa$ based on $L$ random variables:
\begin{equation}\label{eq:kappa}
 \kappa(x, y) = \sum_{i = 1}^{S} \kappa_{i}(x) \zeta_{i}(y)
\qquad x \in D, \ y\in \Gamma,
\end{equation}
where $\kappa_{i}: D \rightarrow \mathbb{R}$ and $\zeta_{i}: \Gamma
\rightarrow \mathbb{R}$.  If $\kappa$ is represented by a truncated
Karhunen--Lo\`eve expansion~\cite{Loeve1977}, then $S = L + 1$ with
$\zeta_{i} = y_{i-1}$ and $y_{0} = 1$.  If a generalised polynomial
chaos expansion~\cite{Xiu2002a} is used, $\zeta_i$ is an $L$-variate
orthogonal polynomial of order $p$ and $S = (L +p)!/(L! p!)$.

\subsubsection{Definitions}

Given the joint probability density function $\rho(y)$ of $\xi$,
with $y \in \Gamma$ and where $\Gamma$ is the support of~$\rho$, a
space~$L^2_\rho(\Gamma)$ equipped with the inner product
\begin{equation}\label{eq:innerL2}
 (v,w)_{L^{2}_{\rho}(\Gamma)} := \int_\Gamma v w \rho \dif y
\qquad v,w \in  L^2_\rho(\Gamma)
\end{equation}
is considered.  The norm $\norm{\cdot}_{L^2(D)\otimes L^2_\rho(\Gamma)}$
is then defined as:
\begin{equation}
  \norm{z}_{L^2(D)\otimes L^2_\rho(\Gamma) }
      := \del{\int_{\Gamma} \int_D |z|^2 \rho \dif x \dif y}^{\frac{1}{2}}.
\end{equation}
Expressing two functions $u_{1}$ and $u_{2}$ as $u_{1}(x, y)
= v_{1}(x)w_{1}(y)$ and $u_{2}(x, y) = v_{2}(x)w_{2}(y)$, where $v_{1},
v_{2} \in H^{1}(D)$ and $w_{1}, w_{2} \in H^{1}_{\rho}(\Gamma)$, an inner
product on $H^{1}(D) \otimes H^{1}_{\rho}(\Gamma)$ is defined by
\begin{equation}\label{eq:innerH1}
 \del{u_{1}, u_{2} }_{H^{1}(D) \otimes H^{1}_{\rho}(\Gamma)}
    := \del{v_{1}, v_{2}}_{H^{1}(D)} \del{w_{1}, w_{2} }_{H^{1}_{\rho}(\Gamma)},
\end{equation}
where $\del{\cdot, \cdot}_{H^{1}(D)}$ is the standard $H^{1}$ inner
product and $\del{\cdot, \cdot}_{H^{1}_{\rho}}$ is the $H^{1}$ inner product
weighted by $\rho$, analogous to~\eqref{eq:innerL2}.  The norm
$\norm{\cdot}_{H^{1}(D)\otimes H^{1}_{\rho}(\Gamma)}$ is that induced
by~\eqref {eq:innerH1}. The $H^{1}$ definition will be used in
Section~\ref{sec:regul} when considering the impact of introducing extra
regularisation terms into the cost functionals.
\subsubsection{Parametric optimal control problem}

Following from the Doob--Dynkin Lemma~\citep{Babuska2004}, the
random field~$z$ can be expressed as a function of the given~$L$
random variables. This enables one to reformulate the stochastic
optimal control problem corresponding to the cost functional
in~\eqref{eq:cost} or~\eqref{eq:cost_alt} and the constraints in
\eqref{eq:constr1}--\eqref{eq:constr3} as a parametric PDE-constrained
optimisation problem. The parametric problem involves solving
\begin{equation}
 \label{eq:miny}
\min_{z,u,g} \mathcal{J}_1\del{z, u, g}
\quad \text{or} \quad
 \min_{z,u,g} \mathcal{J}_2 \del{z, u, g}
\end{equation}
subject to:
\begin{align}
 \label{eq:constr1y}
 -\nabla_{x} \cdot \del{\kappa(x,y) \nabla_{x} z(x,y)}
    &=u(x,y) \quad &x \in D, \ y \in \Gamma,
\\
\label{eq:constr2y}
  z(x,y) &= z_D(x,y) \quad &x \in \partial D_{D}, \ y \in \Gamma,
\\
\label{eq:constr3y}
 \kappa(x,y) \nabla_{x} z(x,y)\cdot n
    &=  g(x,y)\quad &x \in  \partial D_{N}, \ y \in \Gamma,
\end{align}
where
\begin{multline}
\label{eq:costy}
 \mathcal{J}_1\del{z, u, g}
    := \frac{\alpha}{2} \norm{z - \hat{z}}^2_{L^2(D)\otimes L^2_\rho(\Gamma)}
      + \frac{\beta}{2}\norm{\mathrm{std}(z)}^2_{L^2(D)}
      \\
      + \frac{\gamma}{2} \norm{u}^2_{L^2(D)\otimes L^2_\rho(\Gamma)}
      + \frac{\delta}{2} \norm{g}^2_{L^2(\partial D_N)\otimes L^2_\rho(\Gamma)}
\end{multline}
and
\begin{multline}
\label{eq:costy_alt}
 \mathcal{J}_2 \del{z, u, g}
    := \frac{\alpha}{2} \norm{  \bar{z} - \hat{z}}^2_{ L^2(D)}
      + \frac{\beta}{2} \norm{\mathrm{std}(z)}^2_{L^2(D)}
      + \frac{\gamma}{2}\norm{u}^2_{L^2(D)\otimes L^2_\rho(\Gamma)}
      + \frac{\delta}{2} \norm{g}^2_{L^2(\partial D_N)\otimes L^2_\rho(\Gamma)}.
\end{multline}
This parametric form of the problem will be considered in the remainder.

\subsection{One-shot solution approach}
\label{ssec:one_shot}

The PDE-constrained optimisation problem
in~\eqref{eq:miny}--\eqref{eq:constr3y} can be recast as an unconstrained
optimisation problem by following the standard procedure of introducing
Lagrange multipliers~\citep{Gunzburger2003,Troltzsch2010}. By defining
a Lagrangian functional, optimality conditions can be formulated for
the state, control and adjoint variables~\citep{Troltzsch2010} and
these can be solved simultaneously.  This so-called `one-shot' approach
yields a solution for an optimal control problem without having to apply
iterative optimisation routines.
\subsubsection{Lagrangian functionals}

Introducing adjoint variables, or Lagrange multipliers, $\lambda, \chi
\in H_{0}^{1}(D) \otimes L^{2}_{\rho}(\Gamma)$, we define a Lagrangian
associated with the cost functional $\mathcal{J}_{1}$ in~\eqref{eq:costy}
and the constraints \eqref{eq:constr1y}--\eqref{eq:constr3y} as:
\begin{multline}
 \mathcal{L}_{1}(z, u, g, \lambda, \chi)
  :=
   \frac{\alpha}{2} \int_\Gamma \int_D  \del{z-\hat{z}}^2  \rho \dif x \dif y
  + \frac{\beta}{2} \int_\Gamma \int_D  z^2 \rho \dif x \dif y
  - \frac{\beta}{2} \int_D \del{\int_\Gamma z \rho \dif y}^2 \dif x
\\
+ \frac{\gamma}{2} \int_\Gamma \int_D  u^2 \rho \dif x\dif y
  + \frac{\delta}{2} \int_\Gamma \int_{\partial D_N} g^2 \rho \dif s \dif y
  -  \int_\Gamma \int_D \lambda \del{-\nabla_{x} \cdot \del{\kappa \nabla_{x} z}-u} \rho \dif x \dif y
\\
  -  \int_\Gamma \int_{\partial D_N}\chi \del{\kappa \frac{dz}{dn}-g}\rho \dif s \dif y.
\label{eq:lagrangian1}
\end{multline}
A Lagrange multiplier has not been introduced to impose the Dirichlet
boundary condition in~\eqref{eq:constr2y} since this condition will
be imposed by construction via the definition of the function space to
which~$z$ belongs.  For the problem associated with the cost functional
$\mathcal{J}_{2}$ in~\eqref{eq:costy_alt}, we define the Lagrangian
\begin{multline}
 \mathcal{L}_{2}(z, u, g, \lambda, \chi)
  := \frac{\alpha}{2} \int_D \del{\int_\Gamma  z \rho \dif y -\hat{z}}^2  \dif x
   + \frac{\beta}{2} \int_\Gamma \int_D  z^2 \rho \dif x \dif y
   - \frac{\beta}{2} \int_D \del{\int_\Gamma z \rho \dif y }^2 \dif x
\\+  \frac{\gamma}{2} \int_\Gamma \int_D u^2 \rho \dif x \dif y
   +\frac{\delta}{2}  \int_\Gamma \int_{\partial D_N}g^2 \rho \dif s \dif y
   -  \int_\Gamma \int_D\lambda
   \del{-\nabla_{x} \cdot \del{\kappa \nabla_{x} z}-u} \rho \dif x \dif y
 \\
   -\int_\Gamma \int_{\partial D_N}\chi  \del{ \kappa \frac{dz}{dn}-g} \rho\dif s \dif y.
\label{eq:lagrangian2}
\end{multline}

The control problems that we wish to solve involve finding stationary
points of these Lagrangians.  The existence of Lagrange multipliers
for stochastic optimal control problems is proved by \citet{Hou2010}
for $\mathcal{J}_{1}$ with $\alpha =1$, $\beta =0$ and $\delta = 0$. The
result extends to the more general form of $\mathcal{J}_{1}$ because of
the Fr\'echet differentiability of $\norm{\mathrm{std}(z)}^2_{L^2(D)}$
and $\norm{g}^2_{L^2(\partial D_N)\otimes L^2_\rho(\Gamma)}$.

\subsubsection{Optimality system}
\label{sssec:optsys}

To find stationary points of the Lagrangians in~\eqref{eq:lagrangian1}
and~\eqref{eq:lagrangian2}, we consider variations with respect to the
adjoint, state and control variables.  This will lead to the first-order
optimality conditions, which are known as the state, adjoint and
optimality system of equations~\citep{Gunzburger2003}. In the remainder of
this section, these equations are derived by setting the directional
derivative of the Lagrangians with respect to the adjoint, state and control
variable, respectively, equal to zero.

Following standard variational arguments, taking the directional
derivative of~$\mathcal{L}_{1}$ or of~$\mathcal{L}_{2}$ with respect
to $\chi$ and setting this equal to zero for all variations leads to
the recovery of the Neumann boundary condition in~\eqref{eq:constr3y}.
Likewise, taking directional derivatives with respect to $\lambda$ and
setting this equal to zero for all variations leads to the recovery of
the constraint equation in~\eqref{eq:constr1y}. The state system of equations
corresponds thus to the constraint
equations~\eqref{eq:constr1y}--\eqref{eq:constr3y}.

The derivative of~$\mathcal{L}_{1}$ with respect to the state variable
$z$ in the direction of $z^{\star} \in H^1_0(D) \otimes L^2_\rho(\Gamma)$
reads:
\begin{multline}
\label{eq:Lu_poiss}
 \mathcal{L}_{1,z}\sbr{z^{\star}}
  =  \alpha \int_\Gamma \int_D  (z-\hat{z})z^{\star}\rho \dif x \dif y
  +  \beta \int_\Gamma \int_D  z z^{\star} \rho \dif x \dif y
  -\beta \int_\Gamma\int_D  \del{\int_\Gamma z\rho \dif y}
 z^{\star} \rho \dif x\dif y
\\ +  \int_\Gamma \int_D \lambda \nabla_{x} \cdot (\kappa \nabla_{x} z^{\star})
\rho \dif x \dif y
 - \int_\Gamma
 \int_{\partial D_N} \chi  \kappa \frac{d z^{\star}}{d n}\rho \dif s \dif y.
\end{multline}
Setting $\mathcal{L}_{1,z}\sbr{z^{\star}} = 0$ for all $z^{\star}$ and
following standard arguments, the following adjoint system of equations
can be deduced:
\begin{align}
\label{eq:adj1}
 - \nabla_{x} \cdot( \kappa(x,y) \nabla_{x}\lambda(x,y)) &= (\alpha+\beta)z(x,y)
    -\alpha \hat{z}(x,y)
        - \beta \int_\Gamma z \rho \dif y  &x \in D,\: y\in \Gamma,
\\
 \label{eq:adj2}
\kappa(x,y) \nabla_{x} \lambda(x,y) \cdot n
  &= 0 &x \in \partial D_N,\: y\in \Gamma,
\\
\label{eq:adj3}
\lambda(x,y)  &= 0  &x \in \partial D_D,\: y\in \Gamma,
\\
\label{eq:adj4}
\lambda(x,y) &= \chi(x,y) &x \in \partial D_N,\: y\in \Gamma.
\end{align}
The last equation $\lambda = \chi$ is trivial.  For distributive control
via $u$ the directional derivative of $\mathcal{L}_{1}$ with respect to
$u$ in the direction of $u^{\star}\in L^2(D)\otimes L^{2}_{\rho}(\Gamma)$
reads:
\begin{equation}
\label{eq:opt_poiss}
 \mathcal{L}_{1, u}\sbr{u^{\star}}
    = \int_\Gamma \int_D (\gamma u + \lambda) u^{\star} \rho \dif x \dif y.
\end{equation}
Setting the above equal to zero for all $u^{\star} $ implies that
\begin{equation}\label{eq:optcond1}
 \gamma u(x, y) +  \lambda(x, y) = 0 \qquad x \in D, \: y \in \Gamma.
\end{equation}
More specifically, considering the structure of $u$ in~\eqref{eq:u}, in
which only the mean~$\bar{u}$ is unknown, the optimality equation reads:
\begin{equation}
  \gamma \bar{u}(x) +
\int_\Gamma \del{\gamma u'+ \lambda}\rho \dif y =0.
\end{equation}
Since the mean of $u'$ is zero, this optimality condition reduces to
\begin{equation}
\label{eq:optcond2}
 \gamma \bar{u}(x) + \int_\Gamma  \lambda \rho \dif y
  = 0 \qquad x \in D.
\end{equation}
The case of a boundary control is handled in the same fashion, but is
omitted for brevity.

The optimality conditions in~\eqref{eq:optcond1} and~\eqref{eq:optcond2}
permit the expression of the control ($u$ or $\Bar{u}$) as a function
of the adjoint function~$\lambda$.  The control can be eliminated
from the first-order optimality equations, leaving a reduced optimality system
in terms of the state and adjoint variables.  In summary, the
optimal control problem involves solving the parametric constraint
problem~\eqref{eq:constr1y}--\eqref{eq:constr3y}, with $u$ or $\Bar{u}$
eliminated using \eqref{eq:optcond1} or \eqref{eq:optcond2}, and the
parametric adjoint problem~\eqref{eq:adj1}--\eqref{eq:adj3}.

The reduced optimality system corresponding to $\mathcal{L}_{2}$
is constructed in the same fashion as for~$\mathcal{L}_{1}$.
It consists of the parametric constraint problem
in~\eqref{eq:constr1y}--\eqref{eq:constr3y}, with $u$ or $\Bar{u}$
eliminated using \eqref{eq:optcond1} or \eqref{eq:optcond2}, and a
parametric adjoint problem that reads:
\begin{align}
\label{eq:adj_alt1}
 - \nabla_{x} \cdot( \kappa(x,y) \nabla_{x}\lambda(x,y)) &= \beta z(x,y)
    -\alpha \hat{z}(x,y)
        + (\alpha- \beta) \int_\Gamma z \rho \dif y  &x \in D,\: y\in \Gamma,
\\
\label{eq:adj_alt2}
\kappa(x,y) \nabla_{x} \lambda(x,y) \cdot n
  &= 0 &x \in \partial D_N,\: y\in \Gamma,
\\
\label{eq:adj_alt3}
\lambda(x,y)  &= 0  &x \in \partial D_D,\: y\in \Gamma.
\end{align}

\section{Stochastic finite element solution}
\label{sec:sfem}

To compute approximate solutions to the optimality system derived in
Section~\ref{sssec:optsys}, both a stochastic Galerkin and a stochastic
collocation finite element discretisation are formulated.  For brevity,
we present the formulation for distributive control only. The boundary
control case, which is considered in the examples section, is formulated
analogously.

For both the stochastic Galerkin and stochastic collocation methods,
for the spatial domain we define a space $V_{h} \subset H^{1}_{0}(D)$
of standard Lagrange finite element functions on a triangulation
$\mathcal{T}$ of the domain~$D$:
\begin{equation}
  V_{h} := \cbr{v_{h} \in H^{1}_{0}(D): v_{h} \in P_{k}(K)
            \ \forall \ K \in \mathcal{T}},
\end{equation}
where $K \in \mathcal{T}$ is a cell and $P_{k}$ is the space of Lagrange
polynomials of degree~$k$. The space $V_{h}$ is spanned by the basis
functions $\cbr{\phi_{i}}_{i=1}^{N}$.

\subsection{Stochastic Galerkin finite element method}
\label{ssec:gal}

To develop a stochastic Galerkin finite element formulation, we consider
a finite dimensional space $Y_{p} \subset L_{\rho}^{2}\del{\Gamma}$
for the stochastic dimension.  For $Y_{p}$, we adopt a polynomial basis,
also known as a generalised polynomial chaos~\citep{Xiu2009,Ghanem2003}.
We first consider multivariate polynomials $\psi_{q}: \Gamma \rightarrow
\mathbb{R}$ generated via~\citep{Matthies2005}
\begin{equation}
\label{eq:psi}
 \psi_{q} = \prod_{i=1}^{L} \varphi_{q_{i}}(y_i),
\end{equation}
where $q = \del{q_{1}, \hdots, q_{L}} \in \mathbb{N}^{L}$ is a multi-index
that satisfies $|q| \le p$ and $\varphi_{q_{i}} : \Gamma_{q_{i}}
\rightarrow \mathbb{R}$ is a one-dimensional orthogonal polynomial
of degree~$q_{i}$.  Note that $p$ is the prescribed total polynomial
degree and recall that $L$ is the number of random variables in the
problem. The polynomials~$\psi_q$ are orthonormal with respect to
the $L^2_\rho(\Gamma)$-inner product defined in~\eqref{eq:innerL2}.
The space $Y_{p}$ is defined in terms of $\psi_{q}$:
\begin{equation}
  Y_{p} := {\rm span} \cbr{\psi_{q}
               : q \in \mathcal{M}} \subset L_{\rho}^{2}(\Gamma),
\end{equation}
where $\mathcal{M}$ is the set of all multi-indices of length $L$
that satisfy $|q| \le p$ for $q \in \mathcal{M}$.  It holds that
$\dim\del{Y_{p}} = \dim\del{\mathcal{M}} = Q = \binom{L+p}{L} =
(L+p)!/(L!p!)$. Consequently, there exists a bijection
\begin{equation}
 \mu: \cbr{1,\hdots,Q} \rightarrow \mathcal{M}
\end{equation}
that assigns a unique integer $j$ to each multi-index $\mu(j)\in
\mathcal{M}$.

A function $z_{hp} \in V_{h} \otimes Y_{p}$ is represented as
\begin{equation}\label{eq:z_hp}
  z_{hp}(x, y) =
 \sum_{i=1}^{N}\sum_{q \in \mathcal{M}} z_{i,q} \phi_{i}(x)\psi_{q}(y),
\end{equation}
where $z_{i,q}$ is a degree of freedom (recall that $V_{h}$ is spanned
by $\cbr{\phi_{i}}_{i=1}^{N}$).  Using the above, for the case of the
cost functional $\mathcal{J}_1$ with unknown stochastic $u$, a Galerkin
formulation of the optimality conditions reads: find $z_{hp} \in V_{h}
\otimes Y_{p}$ and $\lambda_{hp} \in V_{h}\otimes Y_{p}$ such that
\begin{multline}
\label{eq:weakGs1}
 \int_\Gamma \int_D \kappa \nabla_{x} z_{hp} \cdot \nabla_{x} w_{hp} \rho \dif x \dif y
      + \frac{1}{\gamma}  \int_\Gamma\int_D  \lambda_{hp} w_{hp} \rho \dif x \dif y
\\
      = \int_{\Gamma} \int_{\partial D_N } g w_{hp} \rho \dif s \dif y
  \quad \forall \  w_{hp} \in V_{h} \otimes Y_{p}
\end{multline}
and
\begin{multline}
\label{eq:weakGs2}
\int_\Gamma \int_D \kappa \nabla_{x} \lambda_{hp} \cdot \nabla_{x} r_{hp} \rho \dif x \dif y
      - \del{\alpha + \beta}\int_\Gamma \int_D  z_{hp} r_{hp}  \rho \dif x \dif y
\\
      + \beta \int_D  \del{\int_\Gamma z_{hp} \rho \dif y} \del{\int_{\Gamma} r_{hp} \rho \dif y} \dif x
 = - \alpha\int_\Gamma \int_D \hat{z} r_{hp} \rho \dif x \dif y
  \quad \forall \ r_{hp} \in  V_{h} \otimes Y_{p}.
\end{multline}
For the case that $u$ is decomposed additively according to \eqref{eq:u}
and only $\Bar{u}$ is to be computed, equation~\eqref{eq:weakGs1} is
replaced by:
\begin{multline}\label{eq:weakGs1b}
  \int_\Gamma \int_D \kappa \nabla_{x} z_{hp} \cdot \nabla_{x} w_{hp} \rho \dif x \dif y
      + \frac{1}{\gamma} \int_D  \del{\int_\Gamma \lambda_{hp} \rho \dif y}
      \del{\int_\Gamma w_{hp} \rho \dif y} \dif x
\\
    = \int_\Gamma \int_D u' w_{hp} \rho \dif x \dif y
        +  \int_\Gamma \int_{\partial D_N } g w_{hp} \rho \dif s \dif y.
\end{multline}

It can be helpful to examine the structure of the matrix systems that
result from the above finite element problems.  The space $V_{h}$ is
spanned by $N$ nodal basis functions, with $N$ the number of spatial
degrees of freedom. This enables the construction of a mass matrix $M
\in \mathbb{R}^{N \times N}$ and a set of stiffness matrices $K_{i} \in
\mathbb{R}^{N \times N}$, $i =1,\hdots, S$, with each corresponding to
a deterministic diffusion coefficient $\kappa_{i}$ in~\eqref{eq:kappa}.
The space $Y_{p}$ is spanned by $Q$ polynomials and the stochastic
discretisation of $\kappa$ in~\eqref{eq:kappa} defines a set of matrices
$C_i \in \mathbb{R}^{Q\times Q}$, $i = 1,\hdots, S$, whose elements equal
\begin{equation}
 C_i(j,k) := \int_\Gamma \zeta_i \psi_{\mu(j)} \psi_{\mu(k)} \rho \dif  y
        \quad \mathrm{for}\ j, k = 1, \hdots, Q.
\end{equation}
The resulting matrix formulation of the finite element problem
in~\eqref{eq:weakGs1} and~\eqref{eq:weakGs2} is then given by:
\begin{equation}
\label{eq:mat_kron}
  \del{\sum_{i=1}^S
  \begin{bmatrix}
    C_i & \bm{0}_Q
    \\ \bm{0}_Q & C_i
  \end{bmatrix}
\otimes K_i
    +
  \begin{bmatrix}
      \bm{0}_Q& T\del{\frac{1}{\gamma}, \frac{-\epsilon}{\gamma}} \\
    T(-\alpha, -\beta) &\bm{0}_Q
  \end{bmatrix}
\otimes M}
\begin{bmatrix}
  \bm{z}_1 \\ \hdots \\ \bm{z}_Q \\ \bm{\lambda}_1 \\ \hdots\\ \bm{\lambda}_Q
\end{bmatrix}
=
\begin{bmatrix}
\bm{g}_1+\epsilon \bm{u}^{\prime}_{1}
\\
\hdots
\\\bm{g}_Q+\epsilon \bm{u}'_Q\\ -\alpha \bm{\hat{z}}_1\\\hdots\\
-\alpha \bm{\hat{z}}_Q
\end{bmatrix},
\end{equation}
where $\epsilon = 0 $ in the case of a stochastic control $u$
(see~\eqref{eq:weakGs1}) and $\epsilon = 1$ when only $\bar{u}$ is unknown
(see~\eqref{eq:weakGs1b}). The matrix $\bm{0}_Q \in \mathbb{R}^{Q \times
Q}$ is a zero matrix. The diagonal matrix $T \in \mathbb{R}^{Q \times Q}$
is defined as:
\begin{equation}
\label{eq:T}
 T(a,b):= \mathrm{diag}\del{
  \begin{bmatrix}
    a & a+b  & \hdots & a+b
  \end{bmatrix}^T}.
\end{equation}
The vectors~$\bm{z}_q,\,\bm{\lambda}_q \in \mathbb{R}^N$, with $q =
1, \hdots, Q$, collect the spatial degrees of freedom of $z_{hp}$
and $\lambda_{hp}$, respectively, in~\eqref{eq:z_hp}; the vectors
$\bm{\hat{z}}_q, \bm{g}_q,\,\bm{u}'_q\in\mathbb{R}^N$ represent the finite
element discretisation of $ \int_\Gamma\int_D \hat{z} \psi_q\rho \dif
x \dif y$, $\int_\Gamma\int_{\partial D_N} g \psi_q \rho \dif s \dif y$
and $\int_\Gamma \int_D u' \psi_q\rho \dif x \dif y$, respectively.

\subsection{Stochastic collocation finite element method}
\label{ssec:coll}

In applying the stochastic collocation method~\citep{Babuska2007,Xiu2005},
we solve the optimality system at a collection of collocation points
$\cbr{\Hat{y}^i}_{i=1}^{Q}$, where $\Hat{y}^i \in \Gamma$. Typically,
the collocation points are determined by constructing a sparse grid,
see \citep{Babuska2010,Back2009} for details on the point selection.
Of relevance at this stage is that integrals of the form $\int_\Gamma
\del{\cdot} \rho \dif y$ can be approximated via
\begin{equation}\label{eq:cuba}
 \int_\Gamma z(x,y) \rho dy \approx \sum_{i=1}^Q z(x, \Hat{y}^{i}) w_{i},
\end{equation}
where $w_i$ is an appropriate cubature weight and $Q$ is the number of
cubature points.

From the $Q$ realisations $\cbr{z_i = z(x, \Hat{y}^i)}_{i=1}^{Q}$,
a semi-discrete global approximation of the response can be constructed,
\begin{equation}\label{eq:coll}
 z_p(x,y) = \sum_{i=1}^{Q} z(x, \hat{y}^i) \widetilde{\psi}_i(y),
\end{equation}
where the multivariate polynomials~$\widetilde{\psi}_i$ are commonly
interpolatory Lagrange polynomials, as defined in \citet{Babuska2010}. The
polynomial representation permits an exact evaluation of the expectation
of $z_p$ by the cubature rule~\eqref{eq:cuba} when the order of the
cubature rule is sufficiently high~\citep{Xiu2005}. Other moments of
the response, as well its probability density function, can be easily
extracted using a cubature rule and the expansion in~\eqref{eq:coll}.

The main computational cost of the stochastic collocation method is
associated with solving the optimality system at each collocation
point. That is, for each $\Hat{y}^{i}$, $i= 1,\hdots,Q$, the state
equation~\eqref{eq:constr1y}--\eqref{eq:constr3y},
\begin{align}
\label{eq:coll1}
 -\nabla_{x} \cdot \del{\kappa(x,\hat{y}^i) \nabla_{x} z(x,\hat{y}^i)}
&= u(x,\hat{y}^i)  &\mathrm{on} \  D,
\\ \label{eq:coll2}
 z(x,\hat{y}^i) &= z_D(x,\hat{y}^i)  &\mathrm{on} \  \partial D_D,
\\ \label{eq:coll3}
 \kappa(x,\hat{y}^i) \nabla_{x} z(x,\hat{y}^i)\cdot n
      &= g(x,\hat{y}^i) &\mathrm{on} \  \partial D_N,
\end{align}
and the adjoint problem~\eqref{eq:adj1}--\eqref{eq:adj4},
\begin{align}
 \label{eq:coll4}
 - \nabla_{x} \cdot \del{\kappa(x,\hat{y}^i) \nabla_{x}\lambda(x,\hat{y}^i)}
      &= (\alpha+\beta)z(x,\hat{y}^i)-\alpha \hat{z}(x,\hat{y}^i)
- \beta \sum_{j=1}^{Q} z(x,\hat{y}^j) w_j  & \mathrm{on} \  D,
\\ \label{eq:coll5}
\kappa(x,\hat{y}^i) \nabla_{x} \lambda(x,\hat{y}^i)\cdot n
&= 0   &\mathrm{on} \  \partial D_N,
\\ \label{eq:coll6}
\lambda(x,\hat{y}^i)  &= 0  &\mathrm{on} \  \partial D_D,
\end{align}
are solved.  For the case that the unknown control $u$ is wholly
stochastic, equation~\eqref{eq:optcond1} is used to eliminate $u$ in
favour of~$\lambda$:
\begin{equation}
\label{eq:collop1}
 u(x,\hat{y}^i) = -\frac{1}{\gamma} \lambda(x, \hat{y}^i).
\end{equation}
If the control function has the additive structure of~\eqref{eq:u}, then
\begin{equation}
\label{eq:collop3}
  \bar{u}(x) = -\frac{1}{\gamma} \sum_{j=1}^{Q} \lambda(x, \hat{y}^j) w_{j},
\end{equation}
is used.  This equation corresponds to the case where only the mean
part $\Bar{u}$ in~\eqref{eq:u} is to be computed and the stochastic
part $u'$ is prescribed.

\begin{remark}\label{remark:coupcoll}
When moments of unknown functions, e.g., the expectation or standard
deviation, appear in the control problem, the $Q$ sets of equations in
the stochastic collocation method become coupled.  For the case that $u =
\Bar{u}+ u'$, the state equation~\eqref{eq:coll1} is of the form
\begin{equation}
\label{eq:coll1_det}
 -\nabla_{x} \cdot \del{\kappa(x,\hat{y}^i) \nabla_{x} z(x,\hat{y}^i)}
    = \Bar{u}(x) + u^{\prime}(x, \hat{y}^i),
\end{equation}
where $u^{\prime}(x, \hat{y}^i)$ is known.  Following the usual process
of applying equation~\eqref{eq:collop3} to eliminate the control function
$\Bar{u}$ in the constraint equation~\eqref{eq:coll1_det} in favour of
the adjoint variable $\lambda$, the $Q$ deterministic problems become
coupled via this term. A similar coupling is observed for $\beta\ne
0$ in \eqref{eq:coll4} and $\alpha \ne 0 $ in \eqref{eq:adj_alt1}.
The advantage of decoupled systems of equations that is usually associated
with the stochastic collocation method is therefore lost.
\end{remark}

Once the collocation systems are solved, the cost functional
$\mathcal{J}_{1}$ can be evaluated as follows:
\begin{multline}
 \mathcal{J}_{1}(z,u,g)
    =  \sum_{i=1}^{Q} w_i \int_D z_i \del{\frac{\alpha+\beta}{2} z_{i}
    - \alpha  \hat{z} - \frac{\beta}{2} \sum_{j=1}^{Q} z_j} \dif x
    + \frac{\alpha}{2} \int_D \hat{z}^2 \dif x
\\
    + \frac{\gamma}{2} \sum_{i=1}^{Q} w_i \int_D u_i^2 \dif x
    + \frac{\delta}{2} \sum_{i=1}^{Q} w_i \int_{\partial D_N} g_{i}^{2} \dif s,
\end{multline}
with $u_i= u(x,\hat{y}^i)$ and $g_i = g(x,\hat{y}^i)$. The scalars $w_i$,
$i = 1, \hdots, Q$, represent the cubature weights corresponding to
the cubature points $\cbr{\Hat{y}^i}_{i=1}^{Q}$, see \eqref{eq:cuba}.
The cost functional $\mathcal{J}_{2}$ can be likewise evaluated.

The coupling of the collocation systems can be visualised through a
matrix formulation. Following from the finite element space~$V_{h}$,
we define a mass matrix~$M\in\mathbb{R}^{N\times N}$ and a set of
stiffness matrices $K_i \in \mathbb{R}^{N\times N}$, $i=1,\hdots,Q$,
each corresponding to a~$\kappa(x,\hat{y}^i)$. For a distributive
control function~$u$ ($\delta = 0$ in~\eqref{eq:costy}), the matrix
formulation of the stochastic collocation finite element systems
in~\eqref{eq:coll1}--\eqref{eq:collop3} reads:
\begin{equation}
 \begin{bmatrix}
 K_1 & & & \frac{1}{\gamma} M &
\frac{\epsilon}{\gamma} M \: \hdots& \frac{\epsilon}{\gamma} M \\
& \ddots & & \frac{\epsilon}{\gamma} M& \ddots & \vdots\\
& & K_{Q} &\vdots &\hdots \:\frac{\epsilon}{\gamma} M  & \frac{1}{\gamma} M \\
 &  & & K_1 & & \\
\multicolumn{3}{c}{\del{ -(\alpha+\beta)I_Q + \beta \del{\mathbf{1} \mathbf{w}^T} }\otimes M
 }&  & \ddots & \\
 &   & & & & K_{Q}
\end{bmatrix}
\begin{bmatrix}
  \bm{z}_1 \\ \vdots \\ \bm{z}_{Q} \\ \bm{\lambda}_1 \\ \vdots \\ \bm{\lambda}_{Q}
\end{bmatrix}
=
\begin{bmatrix}
        \bm{g}_1+\epsilon\bm{u}'_1 \\ \hdots \\ \bm{g}_{Q}+\epsilon\bm{u}^{\prime}_{Q} \\
        -\alpha \bm{\hat{z}}_1 \\ \hdots \\ -\alpha \bm{\hat{z}}_Q
\end{bmatrix},
\label{eq:collmat}
\end{equation}
where $\epsilon = 0 $ for the case that $u$ is wholly
stochastic (see~\eqref{eq:collop1}) and $\epsilon = 1$ when
only~$\Bar{u}$ is unknown (see~\eqref{eq:collop3}).  The vectors
$\bm{z}_i,\,\bm{\lambda_i},\,\bm{\hat{z}}_i,\, \bm{g}_i,\, \bm{u}'_i
\in \mathbb{R}^{N}$, with $i =1, \hdots, Q$, are the finite element
representations of $z(x, \hat{y}^i)$, $\lambda(x, \hat{y}^i)$, $\hat{z}(x,
\hat{y}^i)$, $g(x, \hat{y}^i) $ and $u'(x, \hat{y}^i) $, respectively.
The matrix $I_Q \in \mathbb{R}^{Q\times Q}$ is an identity matrix and
the vectors $\mathbf{1}, \mathbf{w} \in \mathbb{R}^Q$ are defined by $[1\
\hdots \ 1]^T$ and $[w_1\ \hdots \ w_Q]^T$, respectively.

In contrast with the stochastic inverse type problems in \citet{Borzi2009}
and \citet{Zabaras2008} to which the non-intrusivity property of the
stochastic collocation method is preserved, deterministic simulation
software cannot readily be reused for $\beta \ne 0$ or $\epsilon \ne 0$.

\begin{remark}
\label{remark:cost_coll}
The stochastic collocation method generally requires more stochastic
degrees of freedom, i.e., a larger~$Q$, than the stochastic
Galerkin method in order to solve a stochastic PDE to the same
accuracy~\citep{Back2009}. Therefore, for the same accuracy, due to the
coupling of the stochastic collocation systems (see \eqref{eq:collmat})
the stochastic collocation method will likely involve a greater
computational cost compared to the stochastic Galerkin method.
When presenting numerical examples, we will therefore only apply the
stochastic collocation method to problems for which the non-intrusivity of
the stochastic collocation method is maintained, i.e., only for the cost
functional $\mathcal{J}_{1}$ with unknown stochastic $u$ and~$\beta = 0$.
\end{remark}
\section{Additional regularisation of stochastic optimal control problems}
\label{sec:regul}

The cost functionals $\mathcal{J}_{1}$ and $\mathcal{J}_{2}$ both contain
a regularisation term of the form $(\gamma/2)\norm{u}^2_{L^2(D)\otimes
L^2(\Omega)}$. This $L^{2}$-regularisation is important for the solvability
of the problem, and as $\gamma$ is increased excessive control values
are penalised~\citep{Hinze2009}.  However, in some applications
$L^2$-penalisation may not be sufficient.  We detail in this section how an
additional $H^{1}$-like regularisation can be included into the cost
functionals and show what computational issues then follow.

Assuming sufficient regularity of the relevant functions, we extend
the functional $\mathcal{J}_{1}$ in \eqref{eq:cost} by including a
$H^1$-penalty on the control variables:
\begin{multline}
\label{eq:cost2}
\mathcal{J}_{1,H^1}\del{z, u, g}
    := \frac{\alpha}{2} \norm{z - \hat{z}}^2_{L^2(D)\otimes L^2_\rho(\Gamma)}
    + \frac{\beta}{2} \norm{\mathrm{std}(z)}^2_{L^2(D)}
\\
  + \frac{\gamma}{2}
      \norm{u}^2_{H^1(D) \otimes H^1_{\rho}(\Gamma)}
  + \frac{\delta}{2}
      \norm{g}^2_{H^1(\partial D_N) \otimes H^1_{\rho}(\Gamma)}.
\end{multline}
A stochastic optimal control problem involving the cost
functional $\mathcal{J}_{1, H^{1}}$ and the PDE constraints
in~\eqref{eq:constr1y}--\eqref{eq:constr3y} can be solved by following the
one-shot strategy described in Section~\ref{ssec:one_shot}.  The state
equations are given by \eqref{eq:constr1y}--\eqref{eq:constr3y} and
the adjoint equations are given by \eqref{eq:adj1}--\eqref{eq:adj4}
(see Section~\ref{sssec:optsys}).  For the distributive control~$u$,
the weak optimality condition for the $\mathcal{J}_{1,H^1}$ case reads:
\begin{equation}
\label{eq:optcondH1}
  \int_{\Gamma} \int_D \lambda u^{\star} \rho \dif x \dif y
         + \gamma \del{u, u^{\star}}_{H^{1}(D) \otimes H^{1}_{\rho}(\Gamma)}  = 0
\end{equation}
for all $u^{\star} \in H^1(D) \otimes H^1_{\rho}(\Gamma)$.
For the additive structure of $u$ in \eqref{eq:u}, in which only the
mean~$\bar{u}$ is unknown, the weak optimality condition reduces to
\begin{equation}
\label{eq:optcondH1det}
  \gamma \int_{D} \bar{u} u^{\star}
              + \nabla_{x} \bar{u} \cdot \nabla_{x} u^{\star} \dif x
  + \int_{\Gamma} \int_{D} \lambda  u^{\star} \rho \dif x \dif y = 0
\end{equation}
for all $u^{\star} \in H^{1}(D)$, in which the zero-mean property
of $u'$ has been used.  In contrast with the optimality conditions
in~\eqref{eq:optcond1} and~\eqref{eq:optcond2}, the optimality
conditions~\eqref{eq:optcondH1} and~\eqref{eq:optcondH1det} are partial
differential equations and do not permit the straightforward expression
of the control function as a function of the adjoint field~$\lambda$. As
a result, no reduced optimality system is constructed and the optimality
equations for the state, control and adjoint variables will be solved
simultaneously.

The Galerkin formulation of the optimality system is composed of three
equations, which are the Galerkin formulation of the constraint
equations~\eqref{eq:constr1y}--\eqref{eq:constr3y}, the Galerkin
formulation of the adjoint equations~\eqref{eq:adj1}--\eqref{eq:adj4}, see
\eqref{eq:weakGs2}, and the optimality condition in~\eqref{eq:optcondH1}
or~\eqref{eq:optcondH1det}.  Some structure of the stochastic Galerkin method
can be revealed from its matrix formulation. When the unknown control $u$ is
wholly stochastic, an algebraic system of size $3NQ\times 3NQ$ results,
with $N$ the number of spatial degrees of freedom and $Q$ the number of
stochastic unknowns.  The algebraic system can be shown to possess a
saddle point structure of the form
\begin{equation}\label{eq:saddle}
 \begin{bmatrix}
  A & B^T \\ B  & 0
 \end{bmatrix},
\end{equation}
where the matrix blocks $A \in \mathbb{R}^{2QN \times 2QN} $ and $B \in
\mathbb{R}^{NQ \times 2QN}$ are given by:
\begin{equation}
 A = \begin{bmatrix}
     T(\alpha,\beta)\otimes M & \mathbf{0}_{NQ}
      \\
     \mathbf{0}_{NQ} & \gamma (I_Q + E) \otimes(M+K)
    \end{bmatrix},
\qquad
  B = \begin{bmatrix}
   -\sum_{i=1}^S C_i \otimes K_i &  I_Q \otimes M
    \end{bmatrix}.
\end{equation}
The matrix~$K\in \mathbb{R}^{N \times N}$ is the stiffness matrix for a
deterministic Laplacian operator. The matrix~$T \in \mathbb{R}^{Q\times
Q}$ is defined in \eqref{eq:T}. The elements of the matrix $E\in
\mathbb{R}^{Q\times Q}$ are defined as
\begin{equation}
\label{eq:matE}
 E(j,k) :=
\int_\Gamma \nabla_y \psi_{\mu(j)} \cdot \nabla_y \psi_{\mu(k)} \rho \dif y
\qquad j,k = 1,\hdots Q,
\end{equation}
where $\psi_{\mu(j)}(y)$ denotes a multivariate, orthogonal polynomial,
as defined in~\eqref{eq:psi}. Analytical expressions for computing the
matrix elements in~\eqref{eq:matE} are presented in Appendix~\ref{app}.

The formulation of a stochastic collocation method follows the process
detailed in Section~\ref{ssec:gal}. Applying the collocation method to
the optimality conditions in this section leads to $Q$
deterministic problems which are coupled via the~$Q$ stochastic degrees
of freedom. In the case of an unknown stochastic control function,
the coupling between collocation points is due to the $\nabla_y$-terms
in~\eqref{eq:optcondH1}. When only the mean control is unknown, the
discussion in Remark~\ref{remark:coupcoll} is applicable. To summarise,
the stochastic collocation discretisation of the optimality system is
coupled in the stochastic degrees of freedom when:
\begin{itemize}
  \item $H^{1}_{\rho}(\Gamma)$-control is applied;

  \item the parameter $\beta$ in \eqref{eq:cost2} is non-zero; or

  \item only the deterministic part of the control, $\bar{u}$ or
  $\bar{g}$, is unknown.
\end{itemize}
In these cases, applying a stochastic collocation method is not attractive
since its usual advantages are lost.
\section{Iterative solvers for one-shot systems}
\label{sec:iterative}

The performance and applicability of stochastic finite element methods
relies on fast and robust linear solvers, and this is perhaps even more
so the case for stochastic one-shot methods for optimal control. In the
context of deterministic optimal control problems, various iterative
solvers have been designed, see for example \citet{Schoberl2007},
\citet{Borzi2009a} and \citet{Rees2010}. These solvers can readily
be applied to the deterministic systems resulting from a stochastic
collocation discretisation when there is no coupling between the collocation
points.

When a stochastic Galerkin method is applied, or in the case of coupled
stochastic collocation systems, new iterative solution methods are
required in order to efficiently solve stochastic optimal control
problems. Such solvers typically consist of a Krylov method combined
with a specially tailored preconditioner. This section presents two
approaches for constructing preconditioners for the stochastic Galerkin
one-shot systems.  Following from Remark~\ref{remark:cost_coll},
solvers for coupled stochastic collocation systems are not considered.
The presented preconditioners are applied in Sections~\ref{sec:num1}
and~\ref{sec:num2} when computing numerical examples.

\subsection{Mean-based preconditioner}
\label{ssec:meanbased}

A straightforward and easy-to-implement preconditioner for stochastic
Galerkin systems is the mean-based preconditioner. Applied to
the high-dimensional algebraic system in~\eqref{eq:mat_kron}, the
preconditioner matrix $P_{\mathrm{mean}}$ is defined as
\begin{equation}
\label{eq:meanbased}
P_\mathrm{mean} :=  I_{2Q} \otimes K_1
    +
  \begin{bmatrix}
    0 & \frac{1}{\gamma}  \\
    -\alpha & 0
  \end{bmatrix} \otimes I_Q \otimes M,
\end{equation}
where $I_{2Q}\in \mathbb{R}^{2Q \times 2Q}$ is an identity matrix.
One application of~\eqref{eq:meanbased} requires solving~$Q$ systems,
each of size~$2N \times 2N$, with $N$ the number of spatial degrees of
freedom.  The solution of these smaller systems can be approximated by a
sweep of the collective smoothing multigrid algorithm for deterministic
PDE-constrained optimisation problems~\citep{Borzi2009}.

The mean-based preconditioner performs well for problems with a
low variance of $\kappa$ and a low polynomial order, as analysed
for stochastic elliptic PDEs by \citet{Powell2009}. Its performance
deteriorates however for small penalty parameters $\gamma$ and $\delta$,
$\beta \neq 0$ or $\epsilon \neq 0$ in~\eqref{eq:mat_kron}.

Most stochastic inverse problem examples in Section~\ref{sec:num2} are
solved using the mean-based preconditioner since it typically leads to
the fastest solution method. The number of spatial degrees of freedom
used in the numerical examples is small so that a direct solver for
the $\del{2N \times 2N}$-subsystems is most appropriate.  To overcome
the lack of robustness of the mean-based preconditioner, a collective
smoothing multigrid preconditioner can be applied to the entire
coupled system~\eqref{eq:mat_kron}. A suitable multigrid
solver is summarised below.

\subsection{Collective smoothing multigrid}
\label{ssec:mg}

A robust iterative solver with a convergence rate independent of the
spatial and stochastic discretisation parameters is obtained by applying a
multigrid method directly to the entire system in~\eqref{eq:mat_kron}. The
multigrid preconditioner is so-called `point-based'. That is, it uses
only a hierarchy of spatial grids. The intergrid transfer operators
therefore obey a Kronecker product representation,
\begin{equation}
 I_{2Q}\otimes P_d,
\end{equation}
where $P_d$ is an intergrid transfer operator based on~$K_{1}$.
A simultaneous update of all unknowns per spatial grid point is key to
the multigrid performance, as discussed for stochastic elliptic PDEs in
\citet{Rosseel2010}. As a consequence, this multigrid preconditioner uses
collective smoothing operators, e.g., a block Gauss--Seidel relaxation
method.

Most stochastic control examples in Section~\ref{sec:num1} use a
collective smoothing multigrid preconditioner because of its optimal
convergence properties. In contrast with the mean-based preconditioner,
this multigrid method does not encounter additional difficulties when
solving the coupled Galerkin system~\eqref{eq:mat_kron} in the case that
only $\bar{u}$ is unknown ($\epsilon =1$).

\section{Stochastic control examples}
\label{sec:num1}

A variety of numerical examples are presented to demonstrate the proposed
formulation for control problems. When referring to control problems, we
imply that the control functions are deterministic or have the additive
structure in~\eqref{eq:u}. Also, in the context of control we consider
deterministic targets $\Hat{z}$ only.  Other scenarios are considered
in Section~\ref{sec:num2} in the context of inverse problems.

\begin{figure}
  \centering
  \includegraphics[width=0.32\textwidth]{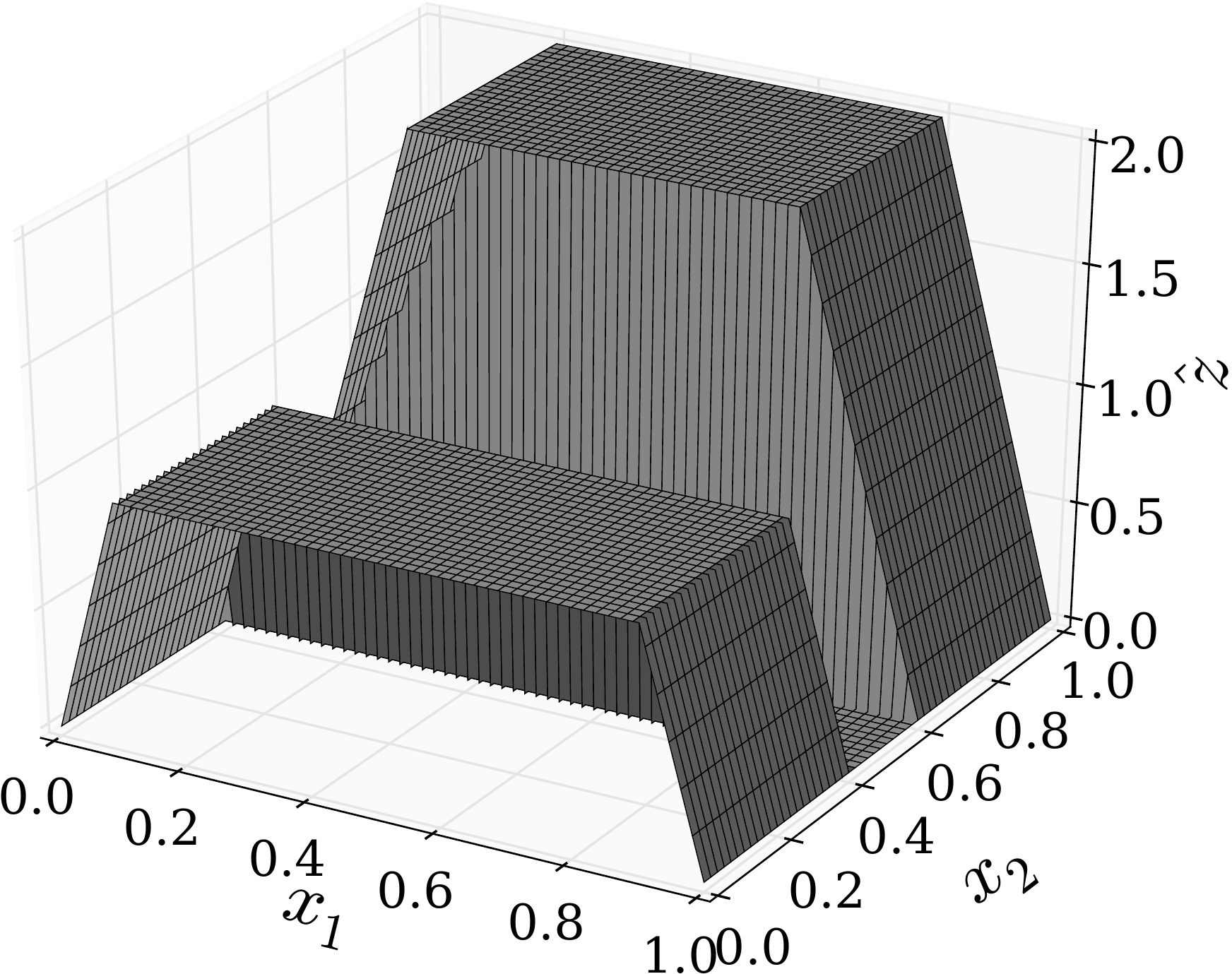}
\caption{Target function $\hat{z}$ for the control examples.}
\label{fig:target}
\end{figure}

In all cases a unit square spatial domain $D = (0, 1)^2$ is considered,
and the boundary of the domain $\partial D$ is partitioned such that
$\partial D_{D} = \{0, 1\} \times (0, 1)$ and $\partial D_{N} = \partial
D \backslash \partial D_{D}$. Unless stated otherwise, zero Dirichlet
and Neumann conditions are applied, i.e., $z_{D} = 0$ on $\partial
D_{D}$ and $g = 0$ on~$\partial D_{N}$.  In all examples, the diffusion
coefficient~$\kappa$ is represented by a Karhunen--Lo\`eve expansion
based on an exponential covariance function with a correlation length of
one and variance $0.25$ (see \eqref{eq:kappa}). The Karhunen--Lo\`eve
expansion is truncated after seven terms; the random variables in this
expansion are assumed to be independent and uniformly distributed on
$\sbr{-\sqrt{3}, \sqrt{3}}$.  A deterministic piecewise-continuous target
function is considered:
\begin{align}
 \hat{z}(x) =
    \begin{cases}
     0            &x \in [0, 1] \times (0.4, 0.6)\\
     1            &x \in (0.1, 0.9) \times [0, 0.4]\\
     2            &x \in (0.2, 0.8) \times [0.6, 1]\\
     10 x_{1}     &x \in [0, 0.1) \times [0,0.4] \cup [0, 0.2) \times [0.6, 1] \\
     10 - 10x_{1} &  \mathrm{otherwise}
    \end{cases}.
\end{align}
This function is illustrated in Figure~\ref{fig:target}.
A spatial finite element mesh of piecewise linear triangular elements
is used. The mesh is constructed as $2^7 \times 2^7$ squares and each
square is subdivided into two triangular finite element cells. This
yields $N=16\,441$ spatial degrees of freedom.  The stochastic Galerkin
discretisation is based on seven-dimensional Legendre polynomials of
order two.  This yields $Q = 36$ (when $u' = 0$ in \eqref{eq:u}) and
the algebraic system in~\eqref{eq:mat_kron} has dimension $1.2 \cdot
10^6 \times 1.2\cdot 10^6$.  The stochastic collocation method employs
a level-two Smolyak sparse grid based on Gauss--Legendre collocation
points with $Q =141$ collocation points.

The computer code used for the examples is built on the library
DOLFIN~\citep{Logg2010}.  The complete computer code for all
presented examples is available as part of the supporting
material~\citep{supporting}.
\subsection{Distributive control with cost functional
\texorpdfstring{$\mathcal{J}_{1}$}{J1}}
\label{ssec:dist_contr_1}

First, we consider a distributive control via $u$ only and using the
cost functional $\mathcal{J}_{1}$ in~\eqref{eq:cost} with $\delta = 0$.
A control function $u$ with the additive structure in~\eqref{eq:u}
is used. The goal is to determine an optimal mean control $\bar{u}$,
which is the deterministic input for the controller response.
\subsubsection{Perfect controller case}
\label{sec:control_deterministic}

\begin{figure}
\centering
\subfloat[mean state]{\includegraphics[width=0.32\textwidth]{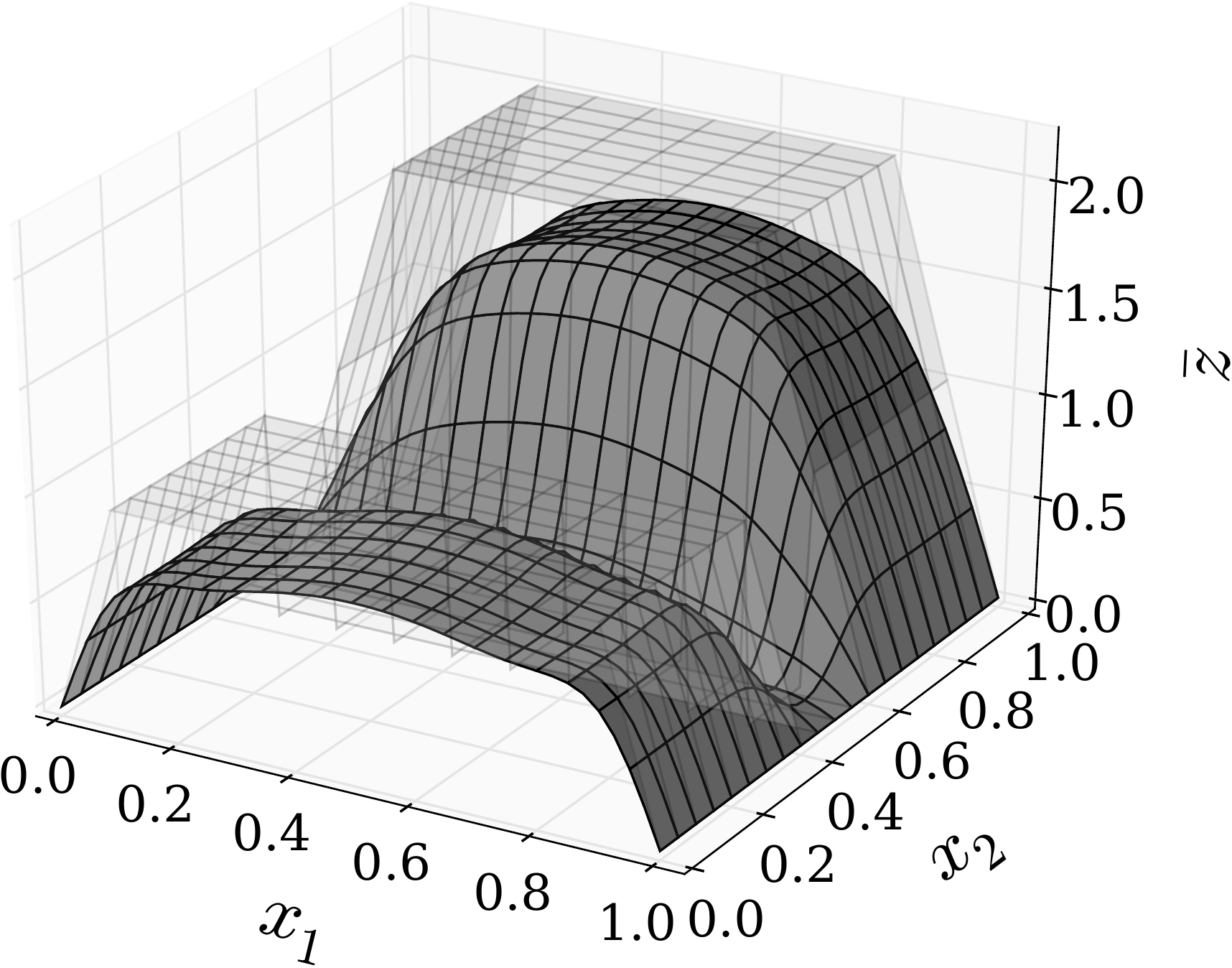}}
\subfloat[variance state]{\includegraphics[width=0.32\textwidth]{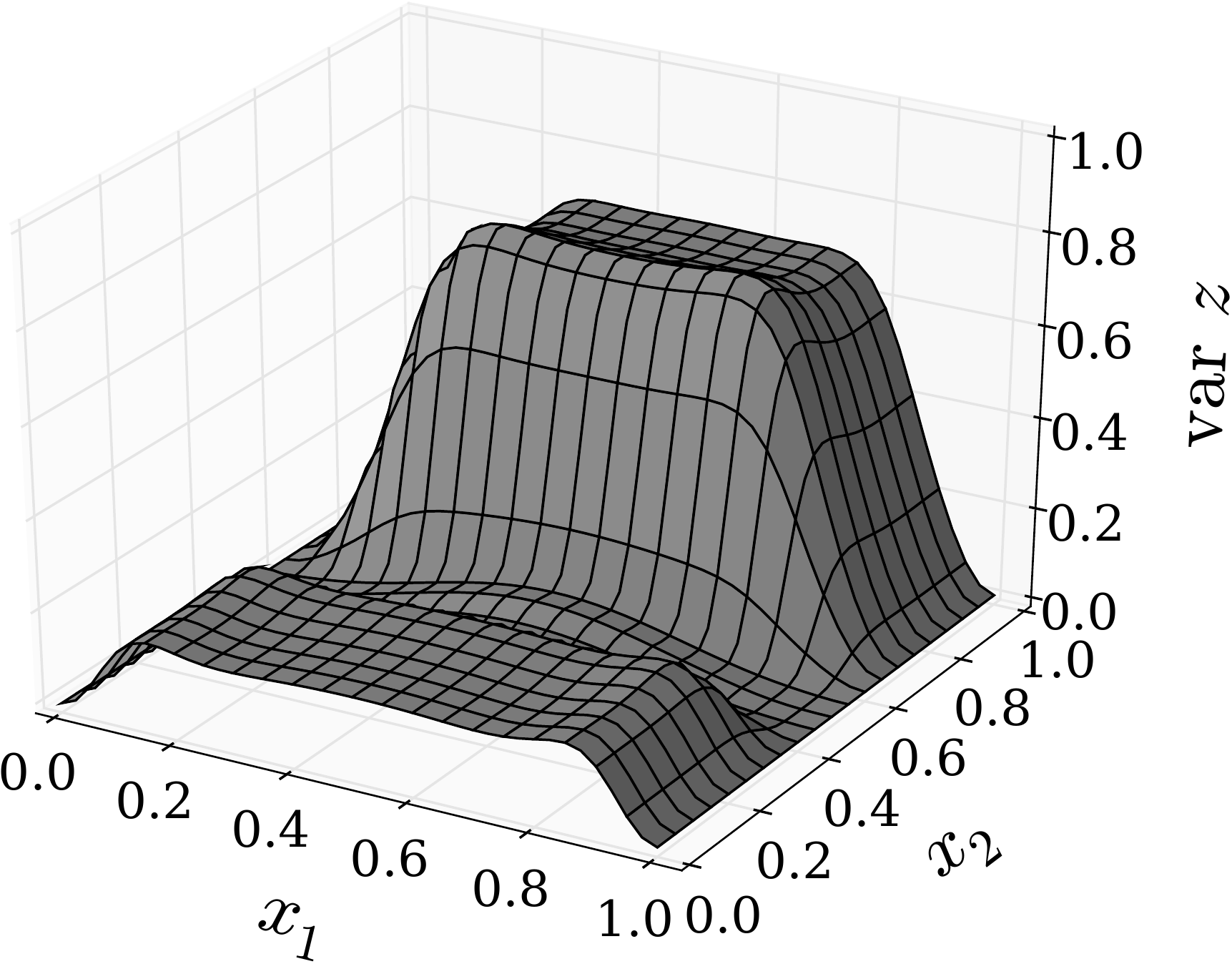}}
\subfloat[deterministic control]{\includegraphics[width=0.32\textwidth]{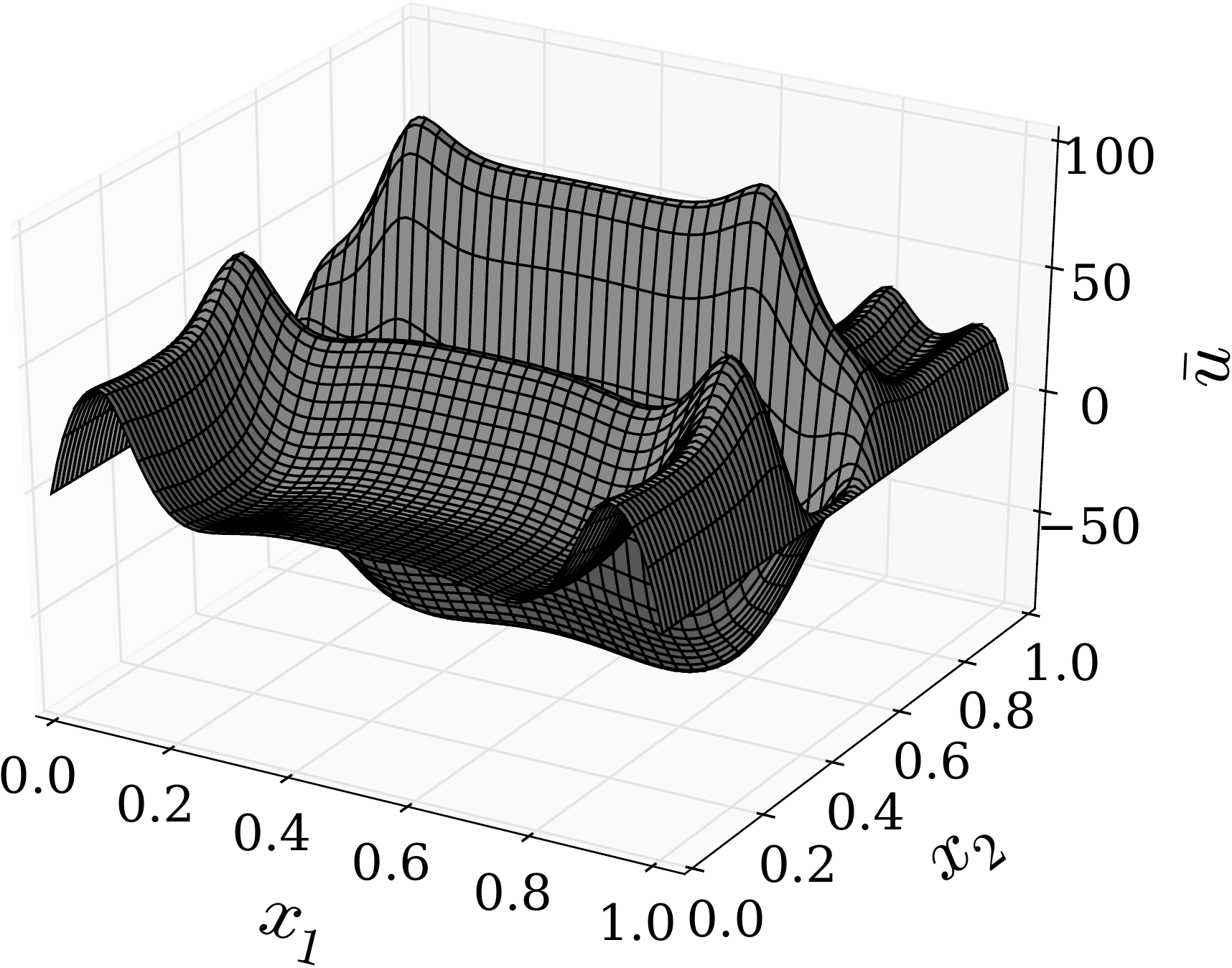}}
\caption{Mean and variance of the optimal state and deterministic
  control variable $u = \Bar{u}$ ($u' = 0$) associated with the cost
  functional $\mathcal{J}_{1}$ and computed with the stochastic Galerkin
  method with $\alpha= 1$, $\beta = \delta = 0$ and $\gamma = 10^{-5}$.
  The target~$\hat{z}$ is illustrated transparently in (a) for reference.}
\label{fig:det_dist_control_A}
\end{figure}

\begin{figure}
\centering
\subfloat[mean state]{\includegraphics[width=0.32\textwidth]{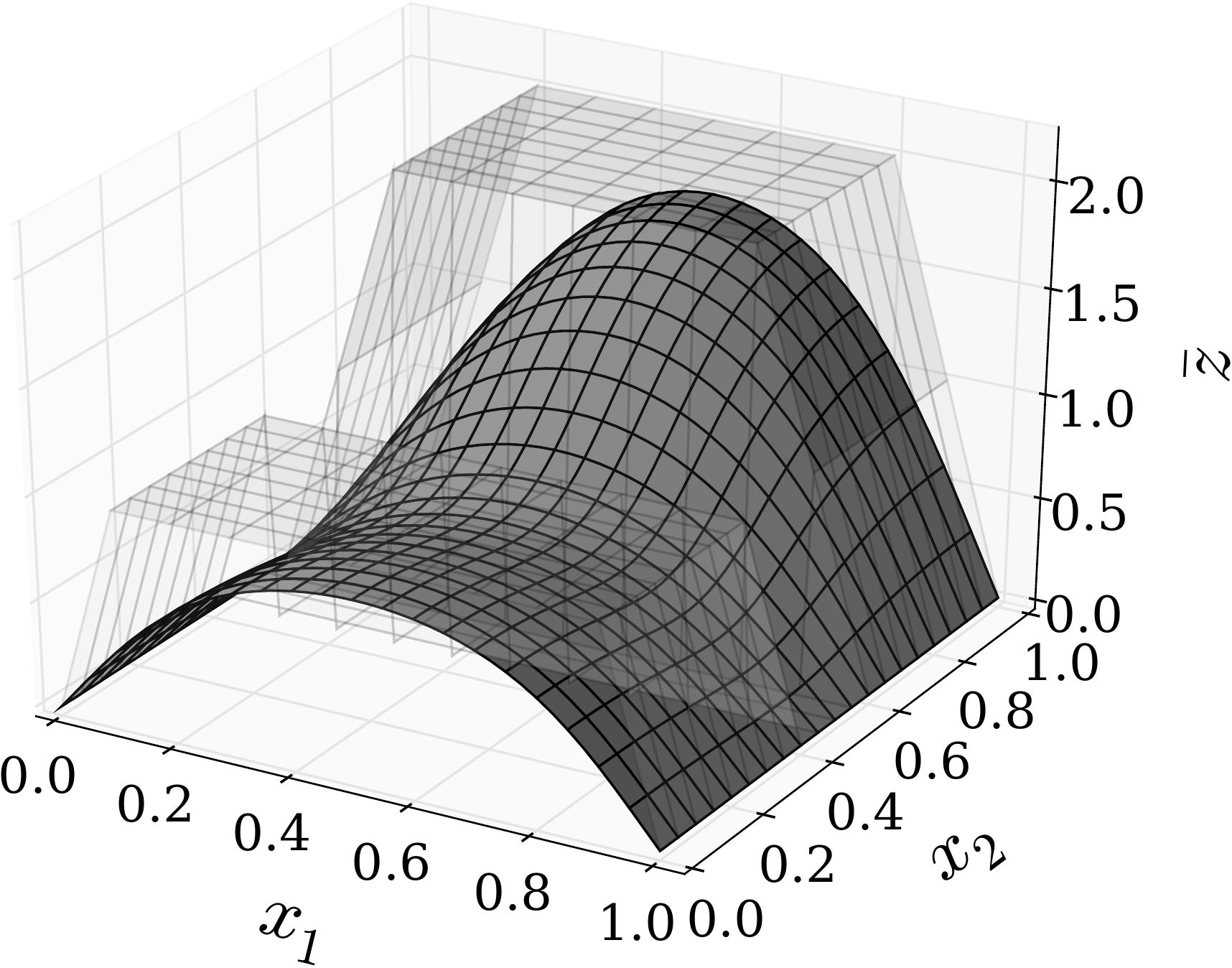}}
\subfloat[variance state]{\includegraphics[width=0.32\textwidth]{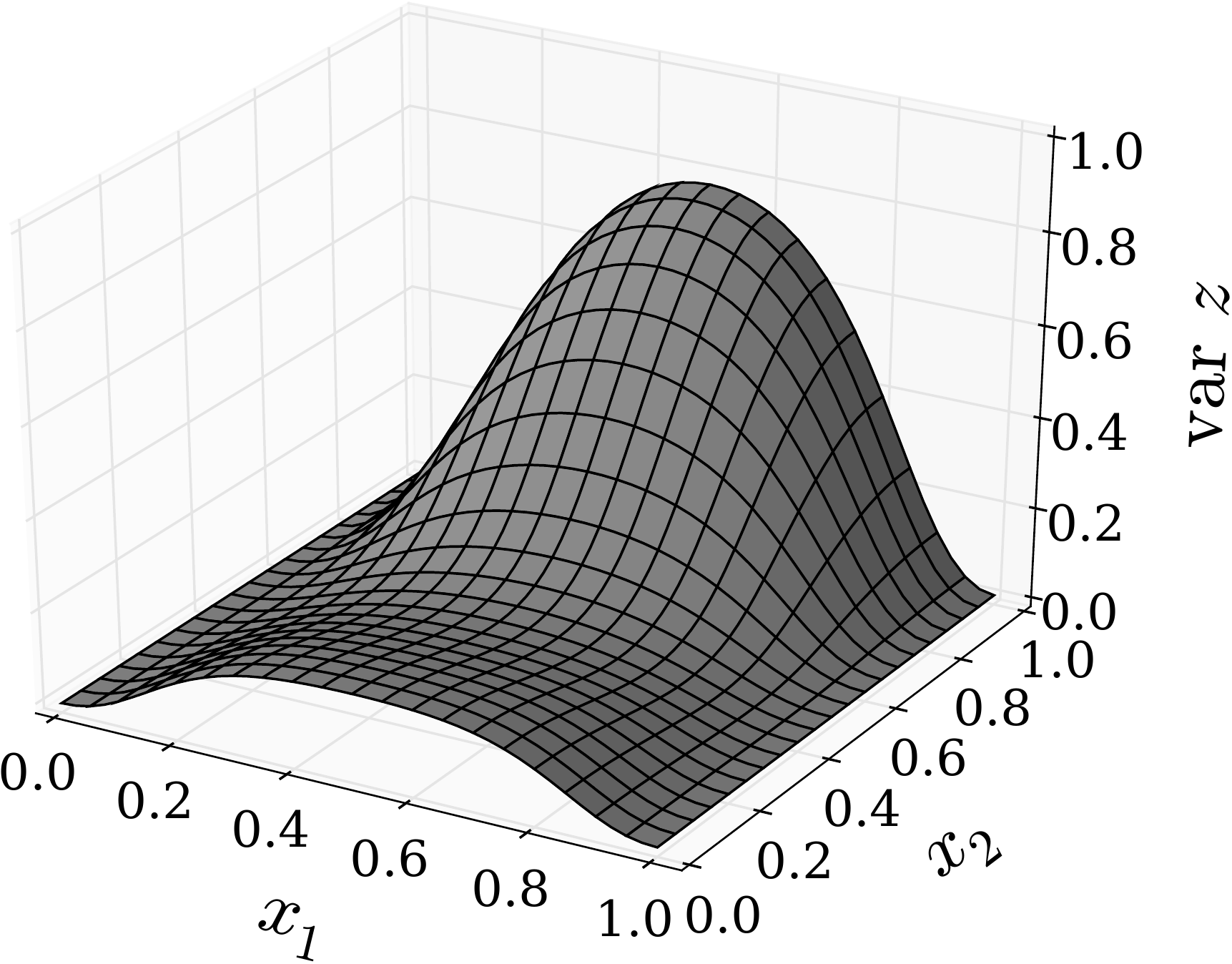}}
\subfloat[deterministic control]{\includegraphics[width=0.32\textwidth]{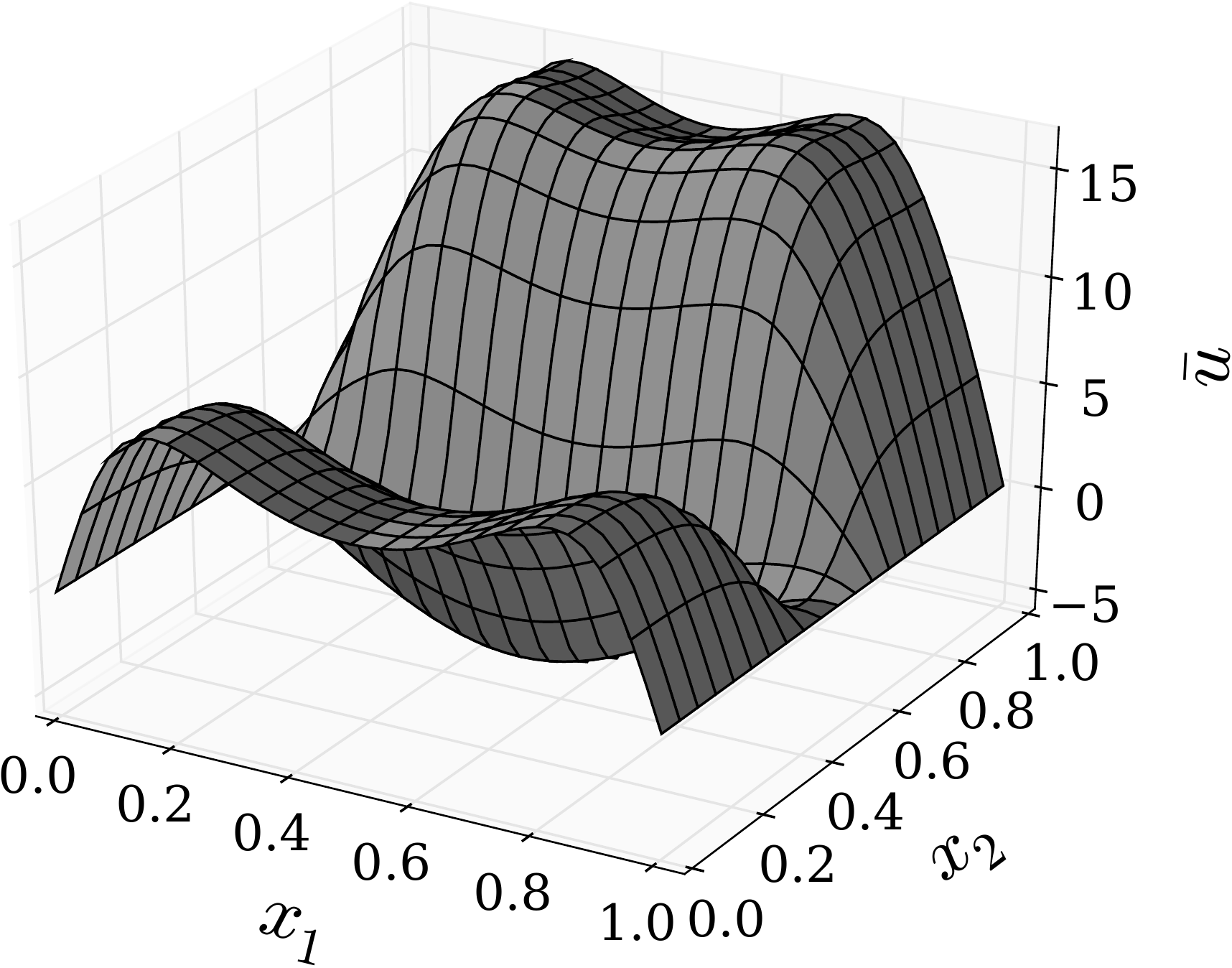}}
\caption{Mean and variance of the optimal state and deterministic
  control variable $u = \Bar{u}$ ($u' = 0$) associated with the cost
  functional $\mathcal{J}_{1}$ and computed with the stochastic Galerkin
  method with $\alpha= 1$, $\beta = \delta = 0$ and $\gamma = 10^{-3}$.
  The target~$\hat{z}$ is illustrated transparently in (a) for reference.}
\label{fig:det_dist_control_B}
\end{figure}

\begin{figure}
  \centering\includegraphics[width=0.5\textwidth]{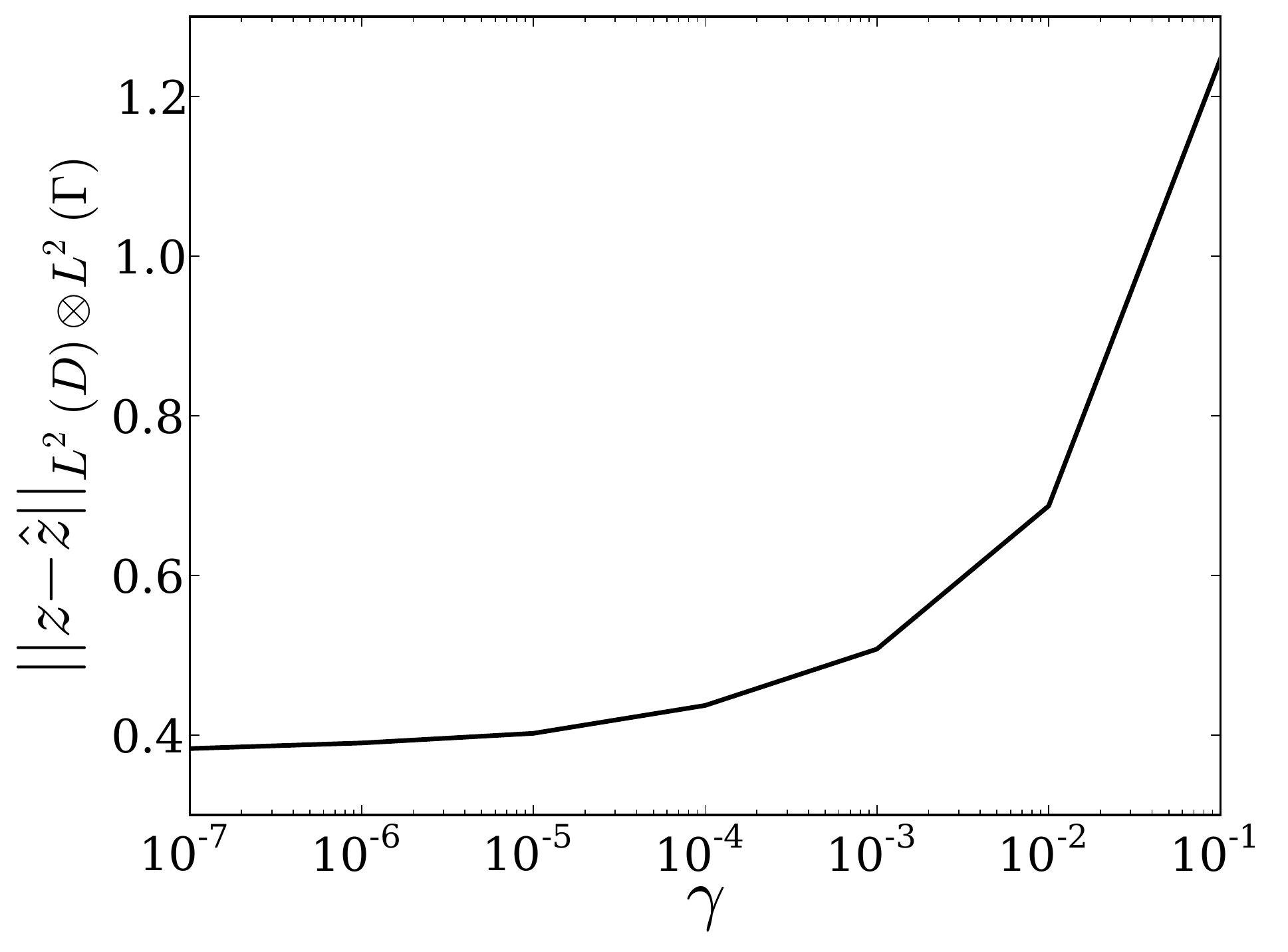}
  \caption{Tracking error $\norm{z - \hat{z}}_{L^2(D)\otimes
    L^2_\rho(\Gamma)}$ as a function of the penalty parameter~$\gamma$ for
    the problem considered in Section~\ref{sec:control_deterministic}.}
\label{fig:gamma}
\end{figure}

The optimal control for the case $u' = 0 $, which corresponds to
the controller action being wholly deterministic, and $\beta = 0$,
which implies no extra control over the variance of~$z$, is computed
and the results are shown in Fig.~\ref{fig:det_dist_control_A}
for the case $\gamma = 10^{-5}$.  The computed mean of
the state function, the variance of the state function and
the control~$\Bar{u}$ are all shown.  The same quantities are
presented in Fig.~\ref{fig:det_dist_control_B} for the case $\gamma =
10^{-3}$. Comparing the results in Figs.~\ref{fig:det_dist_control_A}
and~\ref{fig:det_dist_control_B}, as expected the larger $\gamma$-penalty
term leads to a poorer approximation of the target.  The values of the cost
functional, tracking error $\norm{z - \Hat{z}}^2_{L^2(D)\otimes
L^2_\rho(\Gamma)}$ and standard deviation $\norm{\mathrm{std}(z)}_{L^2(D)}^2$
are summarised in Table~\ref{tab:cost}.  The tracking error as a function of
$\gamma$ is presented in Fig.~\ref{fig:gamma}, illustrating the deterioration
in the quality of the computed result for increasing values of~$\gamma$.

\subsubsection{Imperfect controller case}
\label{sec:dist_contr_imperfect_1}

We now consider the impact of an imperfect controller by introducing
a known non-zero stochastic term $u'$. The stochastic function $u'$ is
modelled by a three-term Karhunen--Lo\`eve expansion based on a zero-mean
Gaussian field with an exponential covariance function, unit variance
and a unit correlation length. For the examples in this section we set
$\gamma = 10^{-3}$.

\begin{figure}
\centering
\subfloat[mean state]{\includegraphics[width=0.32\textwidth]{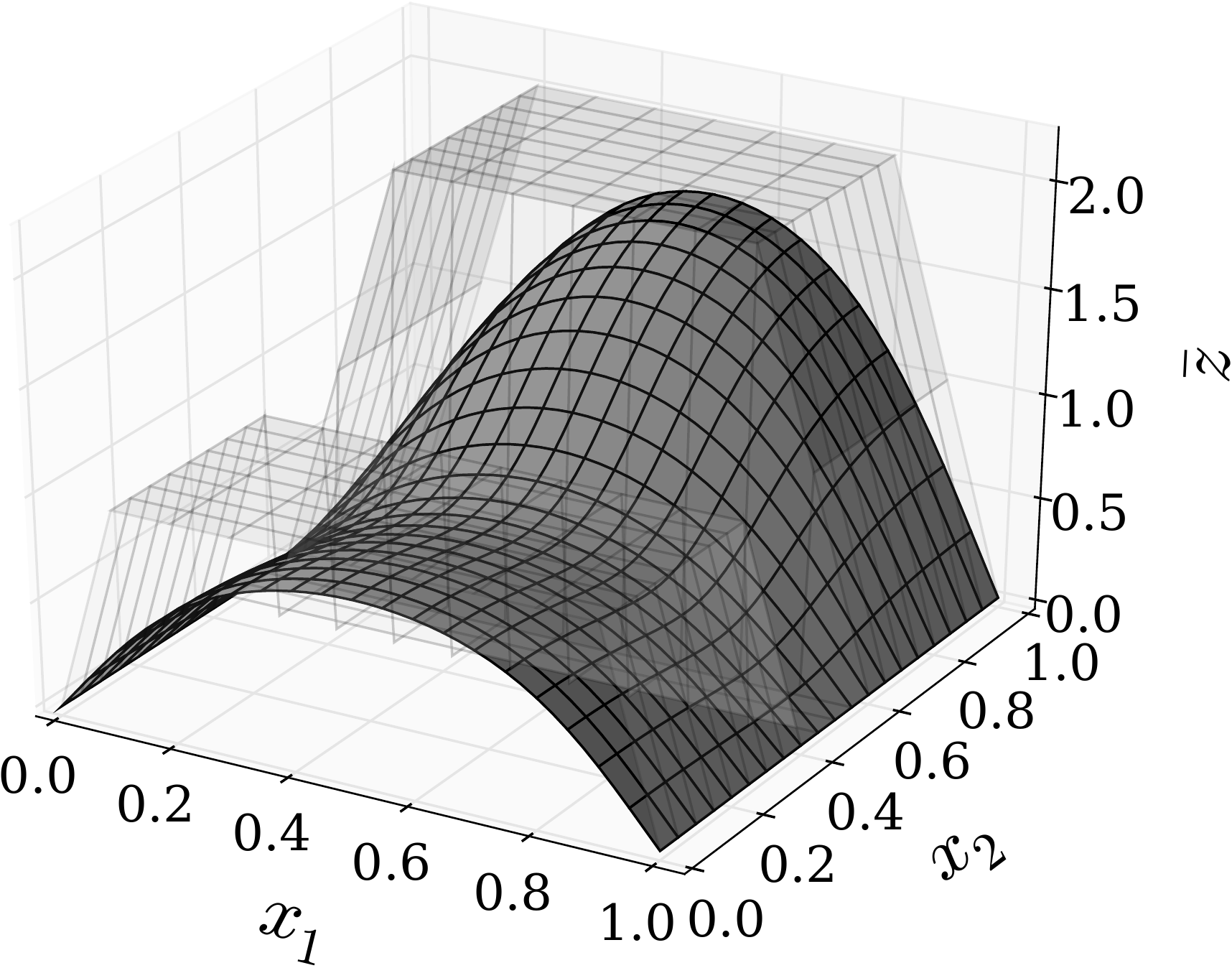}}
\subfloat[variance state]{\includegraphics[width=0.32\textwidth]{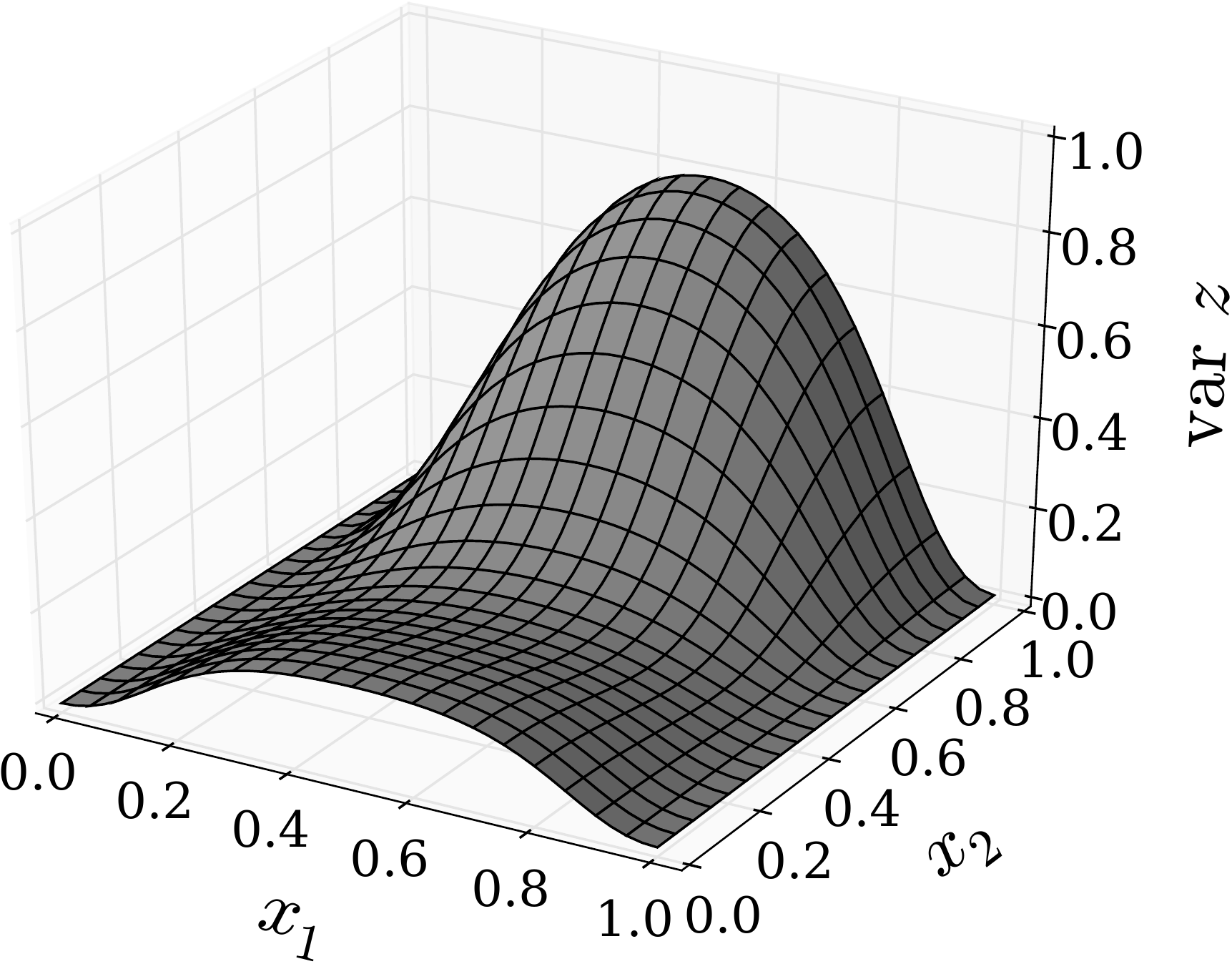}}
\subfloat[mean control]{\includegraphics[width=0.32\textwidth]{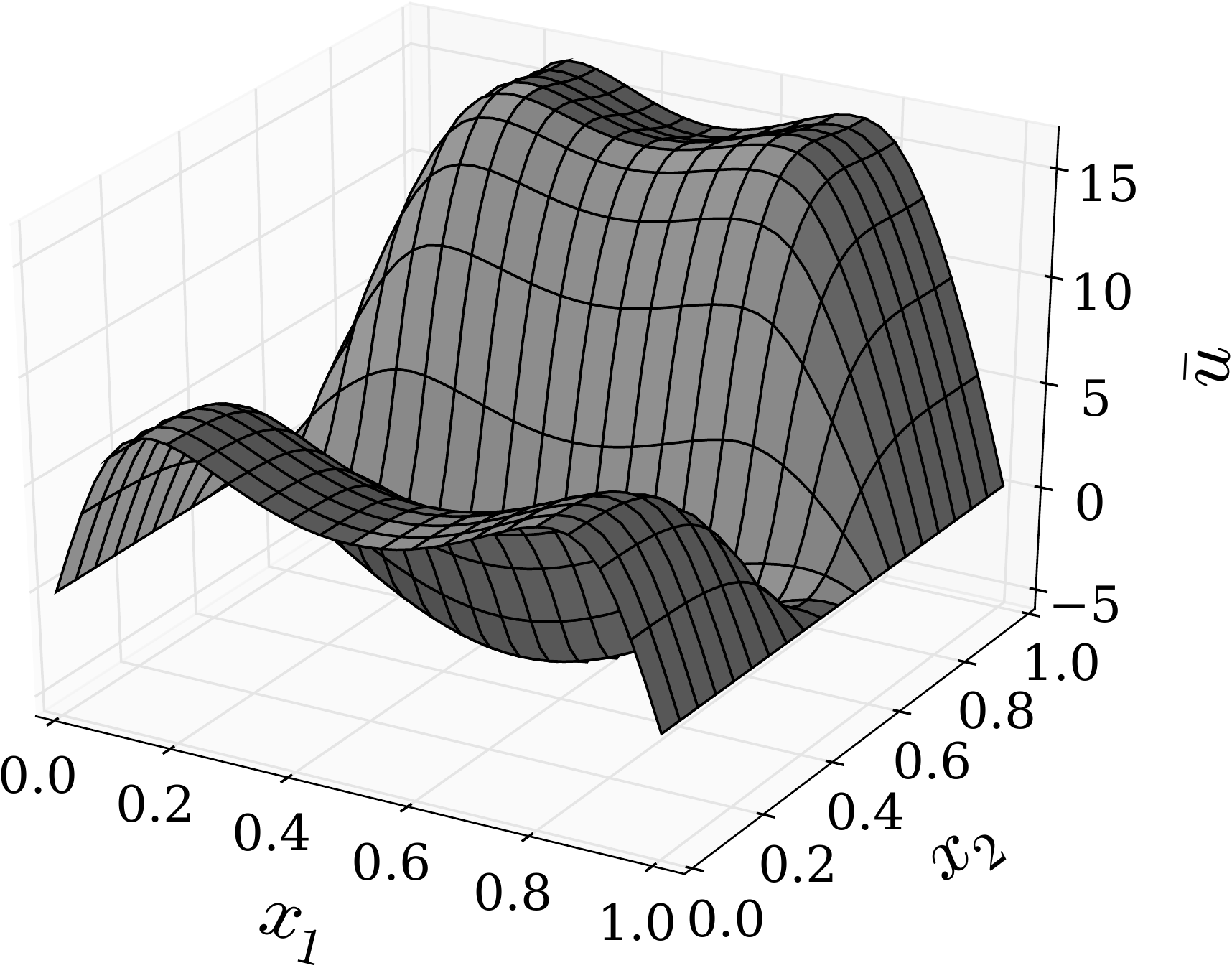}}
\caption{Mean and variance of the optimal state and mean
  control variable $\Bar{u}$ ($u' \ne 0$) associated with the cost
  functional $\mathcal{J}_{1}$ and computed with the stochastic Galerkin
  method with $\alpha= 1$, $\beta = \delta = 0$ and $\gamma = 10^{-3}$.
  The target~$\hat{z}$ is illustrated transparently in (a) for reference.}
\label{fig:det_dist_control_us1_A}
\end{figure}

\begin{figure}
\subfloat[mean state]{\includegraphics[width=0.32\textwidth]{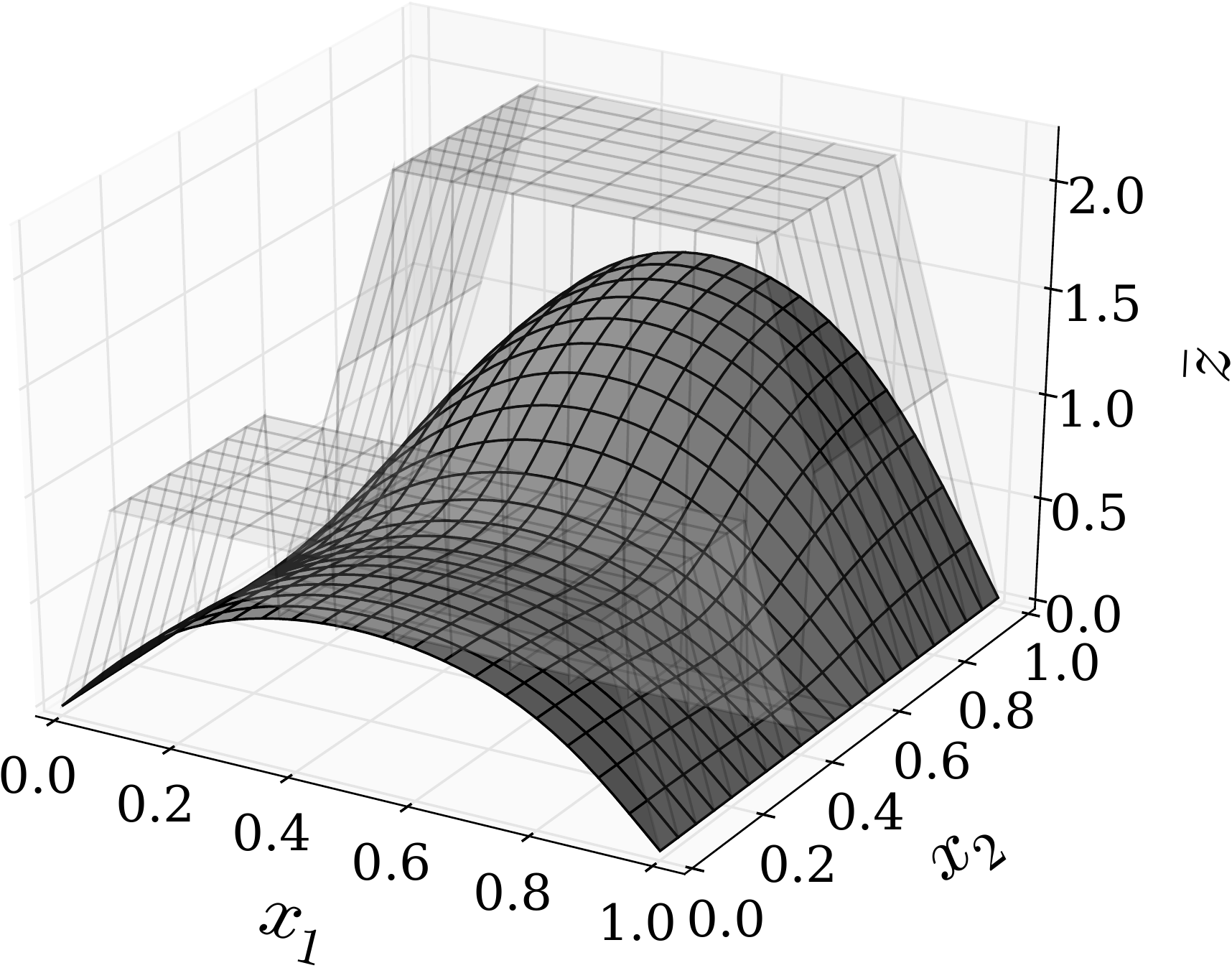}}
\subfloat[variance state]{\includegraphics[width=0.32\textwidth]{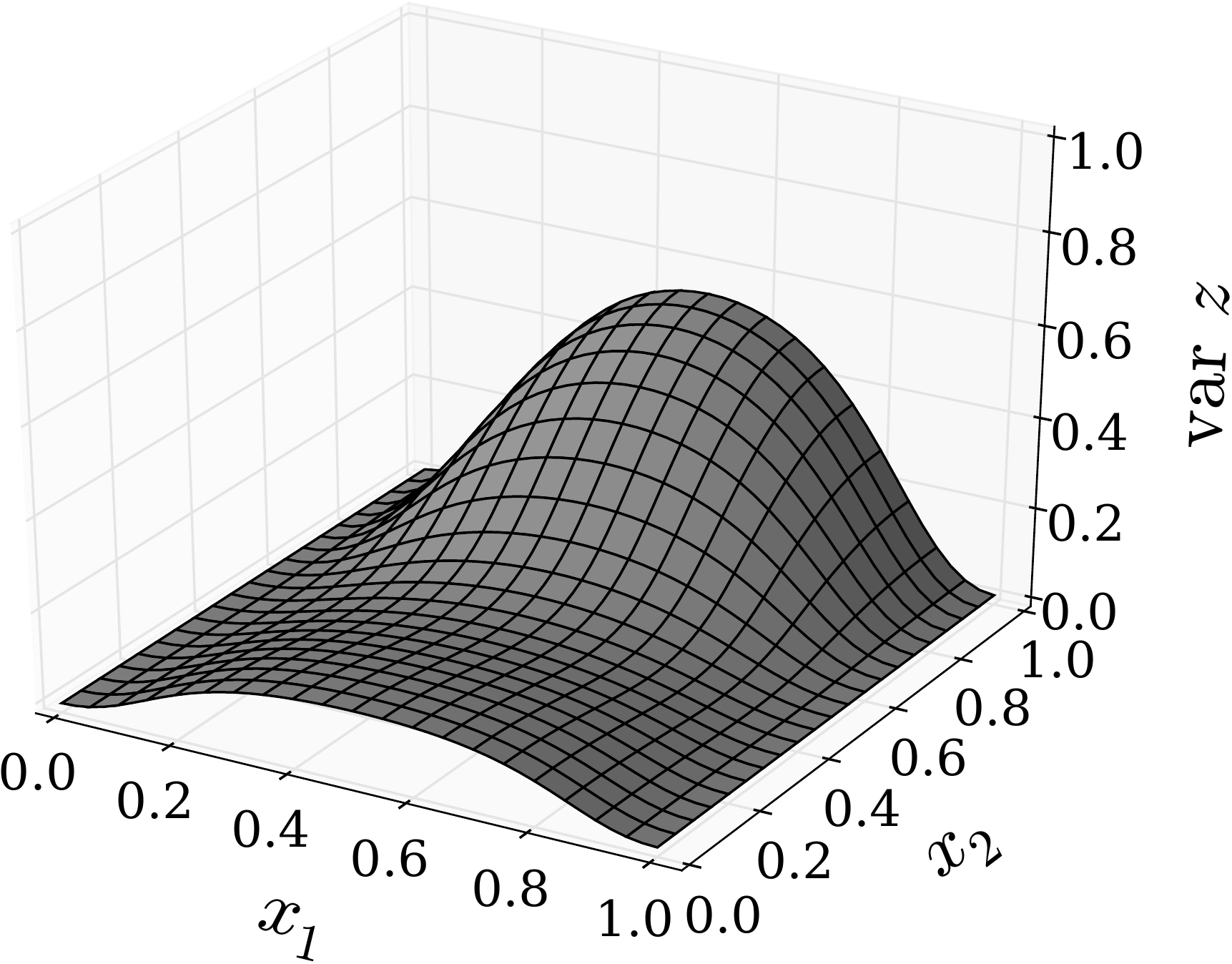}}
\subfloat[mean control]{\includegraphics[width=0.32\textwidth]{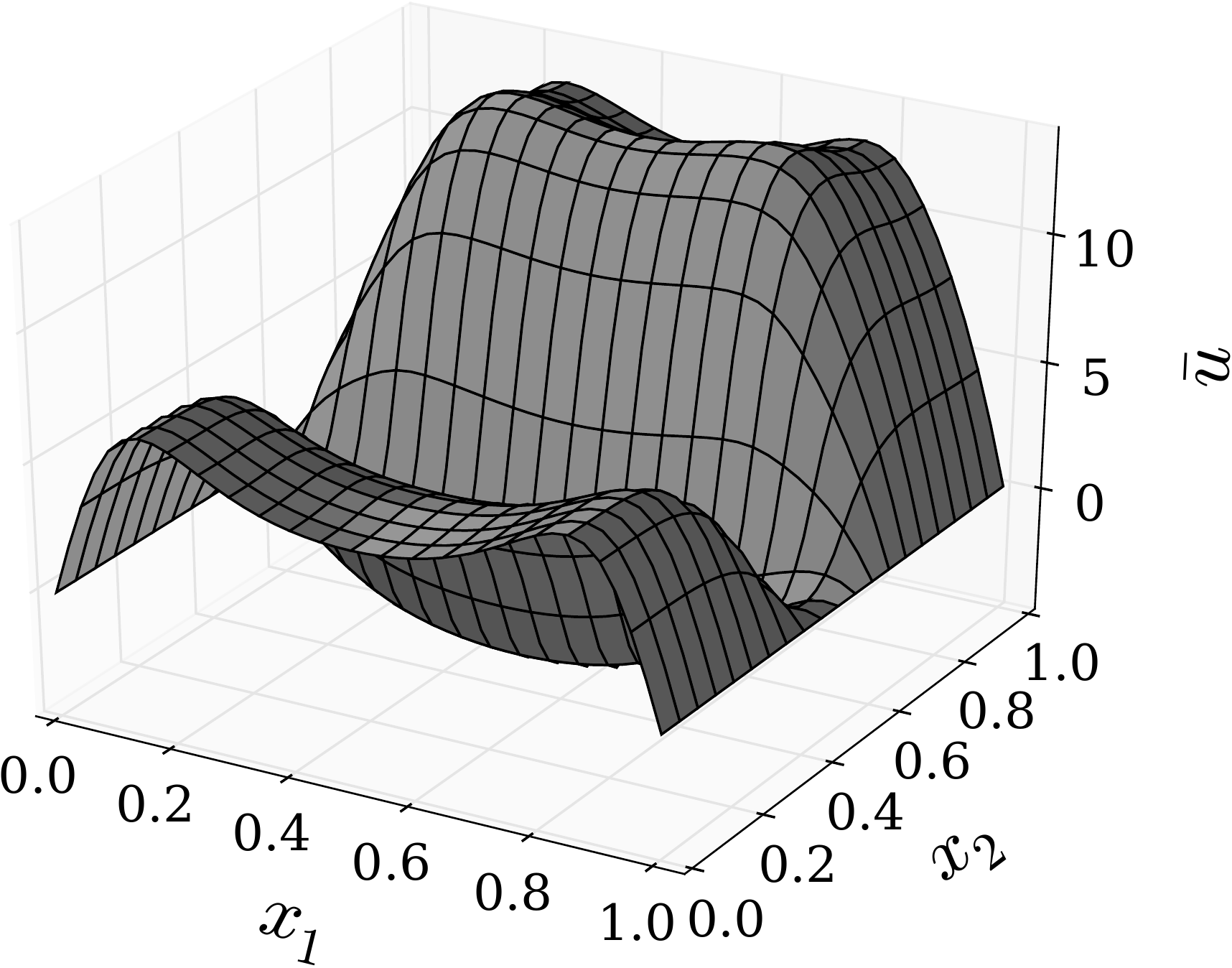}}
\caption{Mean and variance of the optimal state and mean
  control variable $\Bar{u}$ ($u' \ne 0$) associated with the cost
  functional $\mathcal{J}_{1}$ and computed with the stochastic Galerkin
  method with $\alpha= 1$, $\beta =1$, $\delta = 0$ and $\gamma =
  10^{-3}$.  The target~$\hat{z}$ is illustrated transparently in (a)
  for reference.}
\label{fig:det_dist_control_us1_B}
\end{figure}

The computed mean and second central moment of the state~$z$
and the computed mean optimal control~$\bar{u}$ are depicted in
Fig.~\ref{fig:det_dist_control_us1_A} for the case $\beta = 0$, and
in Fig.~\ref{fig:det_dist_control_us1_B} for the case $\beta = 1$.
A comparison between Fig.~\ref{fig:det_dist_control_A}, which corresponds
to the case $u'=0$, and Fig.~\ref{fig:det_dist_control_us1_A} shows that
the prescribed uncertainty on the control has only a limited effect on
the outcome of the optimal control problem for this example.  The mean
control $\Bar{u}$ is visibly unchanged and the variance of the state
variable increases slightly. This observation is also reflected in the
computed values for the cost functional (see Table~\ref{tab:cost}).

To provide control over the variance of the state variable, the
parameter~$\beta$ in the cost functional~\eqref{eq:costy} can be
increased.  Comparing the results in Fig.~\ref{fig:det_dist_control_us1_A}
($\beta = 0$) and Fig.~\ref{fig:det_dist_control_us1_B} ($\beta = 1$),
the peak variance is reduced, but the correspondence between the mean
state and the target is compromised.

\subsection{Distributive control with cost functional
\texorpdfstring{$\mathcal{J}_{2}$}{J2}}
\label{sec:dist_contr_2}

We now mirror the imperfect controller problem considered in
Section~\ref{sec:dist_contr_imperfect_1}, but for the cost functional
$\mathcal{J}_{2}$.  The cost functional $\mathcal{J}_{1}$ provides a
measure of the average distance between the state and target, whereas
the cost functional $\mathcal{J}_{2}$ provides a measure of the distance
between the mean state and the target.

\begin{figure}
\centering
\subfloat[mean state]{\includegraphics[width=0.32\textwidth]{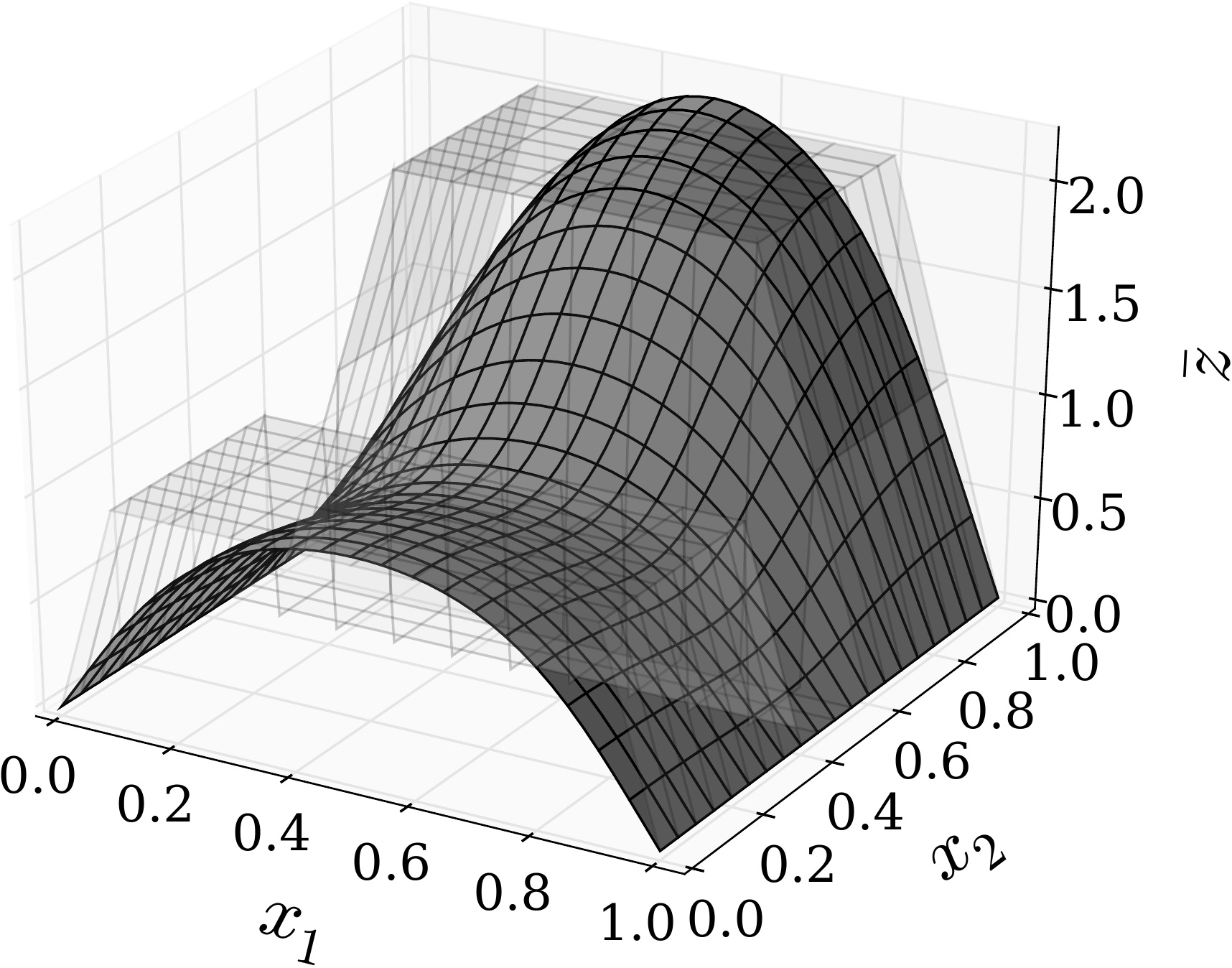}}
\subfloat[variance state]{\includegraphics[width=0.32\textwidth]{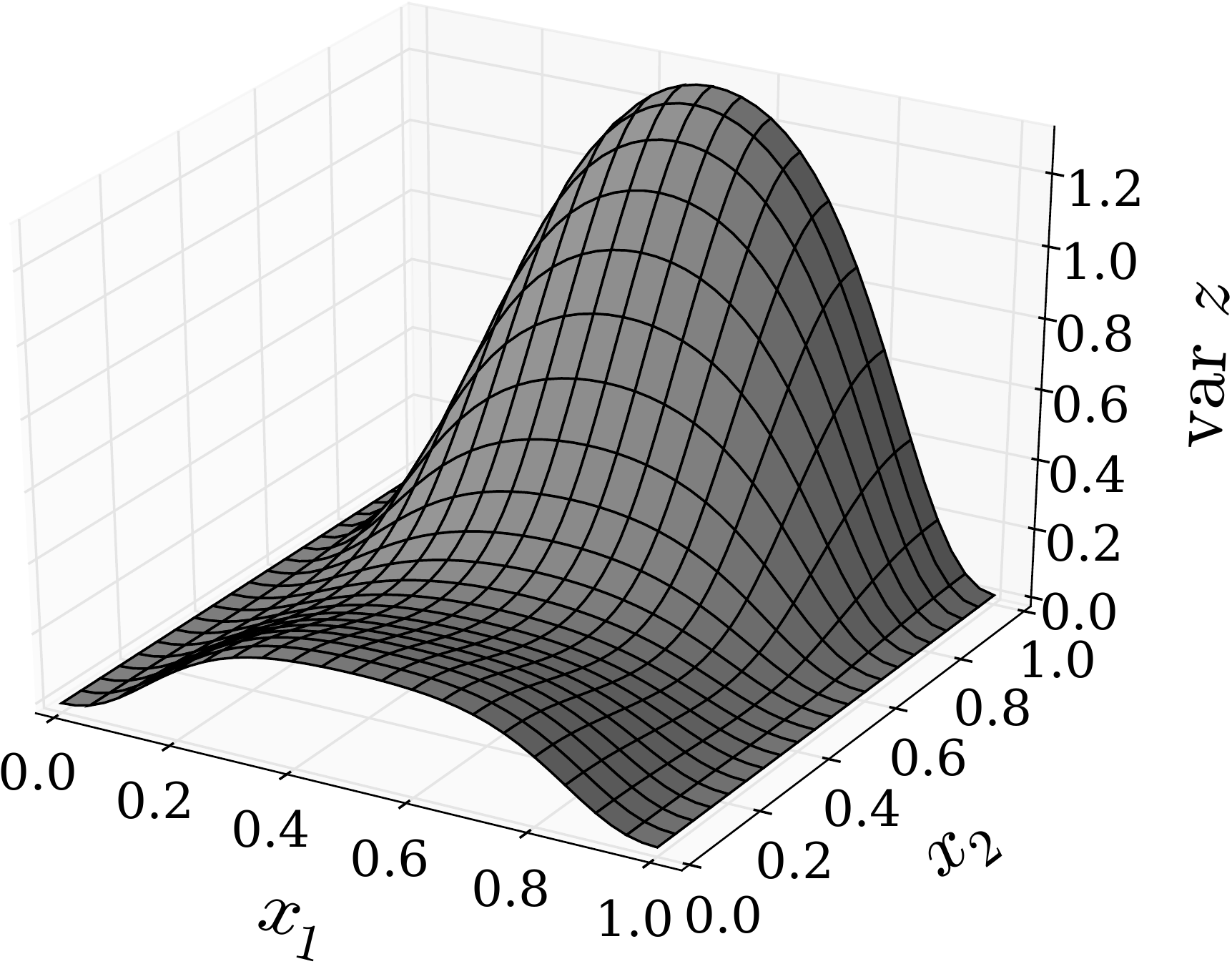}}
\subfloat[mean control]{\includegraphics[width=0.32\textwidth]{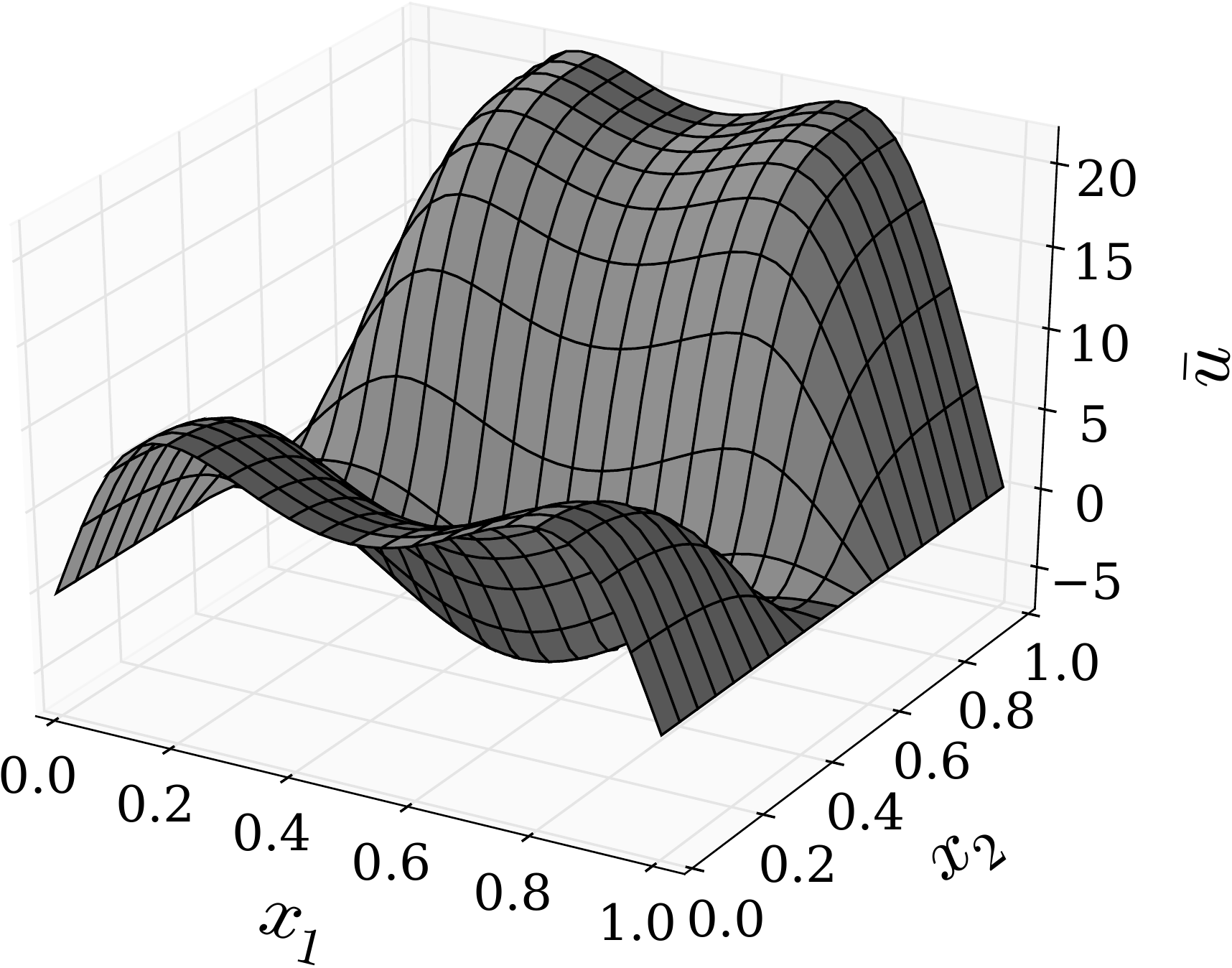}}
\caption{Mean and variance of the optimal state and mean
  control variable $\Bar{u}$ ($u' \ne 0$) associated with the
  cost functional $\mathcal{J}_{2}$ and computed with the stochastic
  Galerkin method with $\alpha= 1$, $\beta = \delta = 0$ and $\gamma
  = 10^{-3}$.  The target~$\hat{z}$ is illustrated transparently in (a)
  for reference.}
\label{fig:det_dist_control_alt_A}
\end{figure}

\begin{figure}
\centering
\subfloat[mean state]{\includegraphics[width=0.32\textwidth]{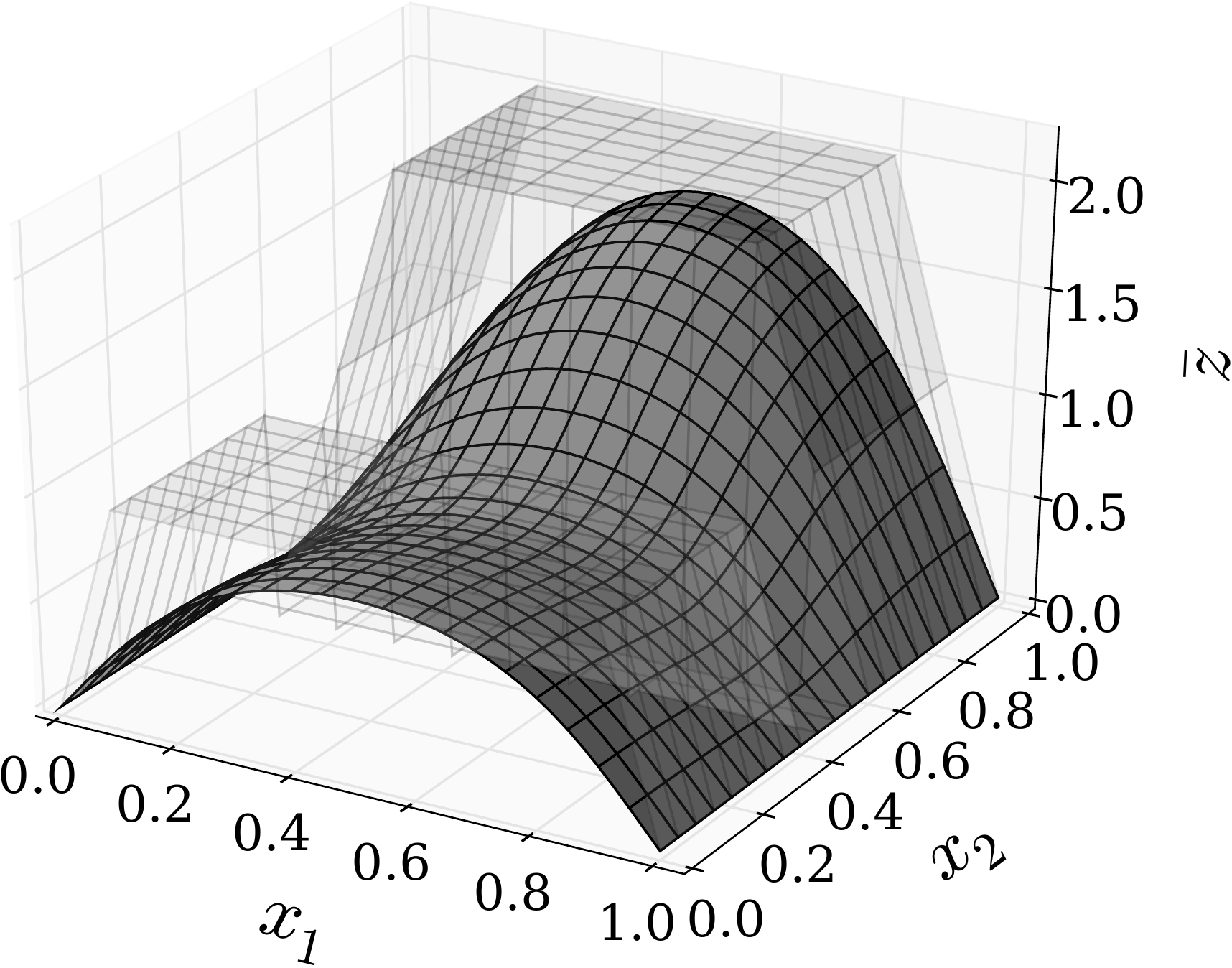}}
\subfloat[variance state]{\includegraphics[width=0.32\textwidth]{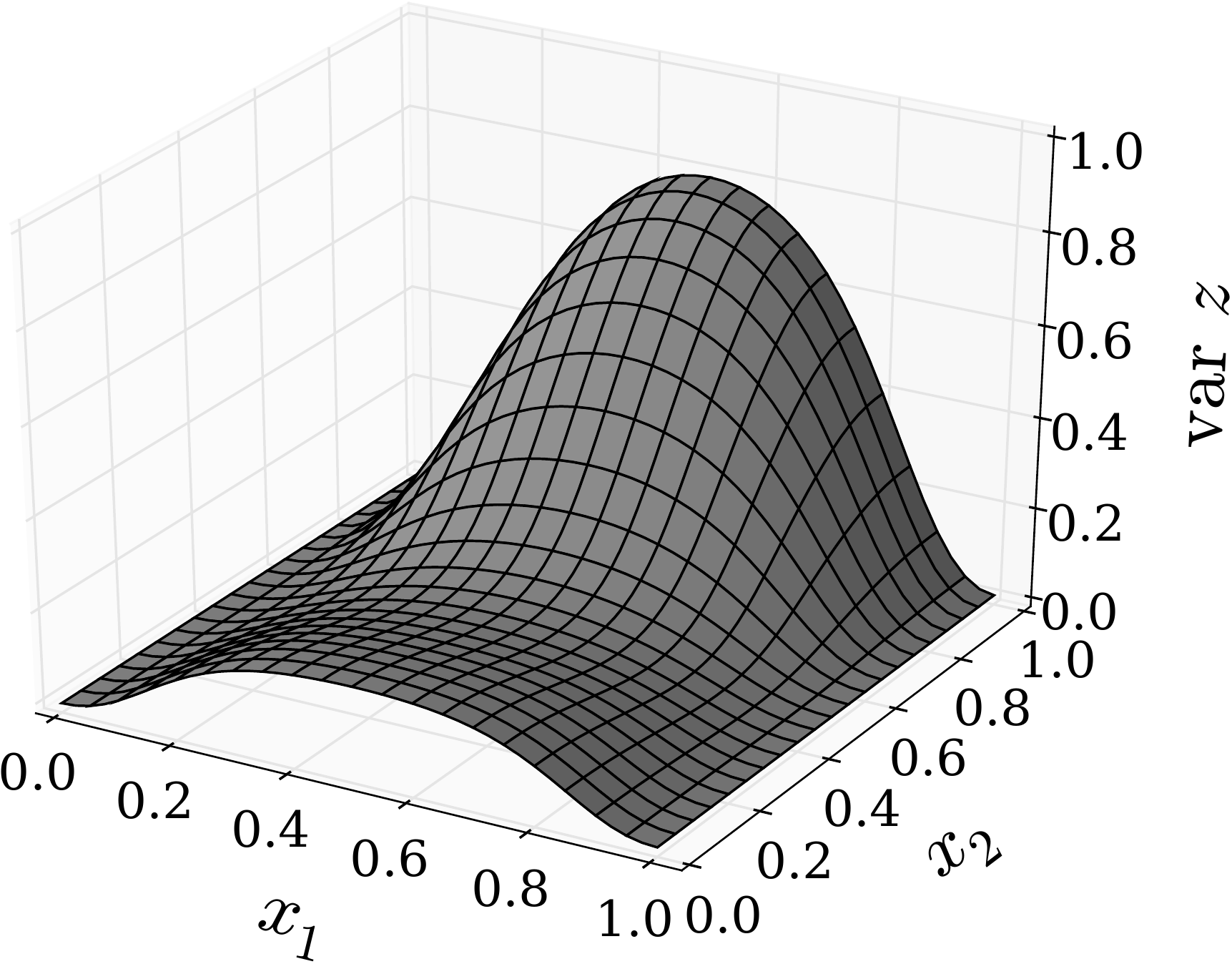}}
\subfloat[mean control]{\includegraphics[width=0.32\textwidth]{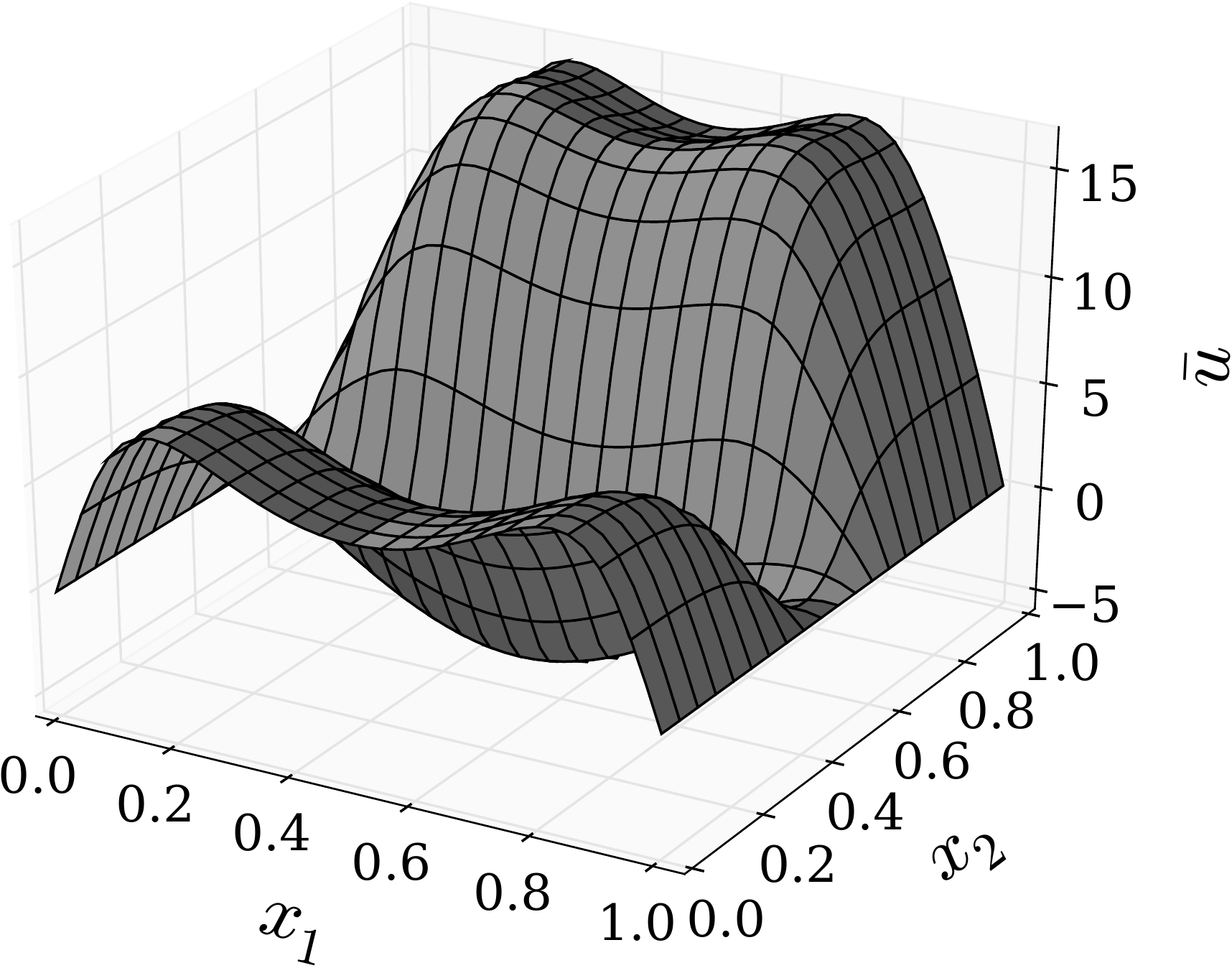}}
\caption{Mean and variance of the optimal state and mean
  control variable $\Bar{u}$ ($u' \ne 0$) associated with the cost
  functional $\mathcal{J}_{2}$ and computed with the stochastic Galerkin
  method with $\alpha= 1$, $\beta = 1$, $\delta = 0$ and $\gamma =
  10^{-3}$.  The target~$\hat{z}$ is illustrated transparently in (a)
  for reference.}
\label{fig:det_dist_control_alt_B}
\end{figure}

For the examples we adopt $\alpha = 1$, $\gamma = 10^{-3}$ and $\delta
= 0$ in $\mathcal{J}_{2}$.  For the case $\beta = 0$, which implies
no extra control over the variance of the response, the mean and
the variance of the computed state variable $z$ and the computed
deterministic part of the control signal $\Bar{u}$ are shown in
Fig.~\ref{fig:det_dist_control_alt_A}.  Compared to the case presented
in Fig.~\ref{fig:det_dist_control_us1_A}, the mean of the computed state
variable is a better approximation of the target, while the variance
of the state variable is significantly larger.  A similar result was
observed for the case~$u' = 0$.  Fig.~\ref{fig:det_dist_control_alt_B}
shows the computed results for the case~$\beta=1$. Increasing $\beta$
reduces the variance of the response, but comes at the expense of the
approximation of the target by the mean state.

The advantage of using the cost functional $\mathcal{J}_{2}$ over
$\mathcal{J}_{1}$ is that it permits a greater tuning of the relative
importance of the mean response versus the variance of the response,
via the parameters $\alpha$ and~$\beta$.  In the case of the cost
functional $\mathcal{J}_{1}$, the variance of the state variable is
already implicitly minimised by the $\alpha$ term as it is a norm over
$L^{2}(D) \times L^{2}_{\rho}(\Omega)$. For the $\mathcal{J}_{1}$ case,
increasing $\beta$ only entails an additional contribution of the variance
to the cost functional.

\begin{table}
\centering \small
\begin{tabular}{|l|ccc|}
\hline
\multicolumn{1}{|c}{}&&&\\[-0.3cm]
\multicolumn{1}{|c}{}& $\mathcal{J}(z,u)$ &
$\norm{z - \hat{z}}^2_{L^2(D)\otimes L^2_\rho(\Gamma)}$ &
$\:\norm{\mathrm{std}(z)}_{L^2(D)}^2$\\[0.1cm]
\hline
&&&\\[-0.3cm]
\textbf{deterministic control $u= \bar{u}(x)$, $u'(x,y)=0$}&&&\\
cost functional~$\mathcal{J}_1$, $\beta = 0$, $\gamma = 10^{-5}$&
 $2.083 \times 10^{-1}$ & $4.022 \times 10^{-1}$ & $ 2.562 \times 10^{-1} $\\
cost functional~$\mathcal{J}_{1}$, $\beta =0$, $\gamma = 10^{-3}$ &
 $2.911\times 10^{-1}$ & $5.078 \times 10^{-1}$ & $1.845  \times 10^{-1} $\\[0.1cm]
\hline
&&&\\[-0.3cm]
\textbf{unknown mean control $\bar{u}(x)$, var($u'$)$ = 1$}&&&\\
  cost functional~$\mathcal{J}_{1}$, $\beta =  0$,
  $\gamma = 10^{-3}$ &
$  2.956\times 10^{-1}$ & $ 5.160\times 10^{-1}$ &$ 1.927\times 10^{-1}$\\
 cost functional~$\mathcal{J}_{1}$, $ \beta = 1$,
$\gamma = 10^{-3}$&
$3.767\times 10^{-1}$ & $ 5.636\times 10^{-1}$ &  $1.367\times 10^{-1}$ \\
cost functional~$\mathcal{J}_2$, $\beta =  0$,
$\gamma = 10^{-3}$ &
$1.764 \times 10^{-1} $ & $2.353\times 10^{-1}$&$ 2.957\times 10^{-1}$\\
cost functional~$\mathcal{J}_2$, $ \beta = 1$,
$\gamma = 10^{-3}$&
$2.956 \times 10^{-1}$ & $3.233 \times 10^{-1}$ & $1.927 \times 10^{-1}$\\
\hline
\end{tabular}
\caption{Summary of the cost functional, tracking error and standard
  deviation of the state variable for the considered optimal control
  problems with a distributive control function.}
\label{tab:cost}
\end{table}

\subsection{Boundary control with cost functional
\texorpdfstring{$\mathcal{J}_1$}{J1}}
\label{ssec:num_bc}

In practical applications it is often more plausible to control boundary
values rather than the source term. In this section, we consider the
same model problem as in Section~\ref{ssec:dist_contr_1}, but now an
optimal boundary flux control is computed (now $\gamma = 0$ in the
cost functionals).  The right-hand side function in the constraint
equation~\eqref{eq:constr1} is set to $u = 5$.  The control boundary
flux is computed on the `upper' ($x_{2} = 1$) and `lower' ($x_{2} =
0$) boundaries.

\subsubsection{Perfect controller case}

\begin{figure}
\centering
\subfloat[mean state]{\includegraphics[width=0.32\textwidth]{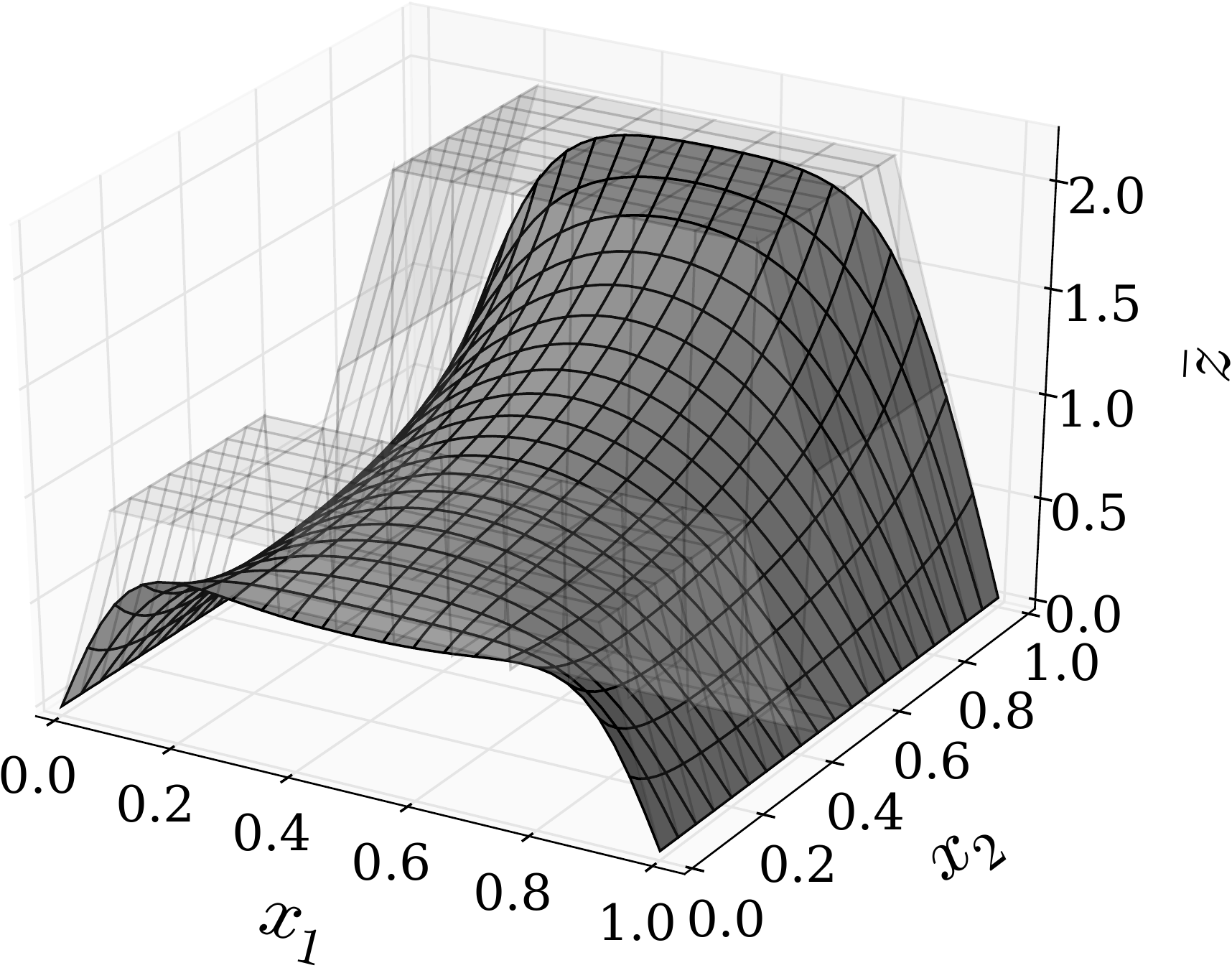}}
\subfloat[variance state]{\includegraphics[width=0.32\textwidth]{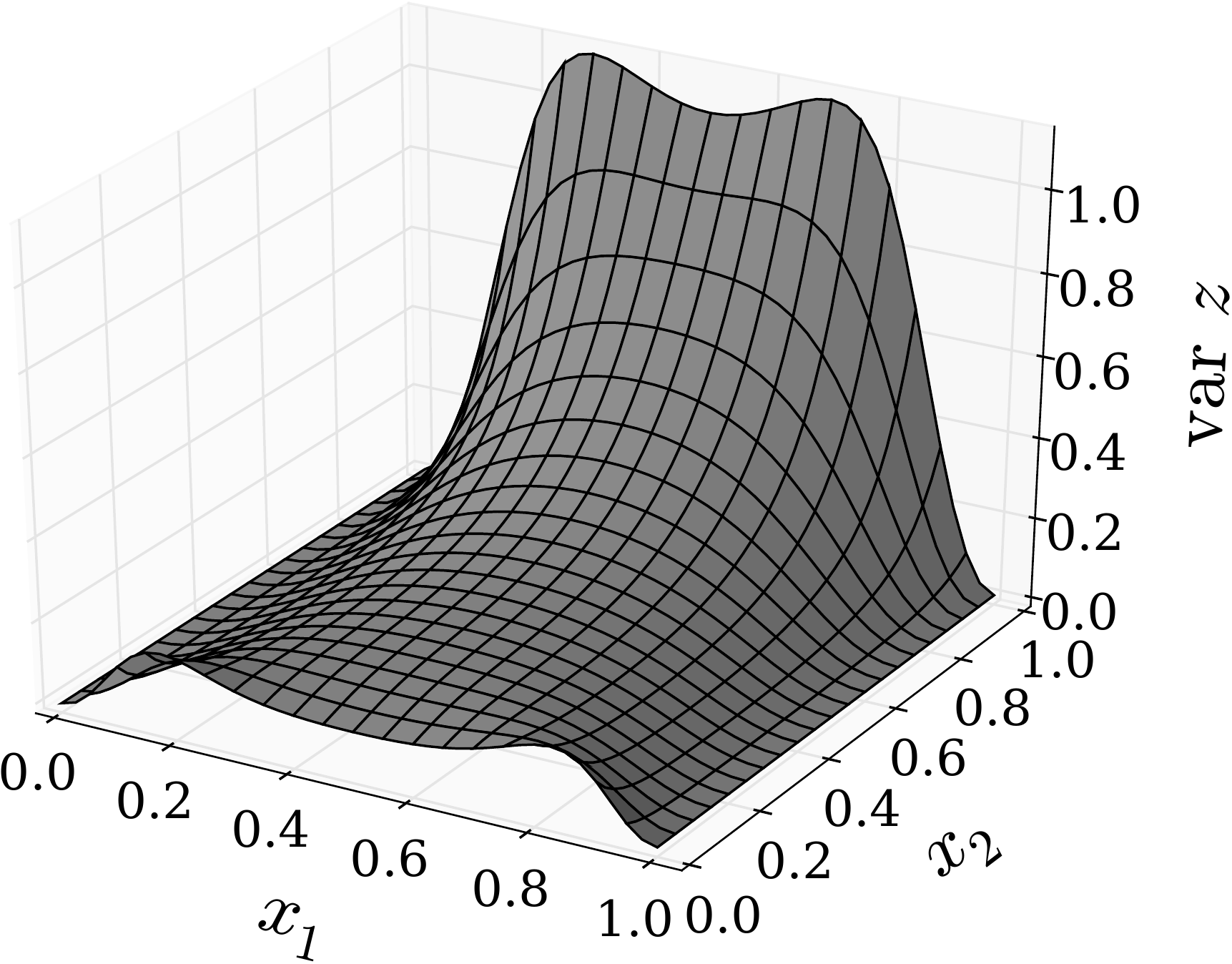}}
\subfloat[deterministic control]{\includegraphics[width=0.32\textwidth]{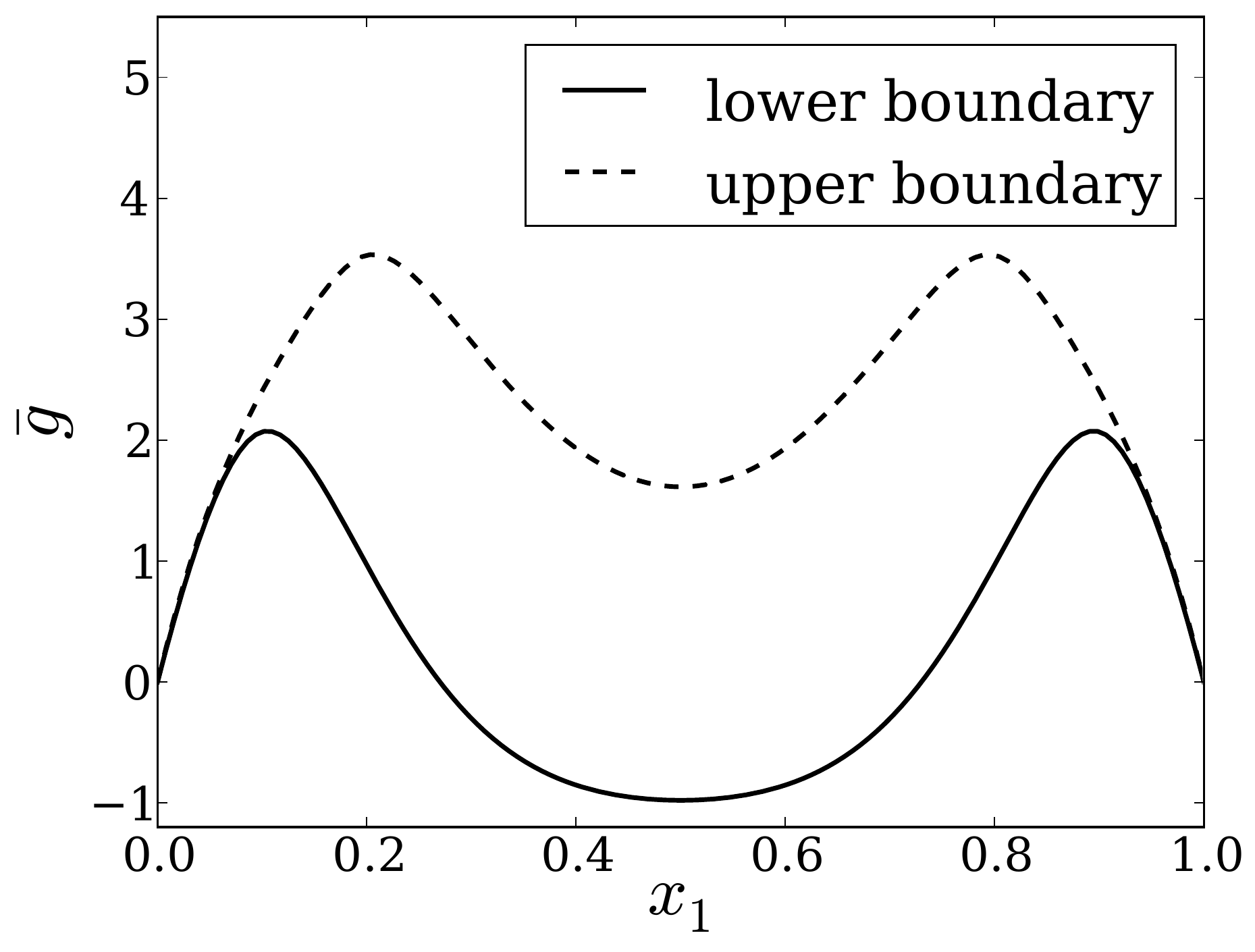}}
\caption{Mean and variance of the optimal state and deterministic
  control variable $g = \bar{g}$ ($g' = 0$) associated with the cost
  functional $\mathcal{J}_{1}$ and computed with the stochastic Galerkin
  method with $\alpha= 1$, $\gamma = 0$, $\delta= 10^{-3}$ and $\beta = 0$.
  The target~$\hat{z}$ is illustrated transparently in (a) for reference.}
 \label{fig:det_bc_controlbeta_A}
\end{figure}

\begin{figure}
\centering
\subfloat[mean state]{\includegraphics[width=0.32\textwidth]{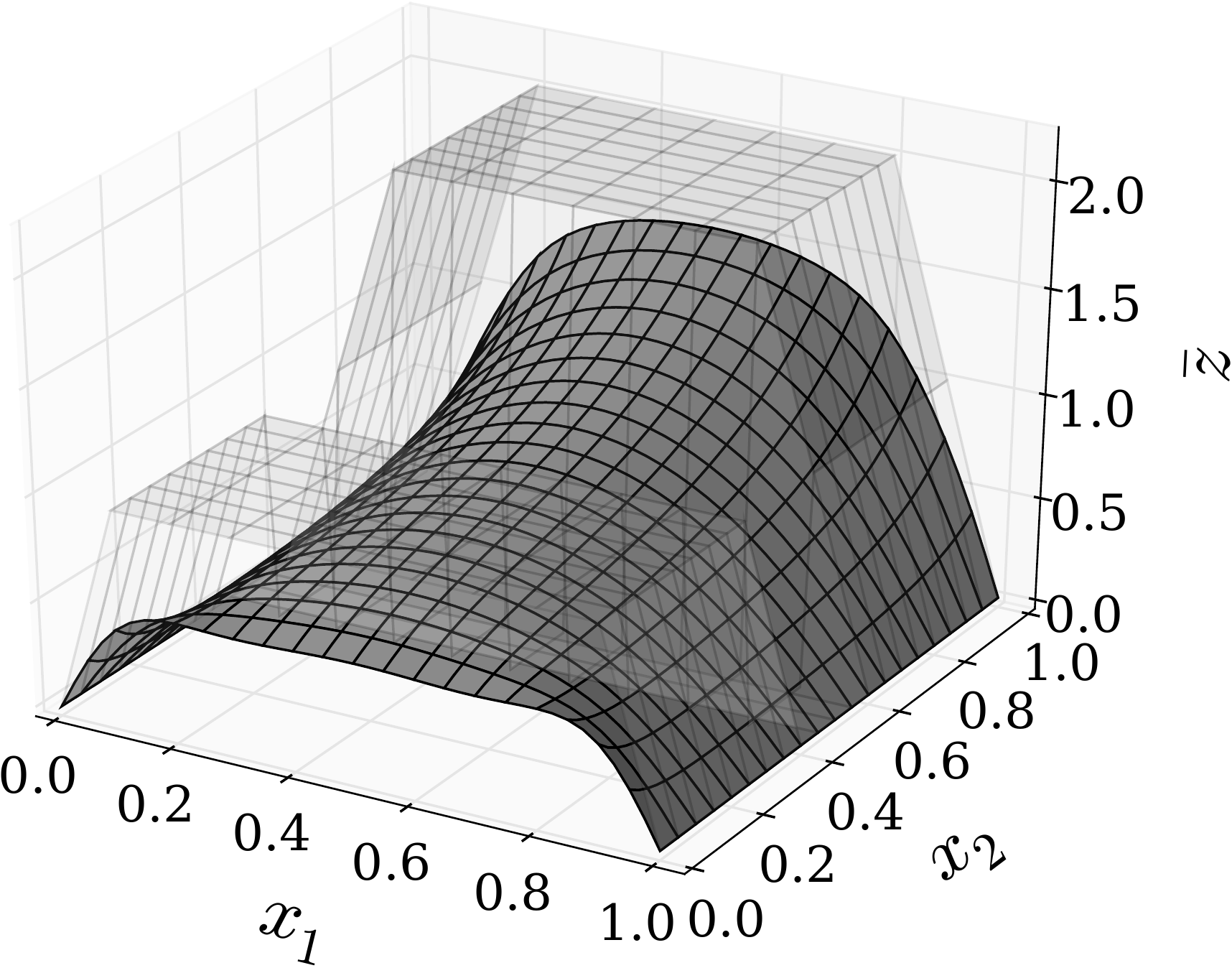}}
\subfloat[variance state]{\includegraphics[width=0.32\textwidth]{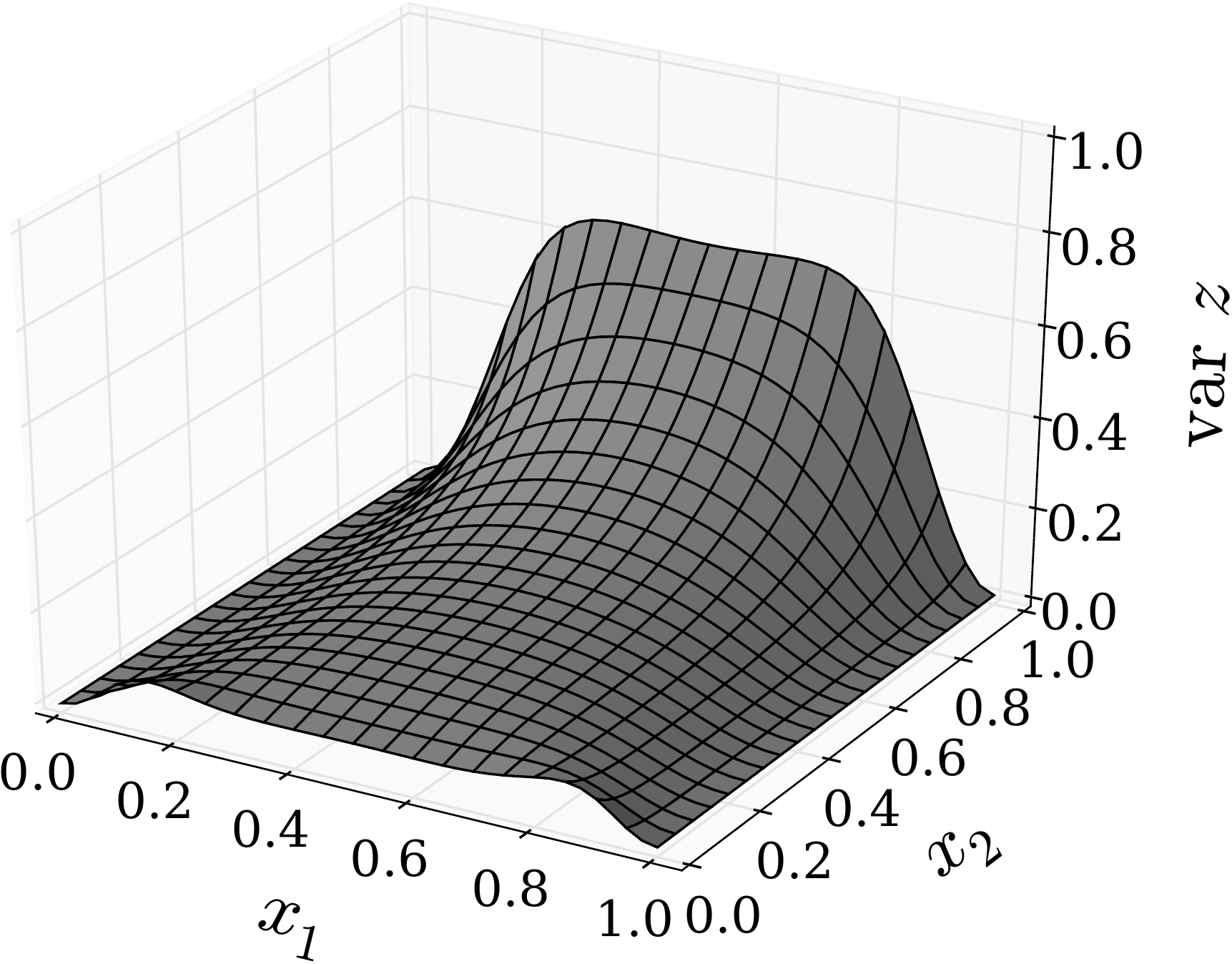}}
\subfloat[deterministic control]{\includegraphics[width=0.32\textwidth]{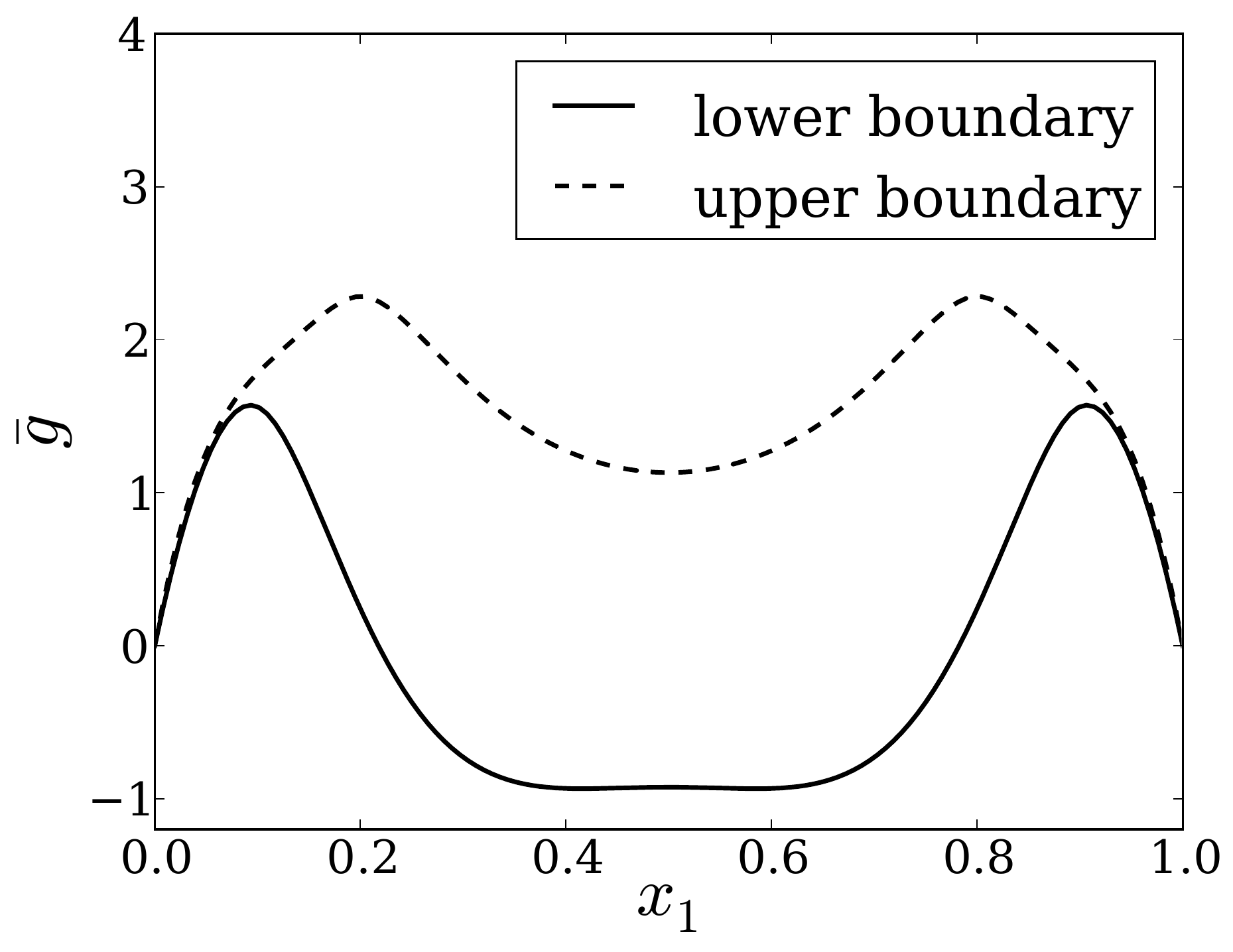}}
\caption{Mean and variance of the optimal state and deterministic
  control variable $g = \bar{g}$ ($g' = 0$) associated with the cost
  functional $\mathcal{J}_{1}$ and computed with the stochastic Galerkin
  method with $\alpha= 1$, $\gamma = 0$, $\delta= 10^{-3}$ and $\beta = 1$.
  The target~$\hat{z}$ is illustrated transparently in (a) for reference.}
\label{fig:det_bc_controlbeta_B}
\end{figure}

Fig.~\ref{fig:det_bc_controlbeta_A} presents the computed mean state,
variance of the state and the boundary flux control~$\Bar{g}$ for the
parameters $\alpha= 1$, $\beta = \gamma = 0$ and $\delta= 10^{-3}$,
and with~$g' = 0$.  The corresponding values of the cost functional
and tracking errors are given in Table~\ref{tab:cost2}.  As can be
anticipated, correspondence between the optimal state and target is poorer
in the case of a boundary control than in the case of a distributive
control (see Fig.~\ref{fig:det_dist_control_B}).  The variance of the
state variable in Fig.~\ref{fig:det_bc_controlbeta_A} can be countered
by increasing the parameter~$\beta$. Results for the case $\beta =
1$ are presented in Fig.~\ref{fig:det_bc_controlbeta_B}. The reduction
in the variance is accompanied by a deterioration in the approximation
of the target function.

\subsubsection{Imperfect controller case}

We now consider the impact of an imperfect controller by decomposing the
boundary control~$g$ according to~\eqref{eq:u}. The zero-mean stochastic
function $g'$ is modelled by a three-term Karhunen--Lo\`eve expansion
based on a Gaussian field with an exponential covariance function,
a variance of~$0.25$ and unit correlation length.

\begin{figure}
\centering
\subfloat[mean state]{\includegraphics[width=0.32\textwidth]{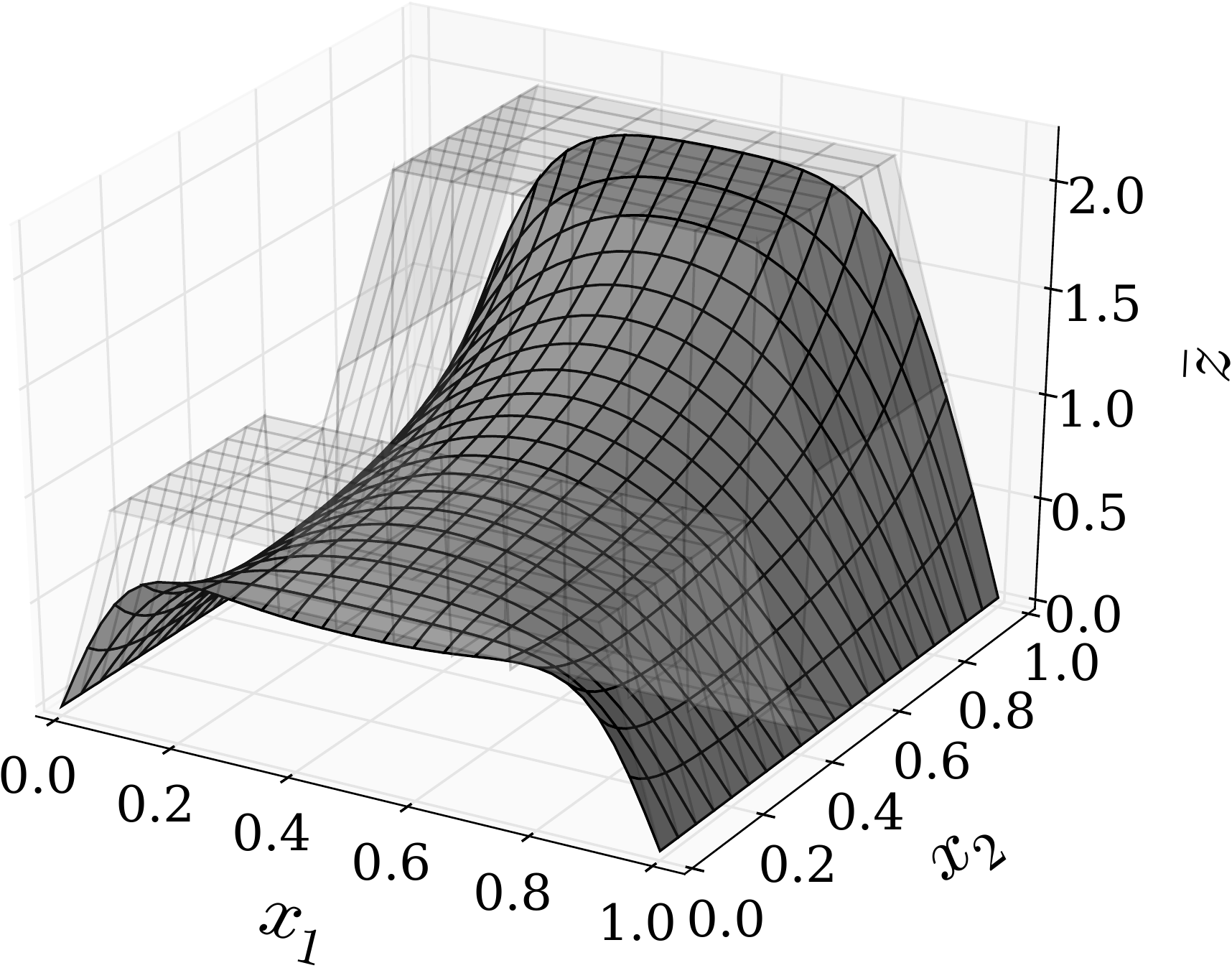}}
\subfloat[variance state]{\includegraphics[width=0.32\textwidth]{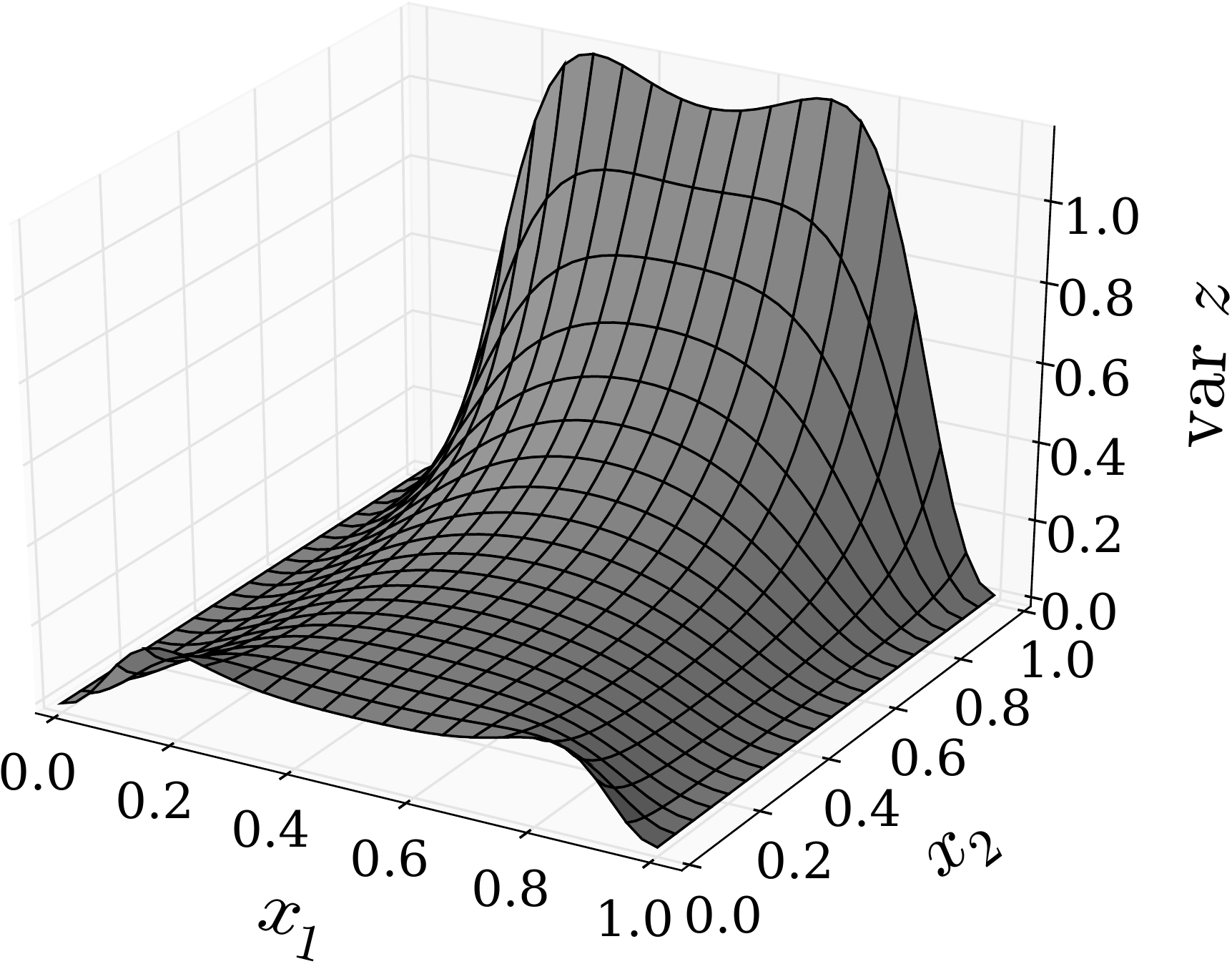}}
\subfloat[mean control]{\includegraphics[width=0.32\textwidth]{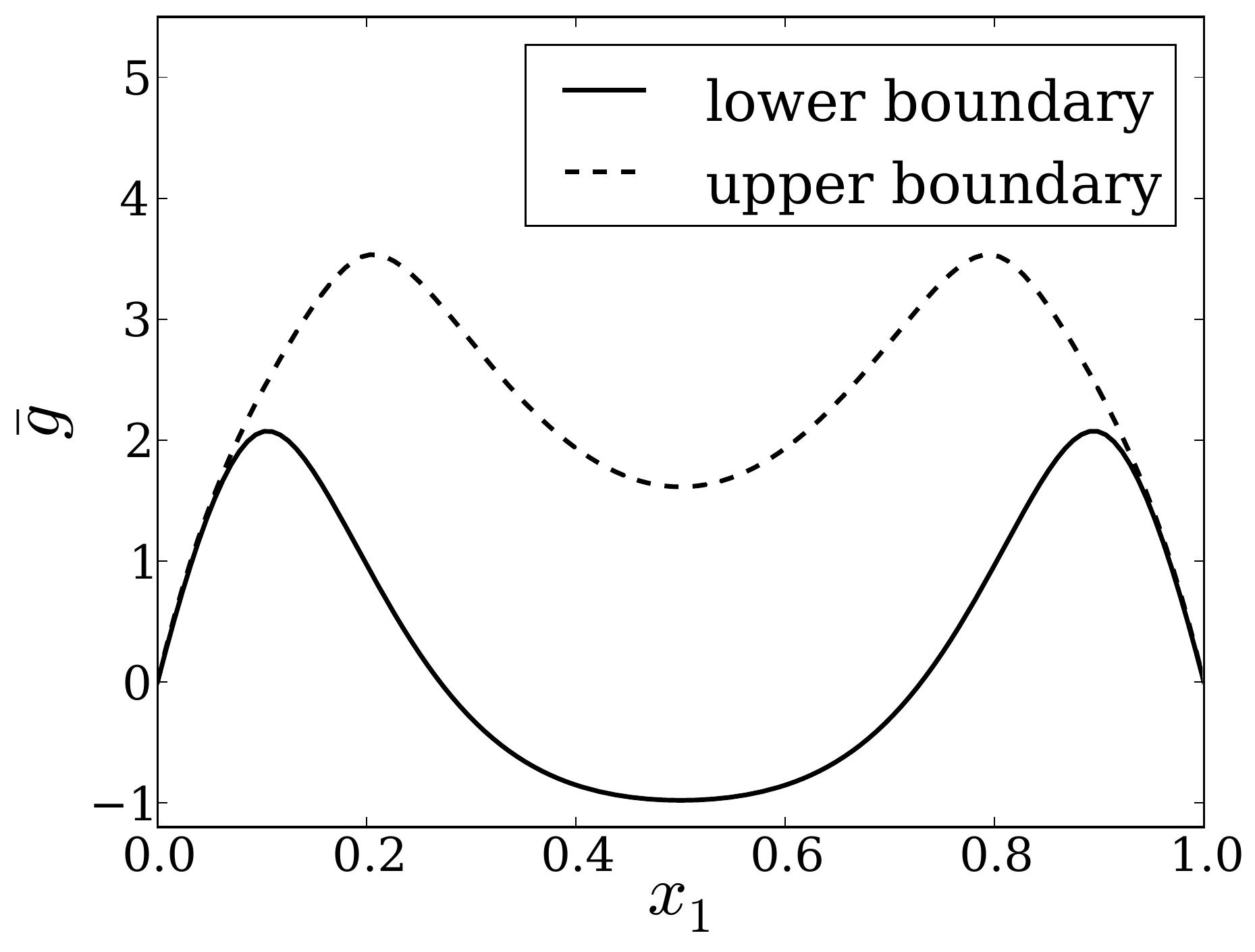}}
\caption{Mean and variance of the optimal state and mean control variable
  $\bar{g}$ ($g' \ne 0$) associated with the cost functional $\mathcal{J}_1$
  and computed with the stochastic Galerkin method with $\alpha= 1$,
  $ \beta =\gamma = 0$ and $\delta= 10^{-3}$. The target~$\hat{z}$ is
  illustrated transparently in (a) for reference.}
 \label{fig:det_bc_control_gs}
\end{figure}

Fig.~\ref{fig:det_bc_control_gs} presents the optimal mean boundary
control $\Bar{g}$ along with the first moments of the state variable.
The values of the cost functional and tracking errors are given
in Table~\ref{tab:cost2}. The mean optimal control~$\Bar{g}$
and the mean state appear visibly equal to the results in
Fig.~\ref{fig:det_bc_controlbeta_A} for the case where a perfect control
device, i.e., with $g'= 0$, is modelled, while the variance of~$z$
is slightly larger.

\begin{table}
\centering\small
\begin{tabular}{|l|ccc|}
\hline
\multicolumn{1}{|c}{}&&&\\[-0.3cm]
\multicolumn{1}{|c}{}& $\mathcal{J}(z,g)$ &
$\norm{z - \hat{z}}^2_{L^2(D)\otimes L^2_\rho(\Gamma)}$ &
 $\:\norm{\mathrm{std}(z)}_{L^2(D)}^2$\\[0.1cm]
\hline
&&&\\[-0.3cm]
\textbf{deterministic control $g = \bar{g}(x)$, $g'(x,y) = 0$}  &&&\\
 cost functional $\mathcal{J}_1$, $\beta =  0$, $\delta = 10^{-3}$ &
$2.711 \times 10^{-1}$ &$ 5.421\times 10^{-1}$ &   $2.091 \times 10^{-1}$\\
cost functional $\mathcal{J}_1$, $\beta = 1$, $\delta = 10^{-3}$&
 $ 3.593\times 10^{-1}$&  $5.757\times 10 ^{-1}$& $  1.428\times 10^{-1}$ \\
\hline
&&&\\[-0.3cm]
\textbf{unknown mean control $\bar{g}(x)$, var($g'$)$=0.25$} &&&\\
cost functional $\mathcal{J}_1$, $\beta = 0$, $\delta = 10^{-3}$ &
$ 2.753 \times 10^{-1}$ & $ 5.499\times 10^{-1}$ & $2.168\times 10^{-1}$\\
cost functional~$\mathcal{J}_1$, $\beta = 1$, $\delta = 10^{-3}$ &
$3.673\times 10^{-1}$ & $5.835 \times 10^{-1}$& $ 1.506 \times 10^{-1}$ \\
\hline
\end{tabular}
\caption{Summary of the cost functional, tracking error and
  standard deviation of the state variable for the considered optimal
  control problems with a boundary control function.}
\label{tab:cost2}
\end{table}

\section{Stochastic inverse examples}
\label{sec:num2}

When the additive structure of the control function is not enforced
and $u$ is permitted to be stochastic and unknown, the optimal control
problem effectively becomes a stochastic inverse problem. In this case,
the stochastic properties of $u$ are unknown, but will be computed.
The problem is: given an observation $\Hat{z}$ of a system that must
obey the constraints in~\eqref{eq:constr1}--\eqref{eq:constr3}, find the
stochastic source term $u$ that would induce the response~$\Hat{z}$.
Computed higher moments of $u$ provide information on the uncertainty
of the driving term. This formulation of the problem is not useful for
control problems, since the properties of a control device are considered
to be known and it is unclear how a stochastic $u$ could be used as a
control signal. The mean of $u$ could be taken as the control, but this will
not in general be the optimal control and computing the uncertainty in
the system response would require an additional computation.

\begin{figure}
\centering
\subfloat[mean target]{\includegraphics[width=0.32\textwidth]{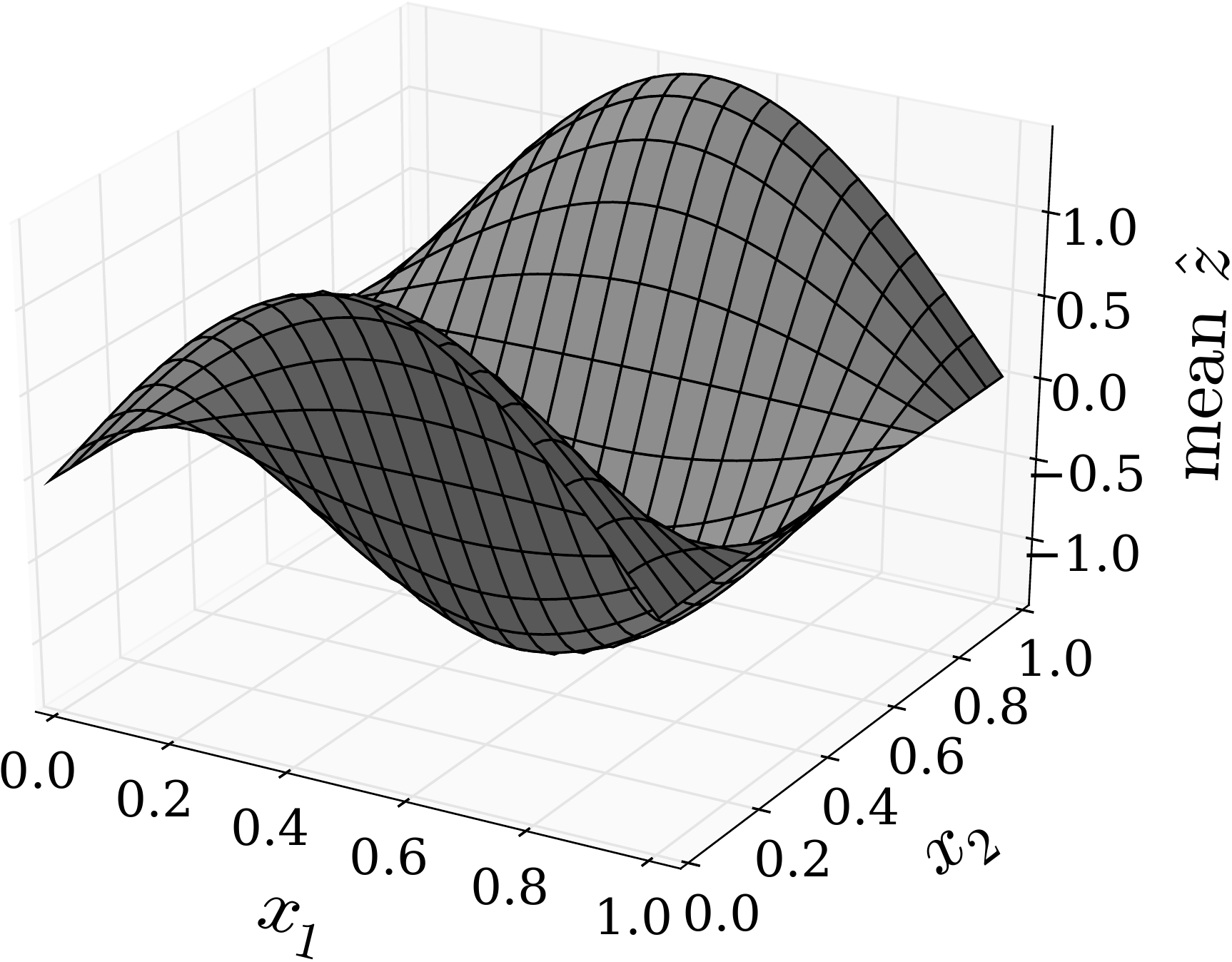}}
\subfloat[variance target]{\includegraphics[width=0.32\textwidth]{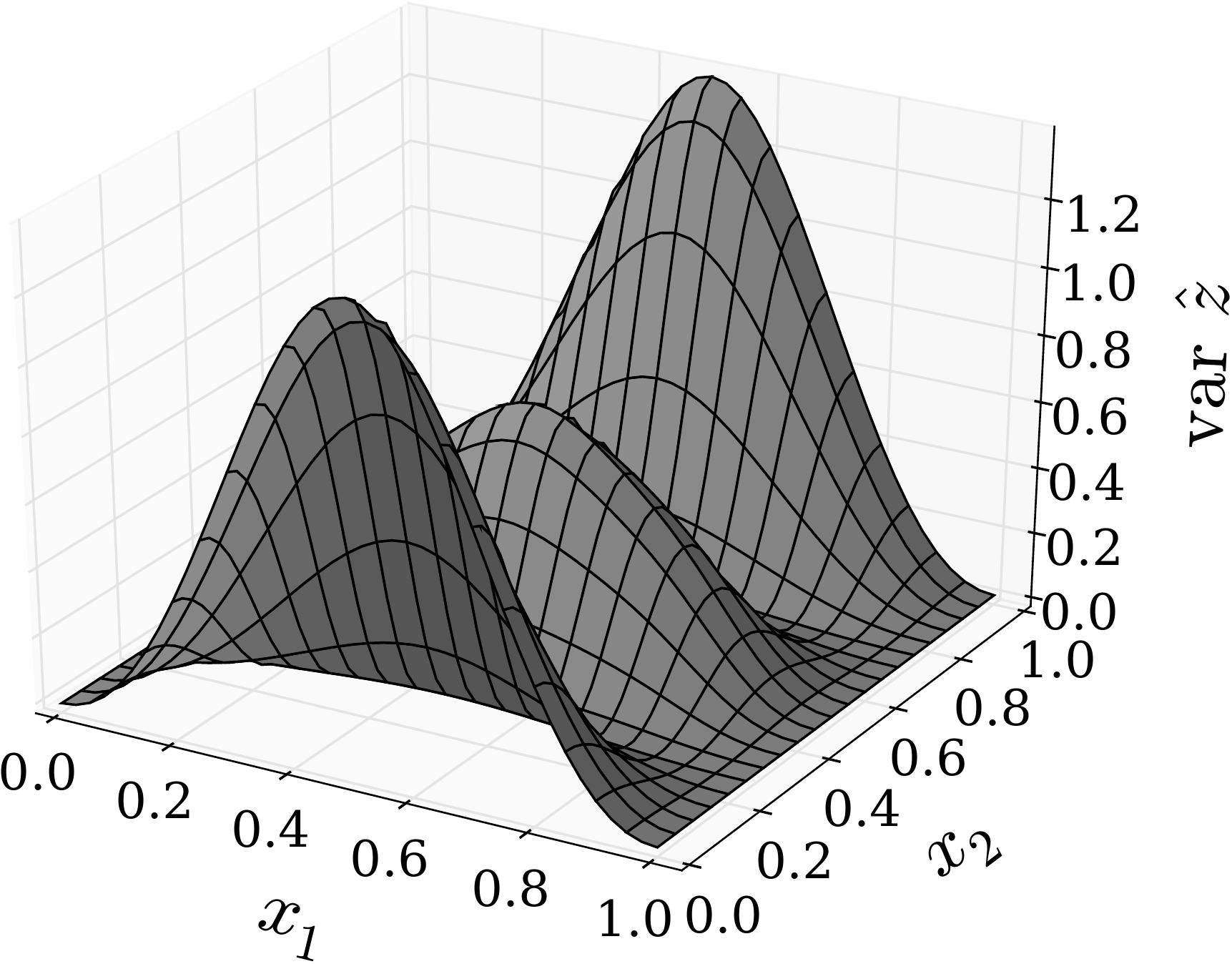}}
\subfloat[deterministic source]{\includegraphics[width=0.32\textwidth]{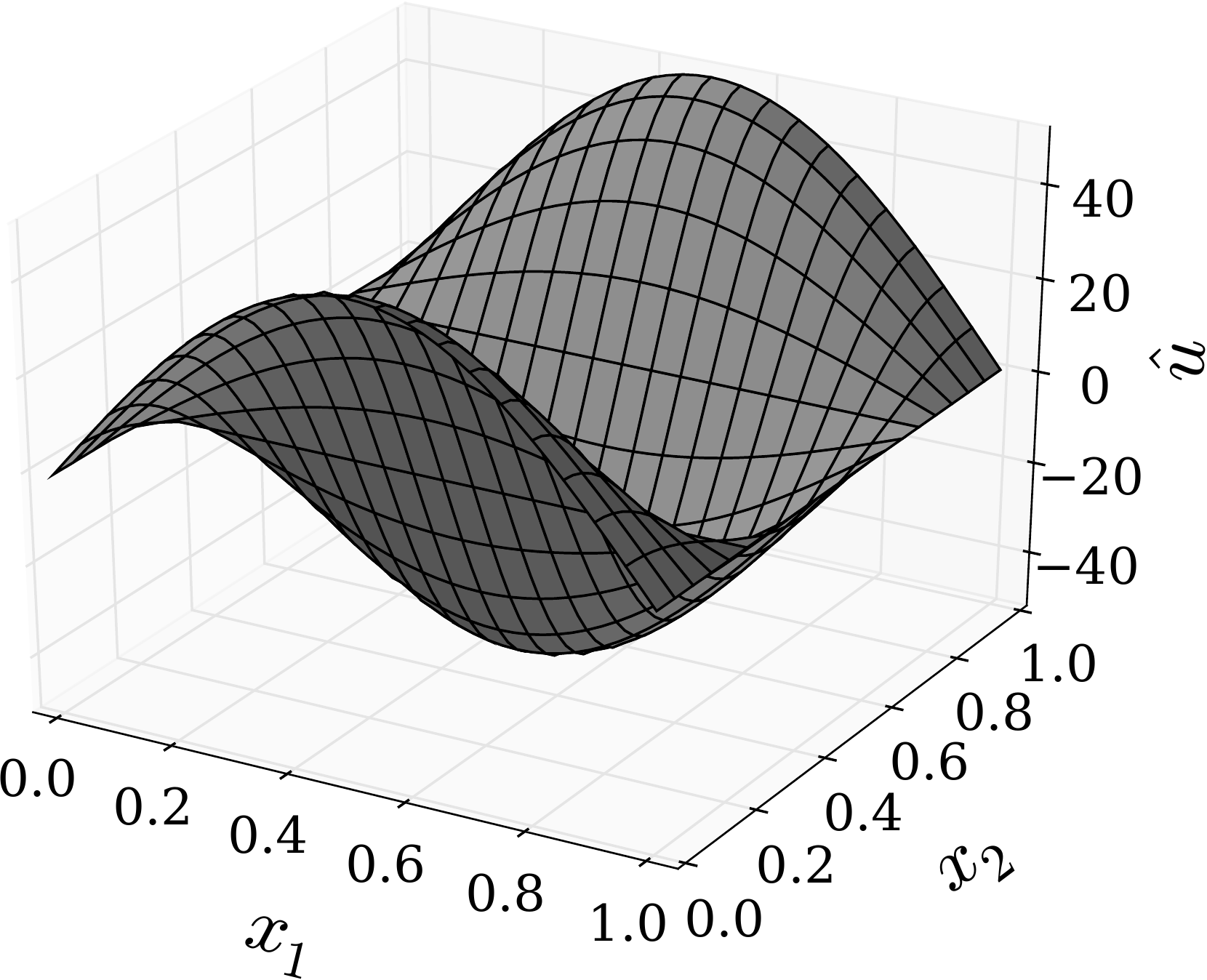}}
\caption{(a)--(b) Mean and variance of the target function $\hat{z}$ for solving
 a stochastic inverse problem. (c) Deterministic source $\hat{u}$ used for
 constructing $\hat{z}$ as the solution of the forward problem.}
\label{fig:target_inverse}
\end{figure}

Stochastic inverse problems associated with the cost
functional $\mathcal{J}_1$ or $\mathcal{J}_2$ and the constraint
equations~\eqref{eq:constr1y}--\eqref{eq:constr3y} are solved using the
same computational domain, boundary conditions, representation of $\kappa$
and discretisation parameters as in Section~\ref{sec:num1}.  For inverse
problems, the collocation systems remain decoupled for $\beta = 0$ in the
cost functional $\mathcal{J}_1$.  We present some simple examples using
a target function~$\Hat{z}$ that is computed by solving the stochastic
forward problem in \eqref{eq:constr1y}--\eqref{eq:constr3y} with the
deterministic source function
\begin{equation}
\label{eq:source_inverse}
 \Hat{u} := 50  \sin(\pi x_1) \cos(2\pi x_2).
\end{equation}
Figure~\ref{fig:target_inverse} illustrates $\Hat{u}$ and the
first moments of the target function that will be used.  We aim to
illustrate how the presented framework can be used for a class of inverse
problems. Application to realistic inverse problems will require deeper
investigations into handling noisy and incomplete data.

\subsection{Determination of the source function using cost functional
\texorpdfstring{$\mathcal{J}_{1}$}{J1}}

We consider the case with $\beta = 0$, which permits the decoupled
solution of problems at collocation points. We therefore apply both the
stochastic Galerkin and collocation methods in this section.

\subsubsection{Deterministic target}
\label{ssec:inverse_det}

As a first case, we consider a stochastic inverse problem with
the cost functional $\mathcal{J}_{1}$ and with the target~$\hat{z}$
taken as the mean of the forward problem. This represents an inverse
problem where only one observation of the system response is available;
some stochastic data computed in the forward problem has been discarded.
Recall that the stochastic properties of $\kappa$ are known.

\begin{figure}
\centering
\subfloat[mean state]{\includegraphics[width=0.32\textwidth]{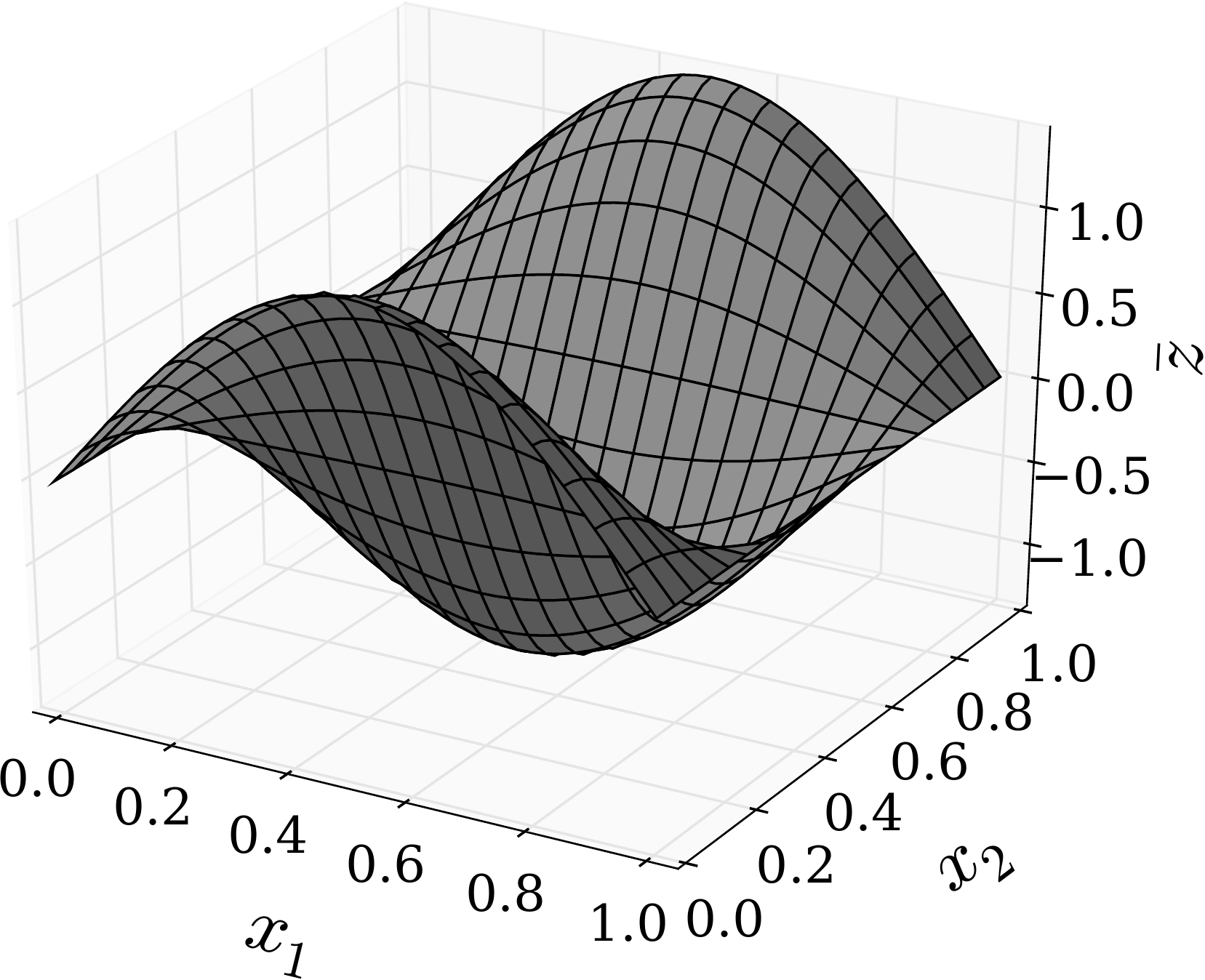}}
\subfloat[variance state]{\includegraphics[width=0.33\textwidth]{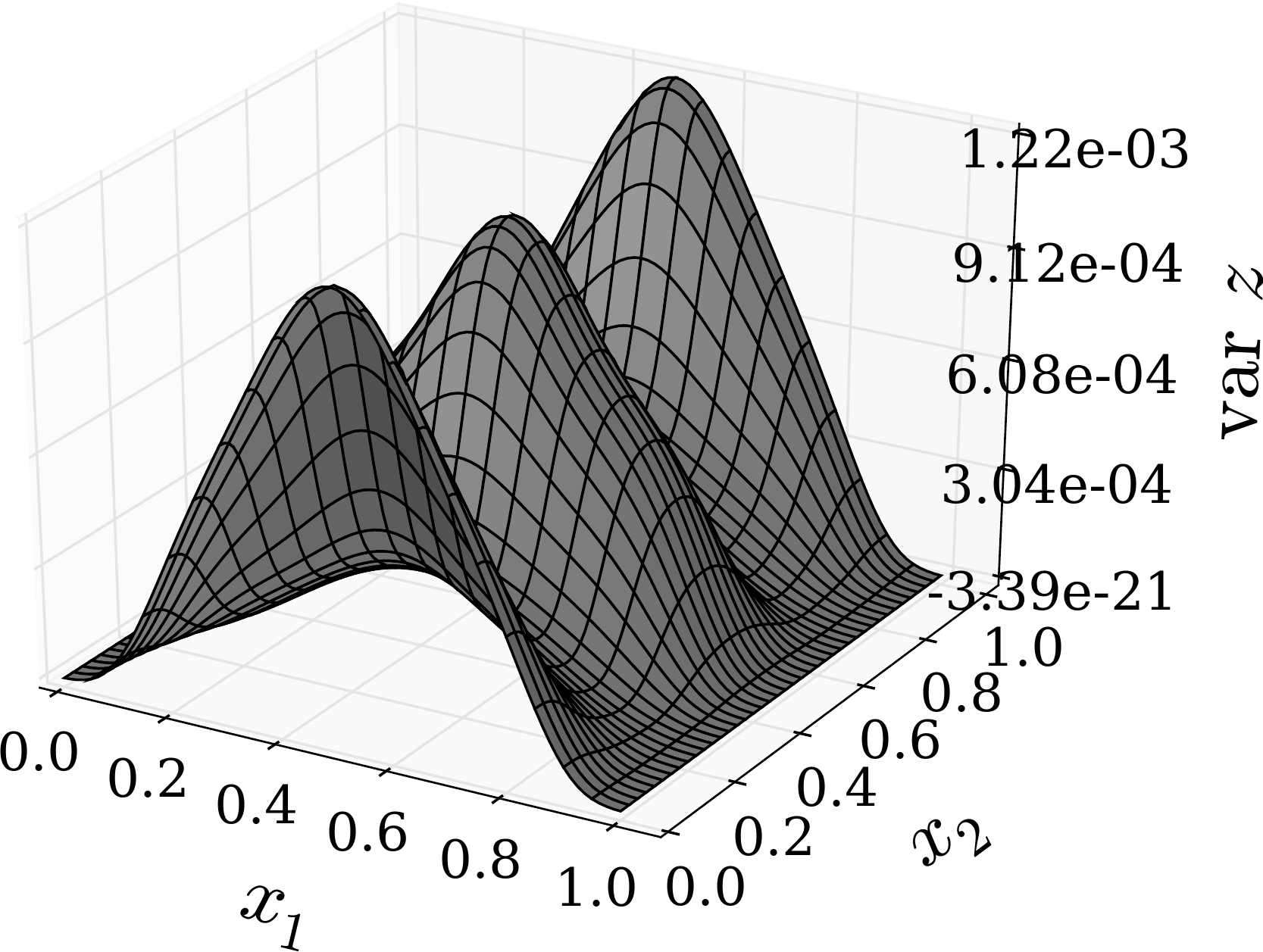}}
\\
\subfloat[mean source]{\includegraphics[width=0.32\textwidth]{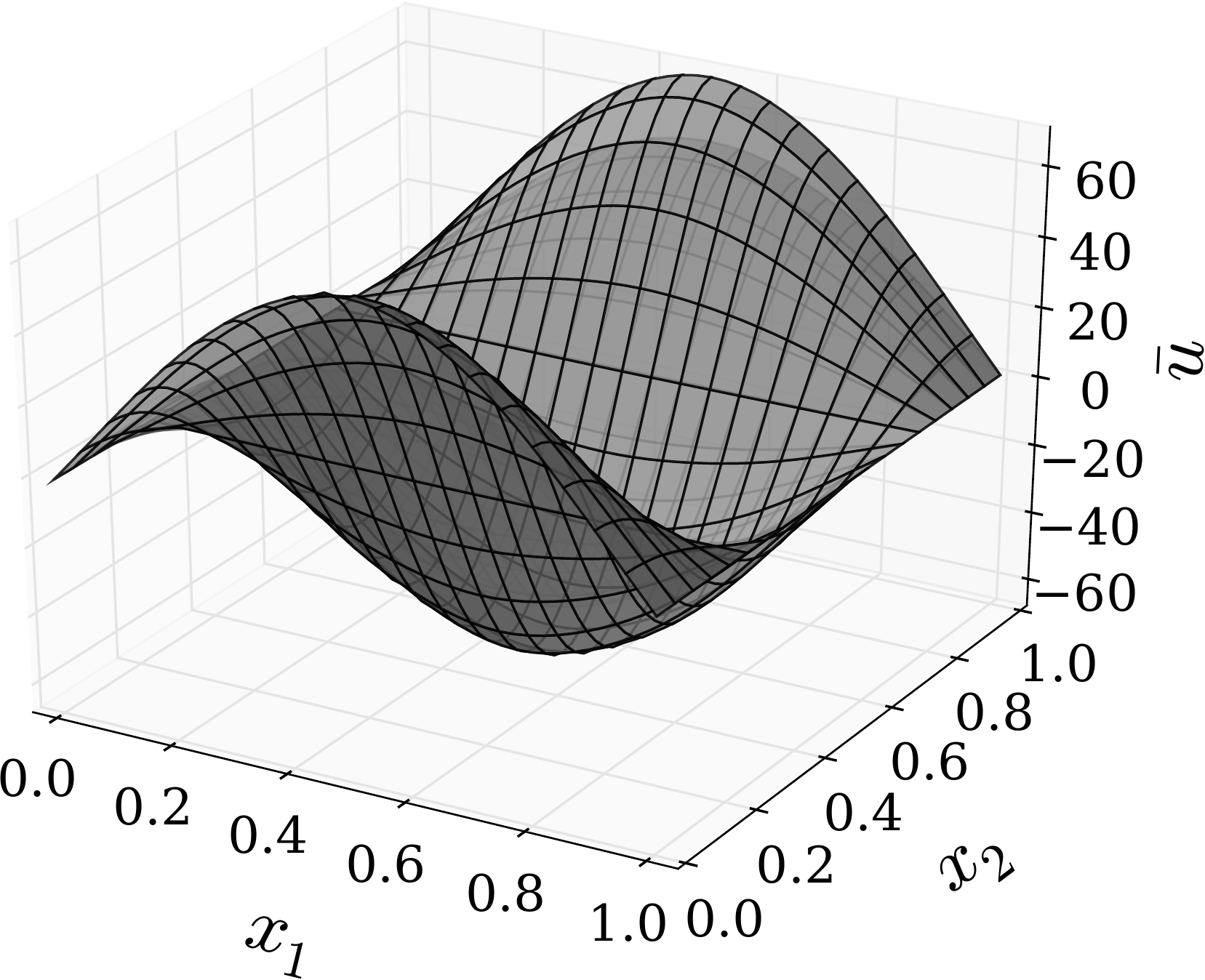}}
\subfloat[variance source]{\includegraphics[width=0.33\textwidth]{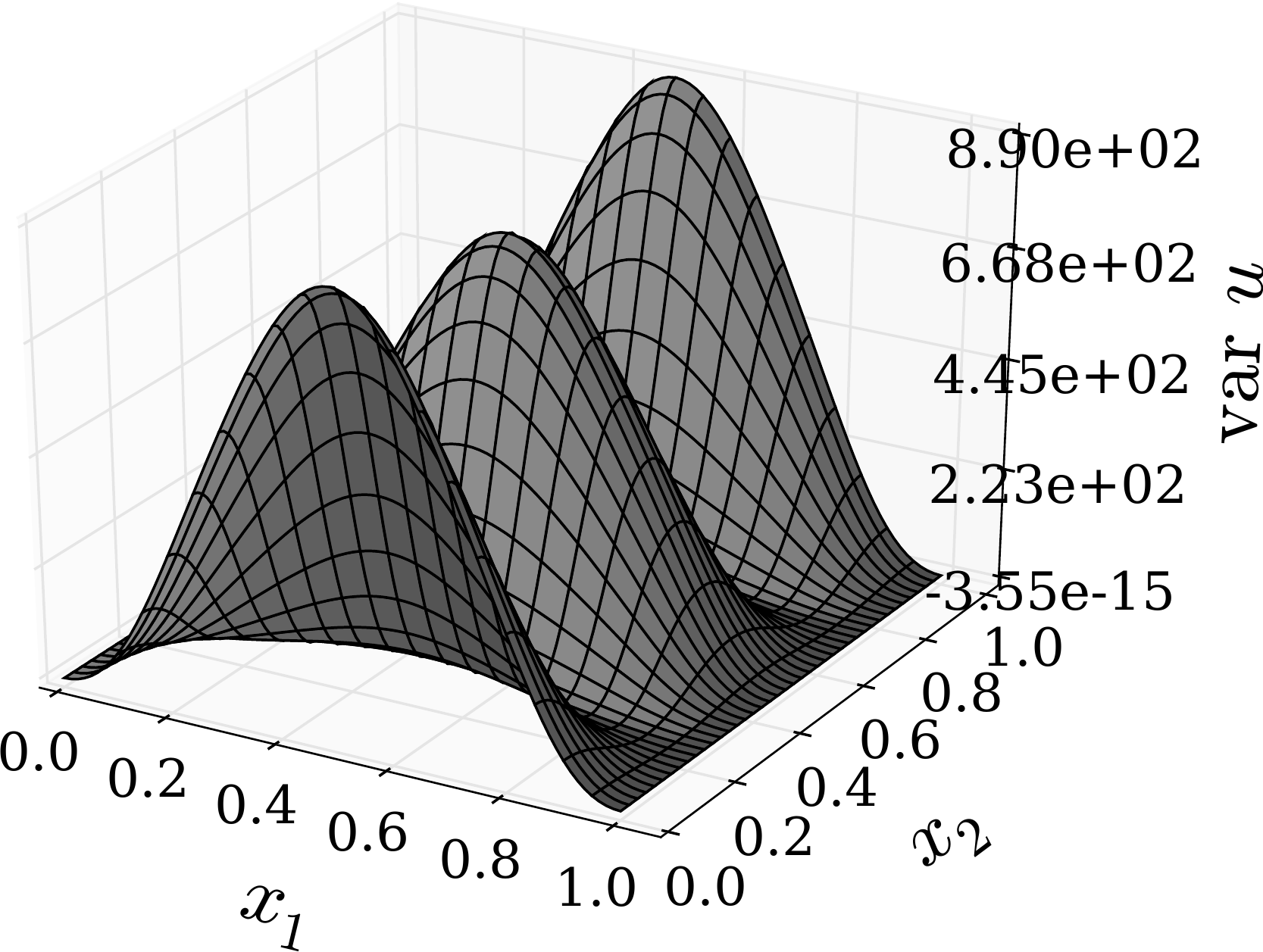}}
\caption{Mean and variance of the state and source variable for an
  inverse problem associated with the cost functional $\mathcal{J}_{1}$
  and computed with the stochastic collocation method with $\alpha = 1$,
  $\beta = \delta = 0$ and $\gamma = 10^{-5}$. The mean of the solution
  to the forward problem is used as the target. The deterministic source
  $\hat{u}$ is illustrated transparently in (c) for reference.}
\label{fig:inverse_det_coll5}
\end{figure}

\begin{figure}
 \centering
  \includegraphics[width=0.5\textwidth]{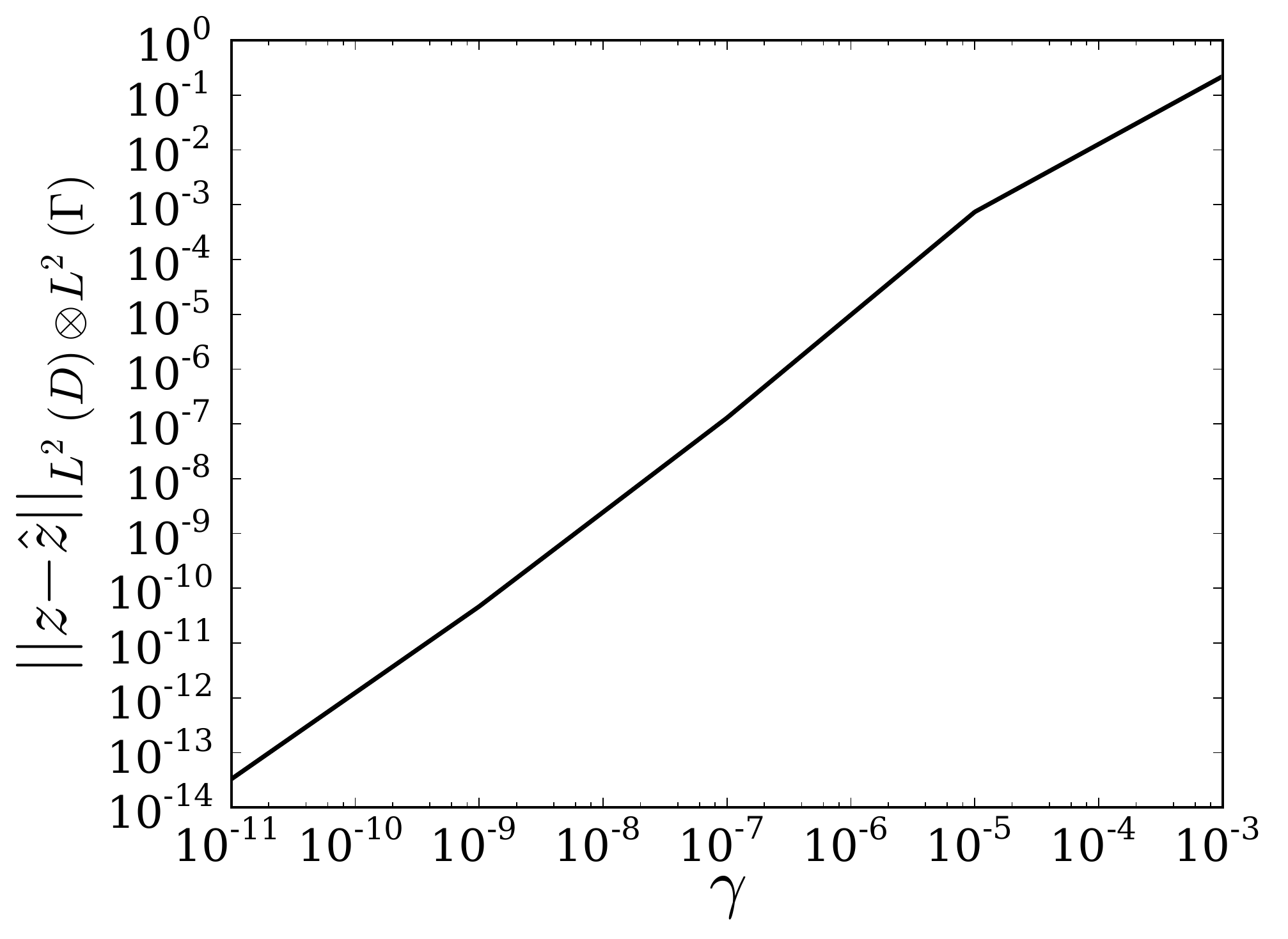}
  \caption{Tracking error $\norm{z - \hat{z}}_{L^2(D)\otimes
    L^2_\rho(\Gamma)}$ as a function of the penalty parameter~$\gamma$
    for the inverse problem considered in Section~\ref{ssec:inverse_det}.}
\label{fig:gamma_inverse}
\end{figure}

The computed first moments of the state and source function, using
$\alpha = 1$, $\beta = \delta =0$ and $\gamma = 10^{-5}$, are shown
in Fig.~\ref{fig:inverse_det_coll5}.  The computed cost functional and
tracking error are presented in Table~\ref{tab:cost_inverse}.  The mean
of the target $\Hat{z}$ and the mean of the state variable $\Bar{z}$
coincide visually, whereas there is a considerable difference between the
actual source $\Hat{u}$ and the mean of the computed source~$\Bar{u}$.
A large variance of the computed source term relative to its mean is also
observed.  This is inherent to the posed problem as limited observation
data is being used.  As $\gamma$ approaches zero, the tracking error
is observed to approach zero, as shown in Fig.~\ref{fig:gamma_inverse}.
This figure contrasts the tracking error in Fig.~\ref{fig:gamma} for a
control problem, in which case the considered target cannot be reached
as the target is not a solution of the forward problem.

\subsubsection{Stochastic target}

\begin{figure}
\centering
\subfloat[mean state]{\includegraphics[width=0.32\textwidth]{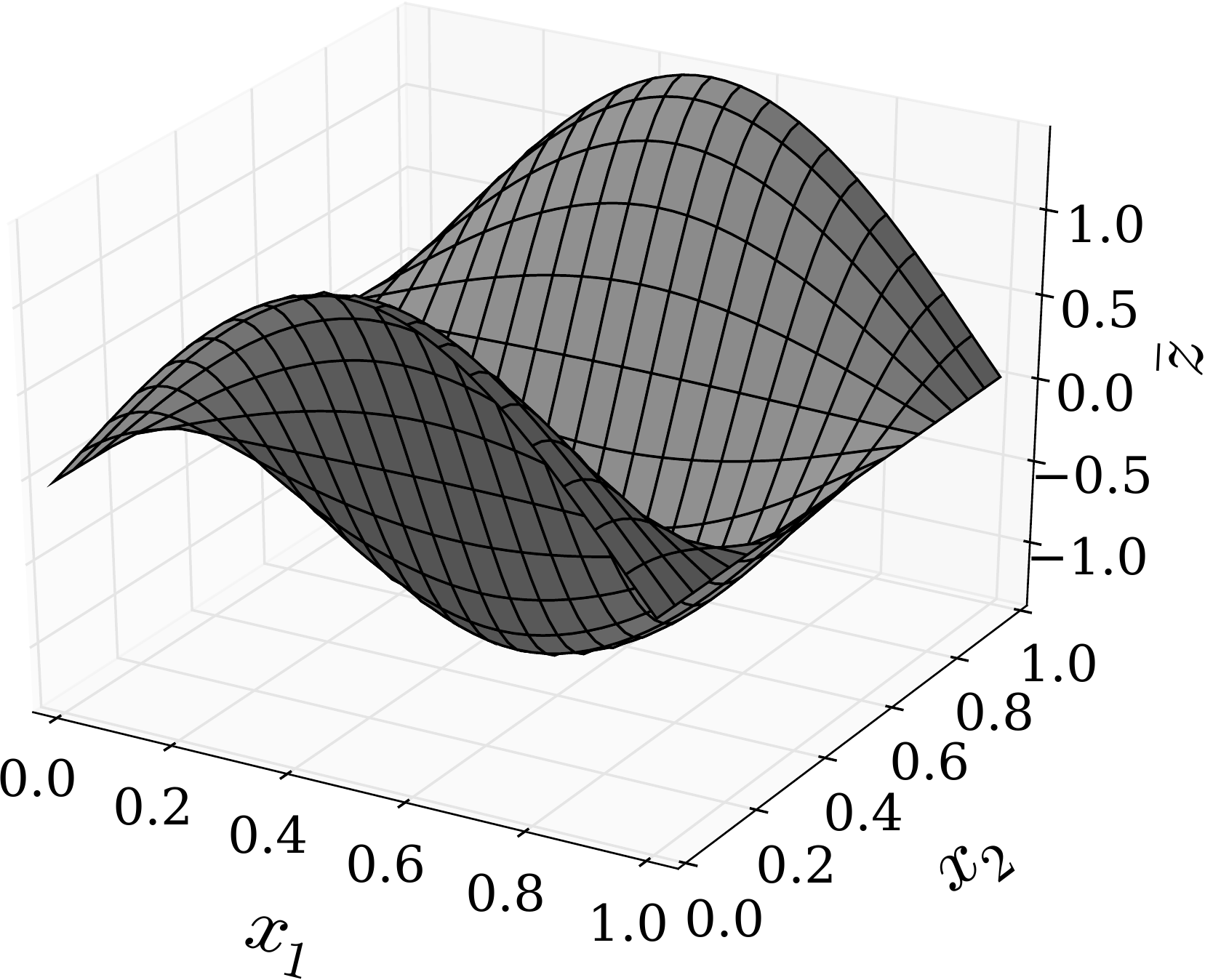}}
\subfloat[variance state]{\includegraphics[width=0.32\textwidth]{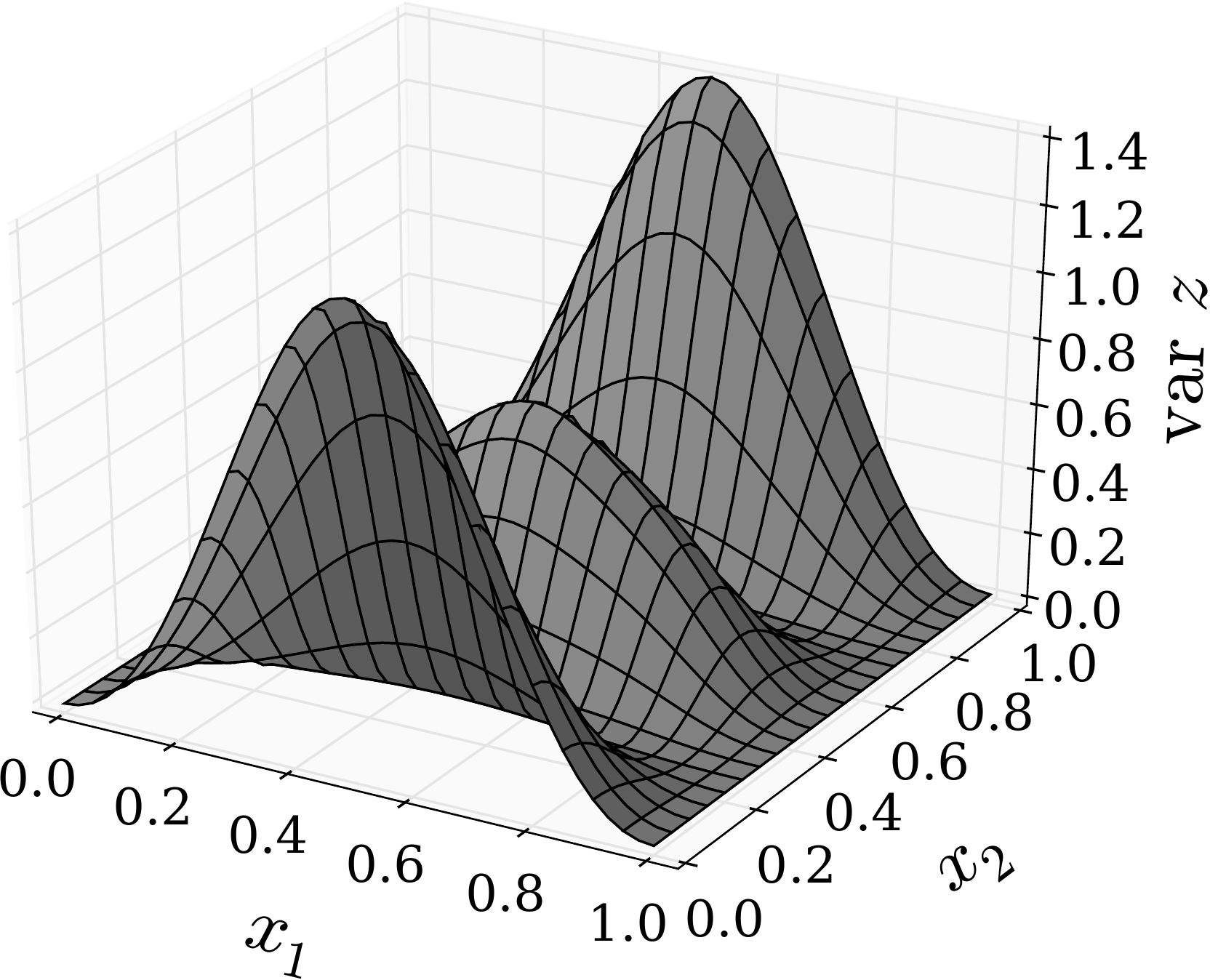}}
\\
\subfloat[mean source]{\includegraphics[width=0.32\textwidth]{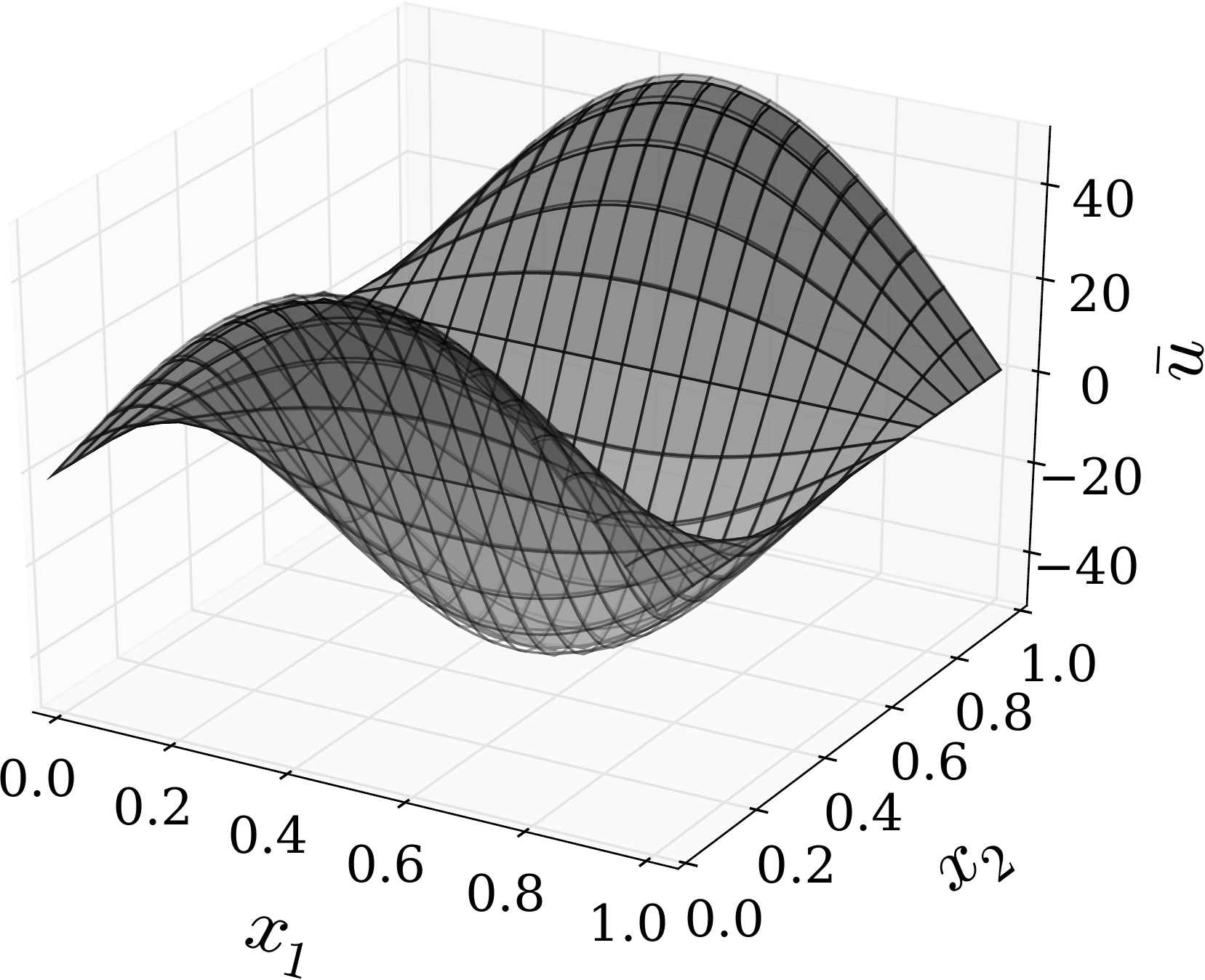}}
\subfloat[variance source]{\includegraphics[width=0.32\textwidth]{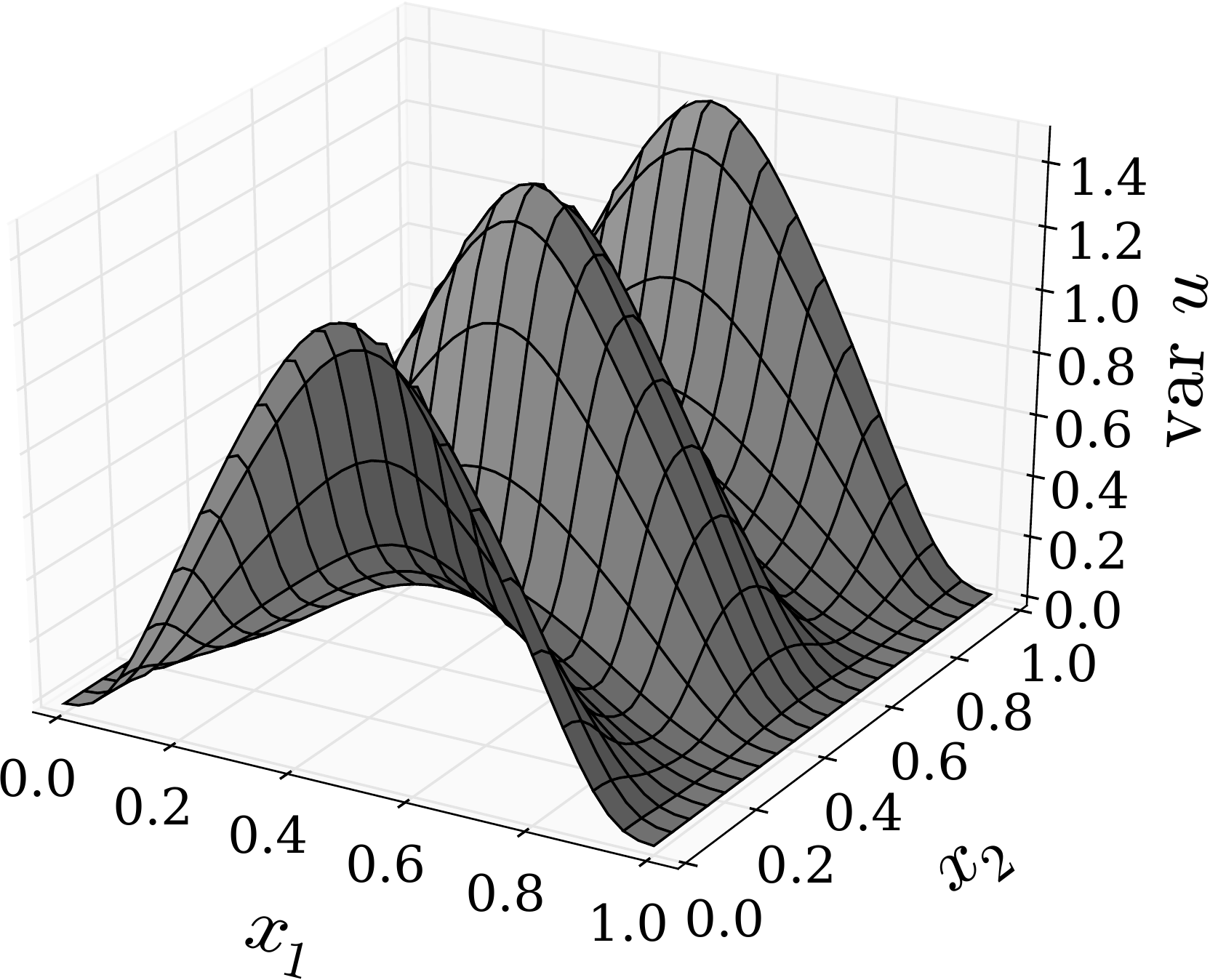}}
\caption{Mean and variance of the state and source variable for an inverse
  problem associated with the cost functional $\mathcal{J}_{1}$ and computed
  with the stochastic collocation method with $\alpha = 1$, $\beta = \delta = 0$
  and $\gamma = 10^{-5}$. The mean and variance of the target~$\hat{z}$
  are illustrated in Fig.~\ref{fig:target_inverse}. The deterministic
  source $\hat{u}$ is illustrated transparently in (c) for reference.}
\label{fig:stoch_dist_control5_A}
\end{figure}

\begin{figure}
\centering
\subfloat[mean state]{\includegraphics[width=0.32\textwidth]{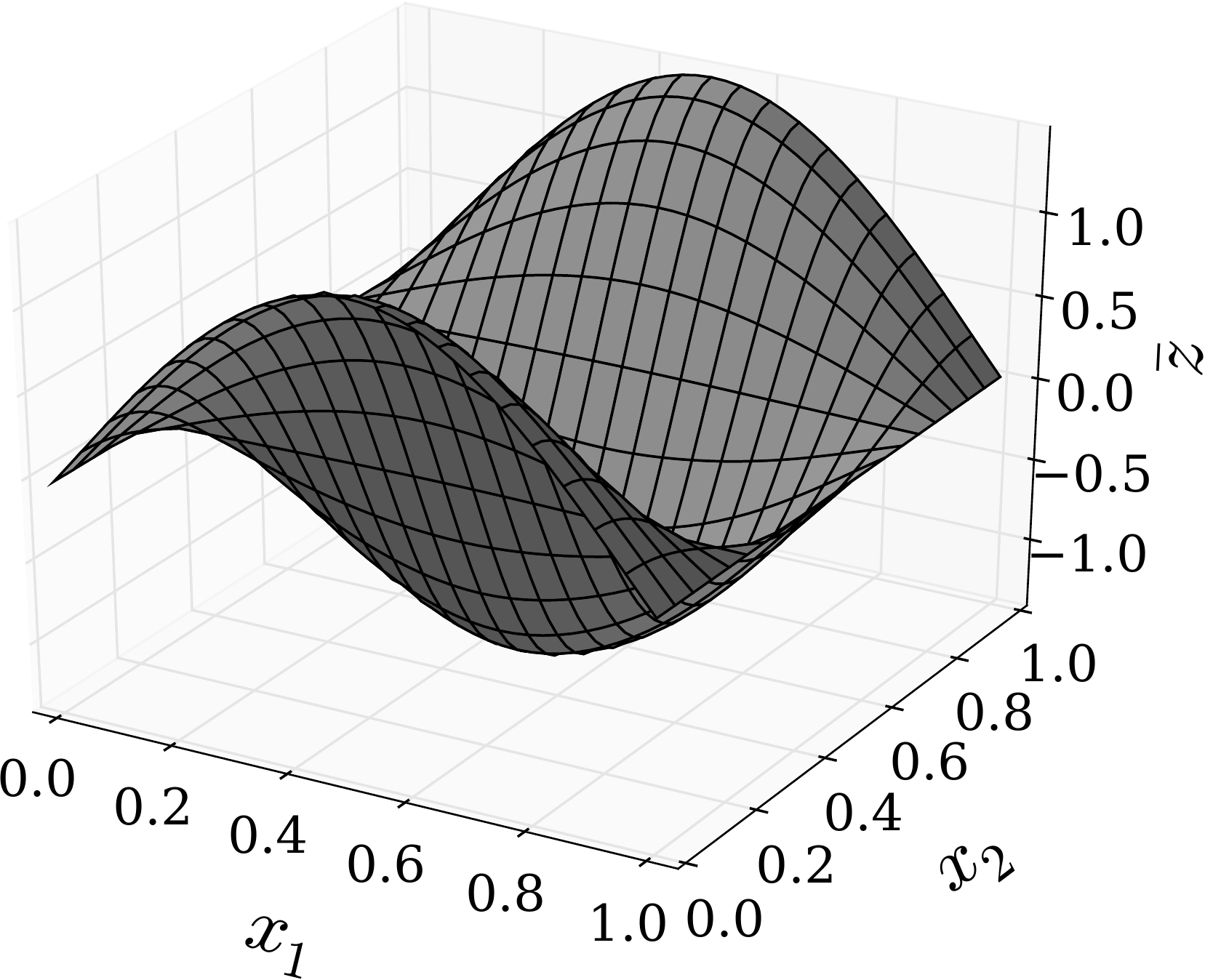}}
\subfloat[variance state]{\includegraphics[width=0.32\textwidth]{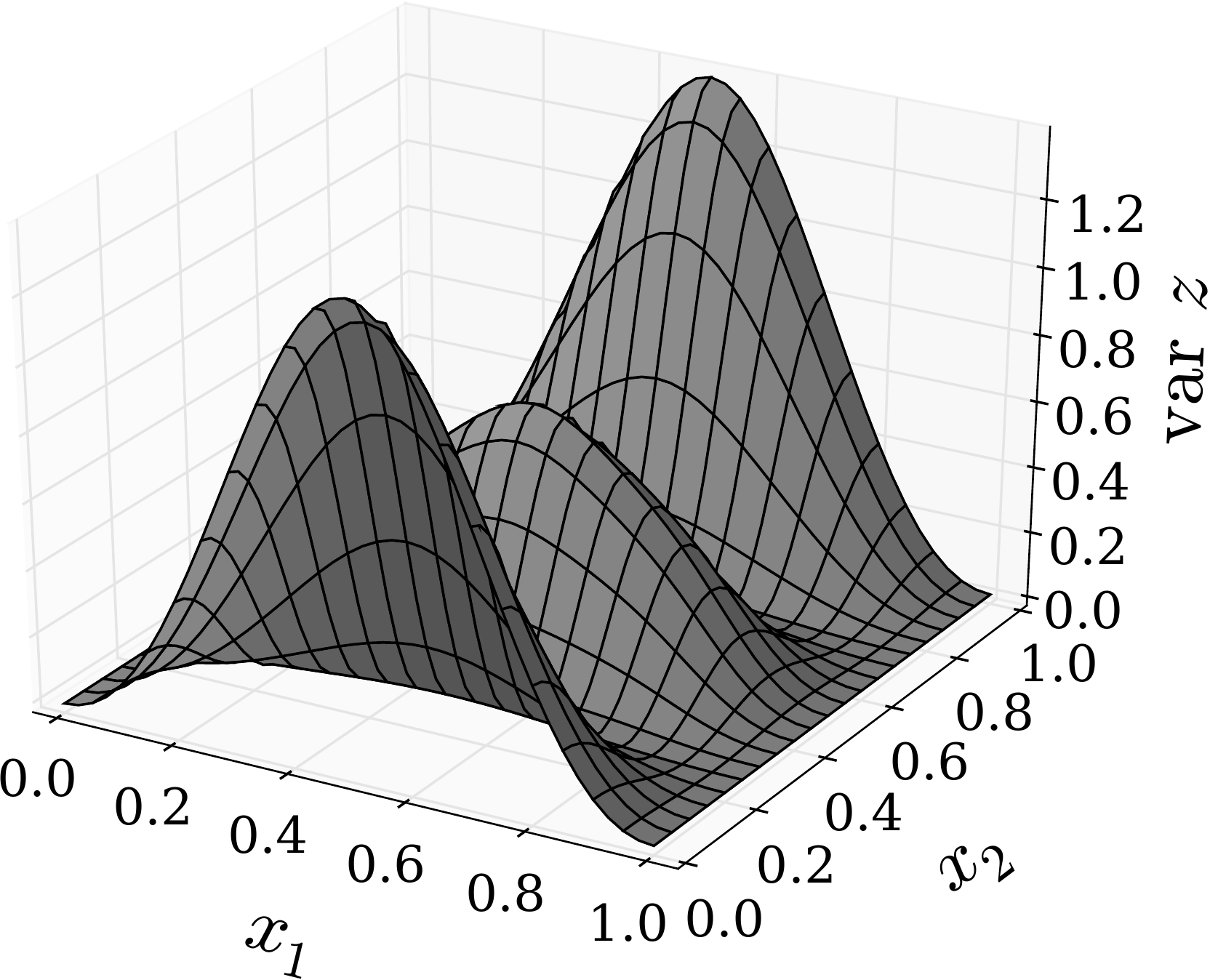}}
\\
\subfloat[mean source]{\includegraphics[width=0.32\textwidth]{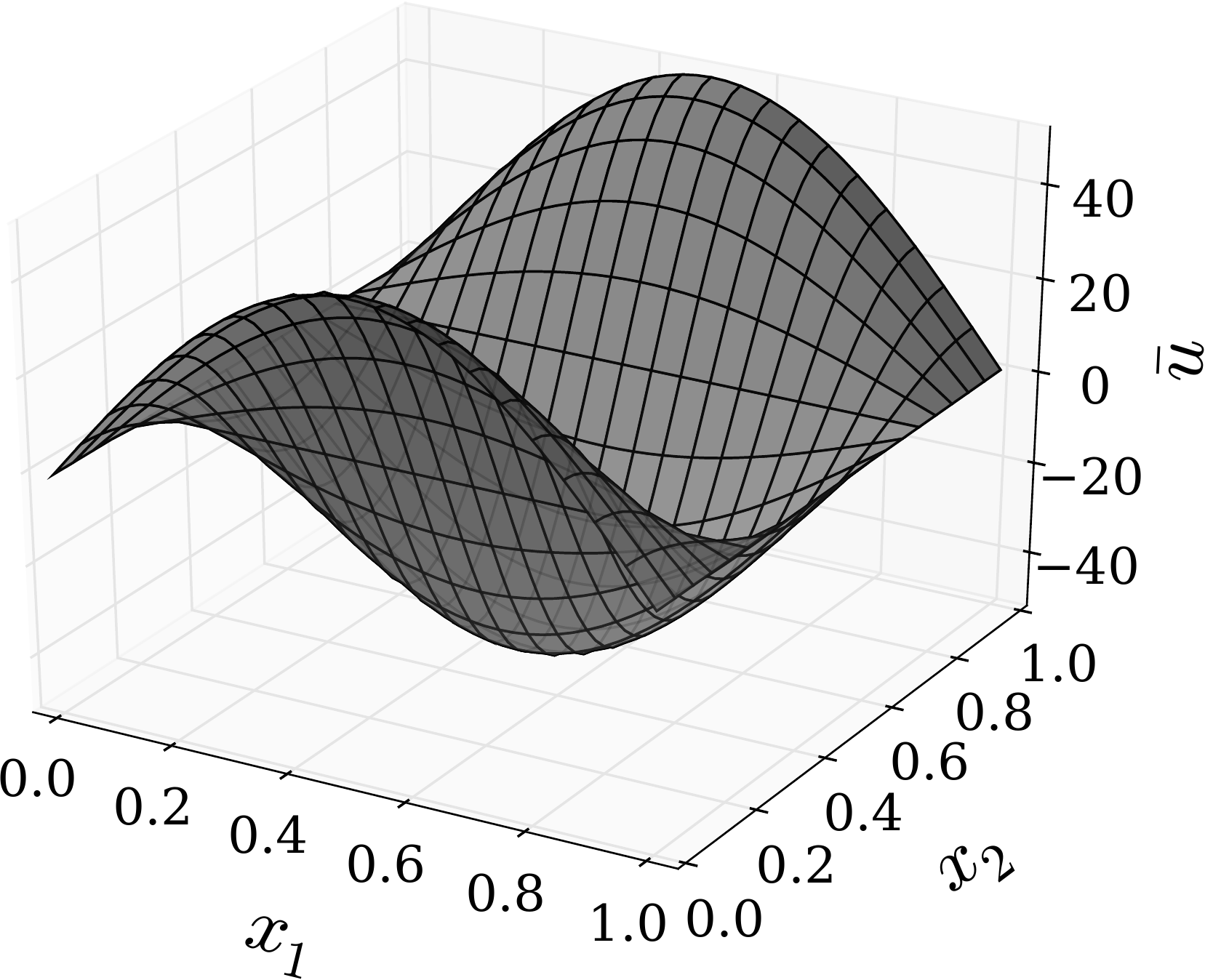}}
\subfloat[variance source]{\includegraphics[width=0.33\textwidth]{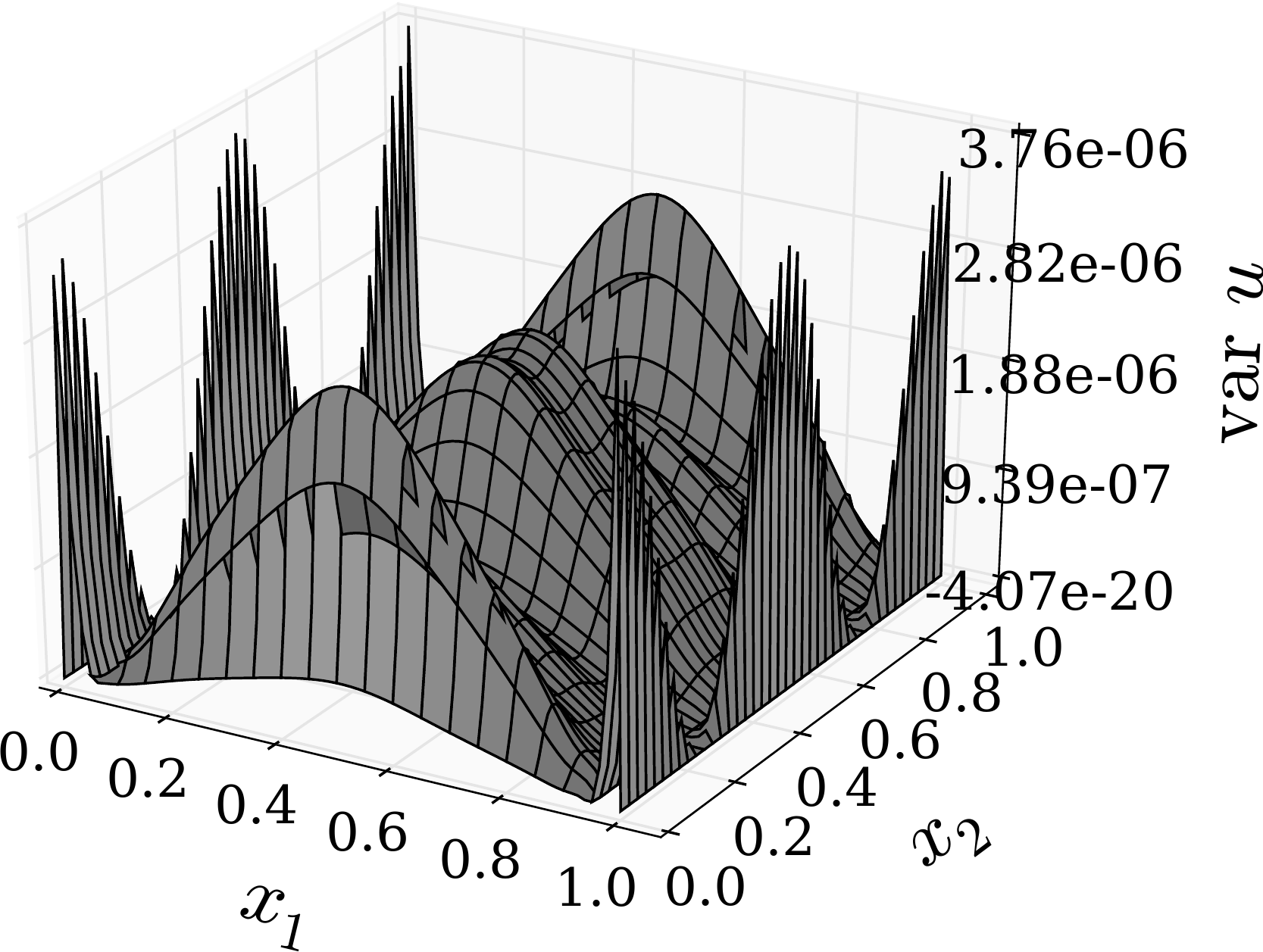}}
\caption{Mean and variance of the state and source variable for an
  inverse problem associated with the cost functional $\mathcal{J}_{1}$
  and computed with the stochastic collocation method with $\alpha = 1$,
  $\beta = \delta = 0$ and $\gamma = 10^{-8}$. The mean and variance of
  the target~$\hat{z}$ are illustrated in Fig.~\ref{fig:target_inverse}. The
  deterministic source $\hat{u}$ is illustrated transparently in (c) for
  reference. }
\label{fig:stoch_dist_control8_A}
\end{figure}

A stochastic target in now considered, which is the complete solution of
the forward problem. Hence, no data has been discarded from the forward
problem.  The computed mean of the state variable, the variance of the
state variable, the mean of the computed source and the variance of
the source function are shown in Fig.~\ref{fig:stoch_dist_control5_A}
for the case $\alpha = 1$, $\beta = 0$, $\delta = 0$ and $\gamma =
10^{-5}$. The results were computed using the collocation method.
The statistics of the state variable and the target visually coincide. A
measure of the quality of the approximation is quantified by the
tracking error in Table~\ref{tab:cost_inverse}. As the penalty
parameter $\gamma$ is decreased the computed state variable
better matches the stochastic target. This effect is illustrated in
Fig.~\ref{fig:stoch_dist_control8_A}, for which $\gamma = 10^{-8}$,
and in Table~\ref{tab:cost_inverse}.  The computed source term in
Fig.~\ref{fig:stoch_dist_control8_A} matches the exact solution $\Hat{u}$ very
well in both the mean and the variance (which is zero for $\Hat{u}$).
The stochastic Galerkin method yields similar results, as indicated in
Table~\ref{tab:cost_inverse} where various quantities computed with the
two methods are summarised.

For small penalty parameters, limited control over the source
function is imposed, which can lead to a more expensive iterative
solution of the one-shot systems. In the stochastic Galerkin
case, the convergence rate of the mean-based preconditioner (see
Section~\ref{ssec:meanbased}) deteriorates severely for small values
of~$\gamma$. The collective multigrid method shows robust convergence
behaviour. A similar observation on the computational complexity was
made by \citet{Zabaras2011}.  Also, observe the non-smooth variance
of the source function on the domain boundaries for $\gamma = 10^{-8}$
in Fig.~\ref{fig:stoch_dist_control8_A}(d).

\subsection{Determination of the source function using cost functional
\texorpdfstring{$\mathcal{J}_{2}$}{J2}}

\begin{figure}
\centering
\subfloat[mean state]{\includegraphics[width=0.32\textwidth]{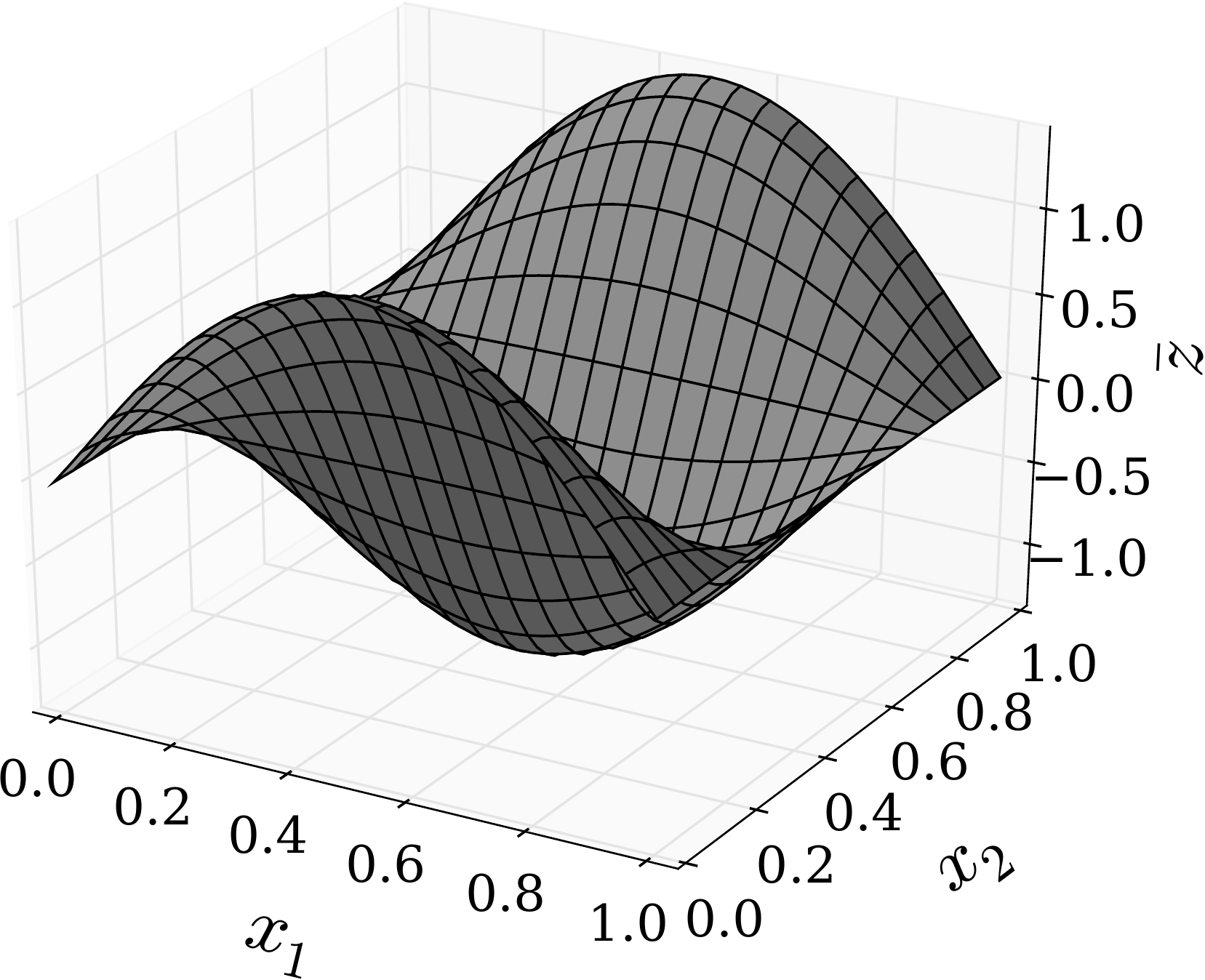}}
\subfloat[variance state]{\includegraphics[width=0.32\textwidth]{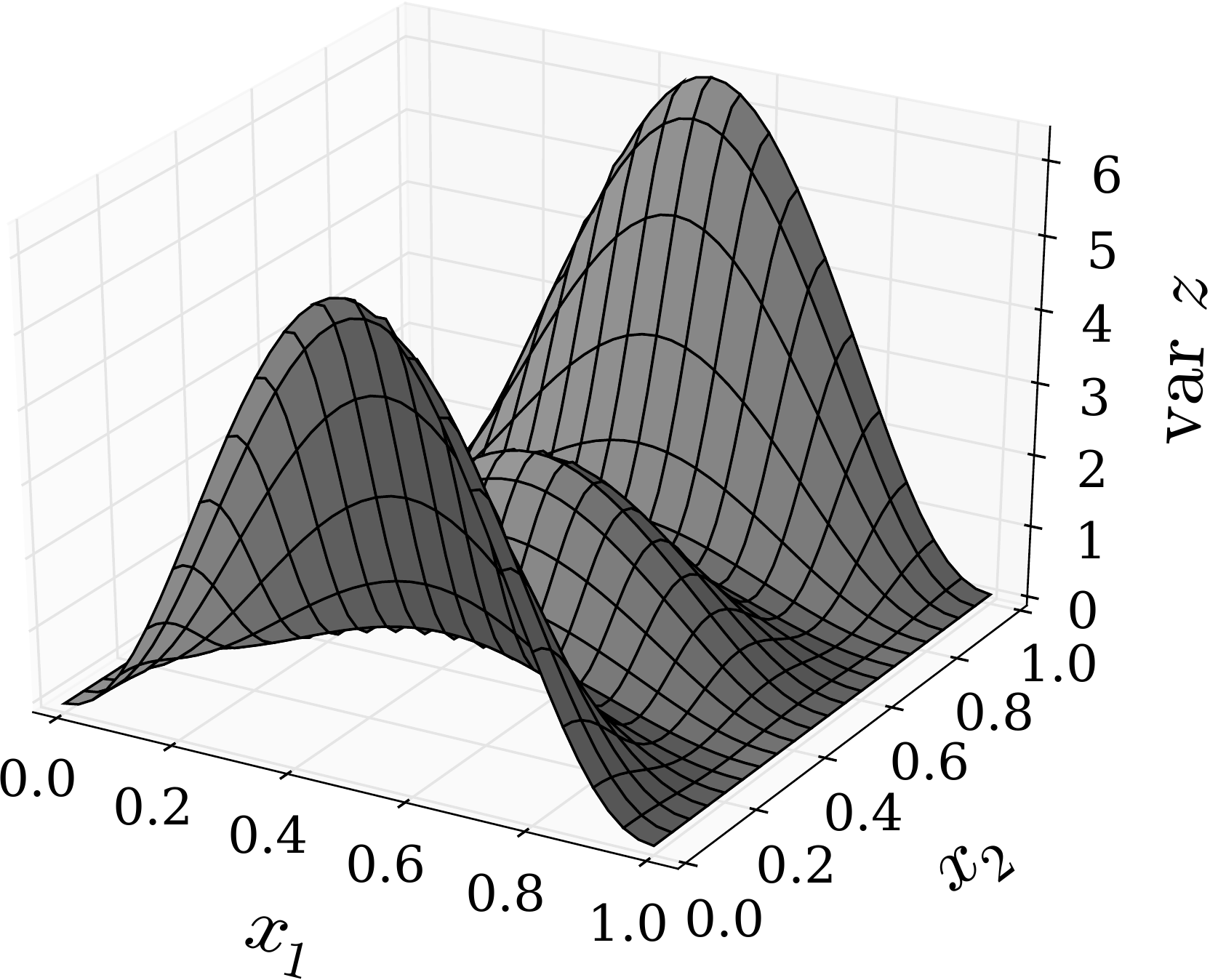}}
\\
\subfloat[mean source]{\includegraphics[width=0.32\textwidth]{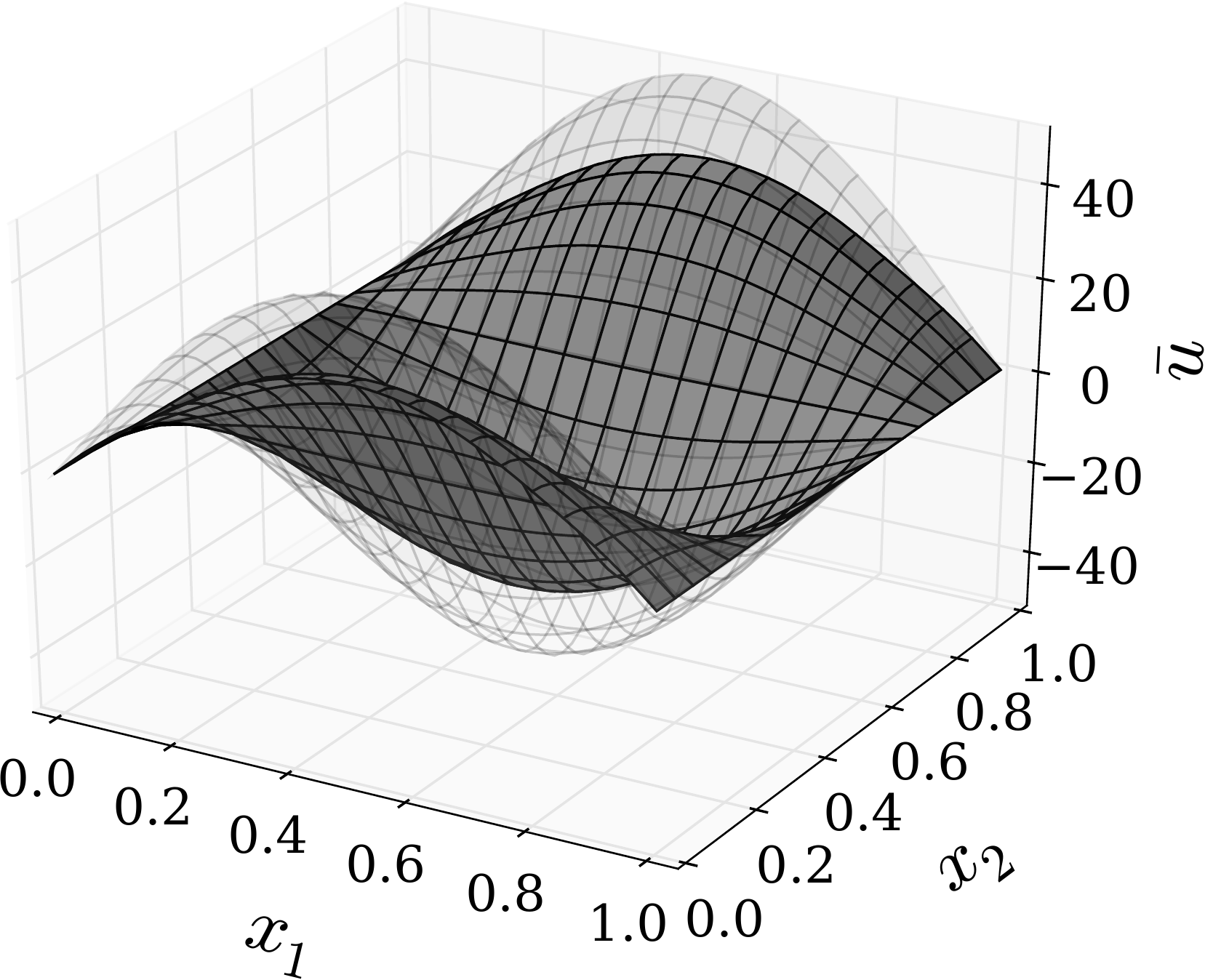}}
\subfloat[variance source]{\includegraphics[width=0.33\textwidth]{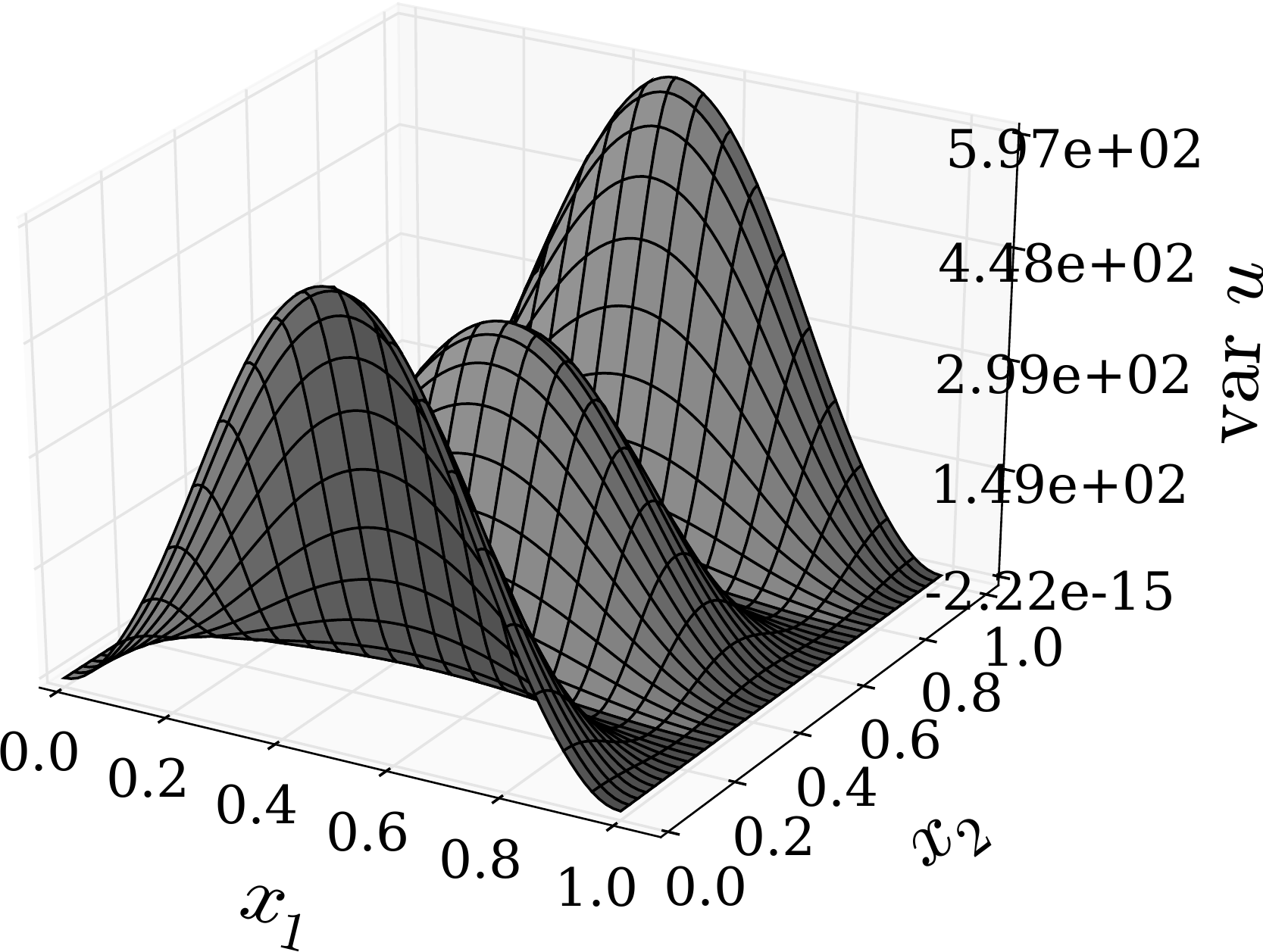}}
\caption{Mean and variance of the state and source variable for an
  inverse problem associated with the cost functional $\mathcal{J}_{2}$
  and computed with the stochastic collocation method with $\alpha = 1$,
  $\beta = \delta = 0$ and $\gamma = 10^{-5}$. The mean and variance of
  the target~$\hat{z}$ are illustrated in Fig.~\ref{fig:target_inverse}.
  The deterministic source $\hat{u}$ is illustrated transparently in (c)
  for reference.}
\label{fig:stoch_dist_inverse2}
\end{figure}

We mirror the stochastic inverse problem considered
in Section~\ref{ssec:inverse_det}, but now for the cost
functional~$\mathcal{J}_{2}$.  The $\mathcal{J}_{2}$ formulation involves
the mean of an unknown field, which means that a stochastic collocation
formulation will not lead to decoupled problems when~$\alpha \neq 0$.
Fig.~\ref{fig:stoch_dist_inverse2} shows the computed mean and variance
of the state and source functions, computed with the stochastic Galerkin
method. As in Fig.~\ref{fig:inverse_det_coll5}, the mean of the state
variable and the mean target visually coincide. Since the variance of the
state variable does not contribute to $\mathcal{J}_{2}$, a larger variance
is observed than in Fig.~\ref{fig:inverse_det_coll5}(b). The computed
source function clearly differs (note the magnitudes) from that computed
when using the cost functional $\mathcal{J}_1$. The corresponding
values of the cost functional and tracking error are summarised in
Table~\ref{tab:cost_inverse}.

\begin{table}
\centering \small
\begin{tabular}{|l|cr@{$\ \times\ $}lr@{$\ \times\ $}l|}
\hline
\multicolumn{1}{|c}{}&\multicolumn{5}{c|}{}\\[-0.3cm]
\multicolumn{1}{|c}{}& $\mathcal{J}(z,u)$ &
\multicolumn{2}{c}{$\norm{z - \hat{z}}^2_{L^2(D)\otimes L^2_\rho(\Gamma)}$}&
\multicolumn{2}{c|}{$e_{u}$}\\[0.1cm]
\hline
\multicolumn{1}{|c}{}&\multicolumn{5}{c|}{}\\[-0.3cm]
\multicolumn{6}{|l|}{Stochastic Galerkin finite element solution
($ Q = 36$)}\\
\hline
\multicolumn{1}{|c|}{}&\multicolumn{5}{c|}{}\\[-0.3cm]
cost functional $\mathcal{J}_{1}$, $\hat{z}$ deterministic,
 $\gamma = 10^{-5}$ & $ 6.786\times 10^{-3}$ & $7.225$ & $10^{-4}$ &
 4.368& $10^{-1}$\\
cost functional $\mathcal{J}_{1}$, $\hat{z}$ stochastic,
$\gamma = 10^{-5}$ &$3.035 \times 10^{-3}$ & $1.678$&$ 10^{-4}$
& 1.505 & $10^{-3}$\\
cost functional $\mathcal{J}_{1}$, $\hat{z}$ stochastic,
 $\gamma = 10^{-8}$&$3.123\times 10^{-6}$ & $1.882$ &$ 10^{-10}$
& 2.339& $10^{-9}$\\
cost functional $\mathcal{J}_{2}$,  $\gamma = 10^{-5}$ &
$2.107\times 10^{-3}$ & $3.784 $ & $10^{-5} $&3.193& $10^{-1}$\\
cost functional $\mathcal{J}_{2}$,  $\gamma = 10^{-8}$ &
$2.127 \times 10^{-6}$ & $9.833$&$ 10^{-11}$&3.190& $10^{-1}$ \\ 
\hline
\multicolumn{1}{|c}{}&\multicolumn{5}{c|}{}\\[-0.3cm]
\multicolumn{6}{|l|}{Stochastic collocation finite element solution
($Q = 141$)}\\
\hline
\multicolumn{1}{|c|}{}&\multicolumn{5}{c|}{}\\[-0.3cm]
cost functional $\mathcal{J}_{1}$, $\hat{z}$ deterministic,
 $\gamma = 10^{-5}$ &
 $6.957\times 10^{-3}$&$ 7.406$&$ 10^{-4}$ & $4.556$ &$ 10^{-1}$
\\
cost functional $\mathcal{J}_{1}$, $\hat{z}$ stochastic,
 $\gamma = 10^{-5}$ &
$3.035\times 10^{-3}$ & $1.678$&$ 10^{-4}$&  $1.506$&$ 10^{-3}$ \\
cost functional $\mathcal{J}_{1}$, $\hat{z}$ stochastic,
 $\gamma = 10^{-8}$&
$3.123\times 10^{-6}$ & $1.882 $&$ 10^{-10}$& $2.334$&$ 10^{-9}$ \\ 
\hline
\end{tabular}
\caption{Summary of the cost functional, tracking error of the state variable
  and relative error $e_{u} :=\norm{u - \hat{u}}^2_{L^2(D)\otimes
 L^2_\rho(\Gamma)}/\norm{\hat{u}}^2_{L^2(D)} $  of the computed source for the
 considered inverse problems with unknown source function.}
\label{tab:cost_inverse}
\end{table}

\section{Conclusions}
\label{sec:conclu}

A one-shot solution approach for stochastic optimal control problems
with PDE constraints has been presented. The problem was formulated
as an optimisation problem constrained by a stochastic elliptic PDE,
and the framework is sufficiently general to also address a class of
inverse problems that involve uncertainty.  Statistical moments have
been included in the cost functional and uncertainty in the controller
response has been accounted for.  To compute solutions, a one-shot method
is combined with stochastic finite element discretisations.  It is shown
that the non-intrusivity property of the stochastic collocation method
is lost when moments of the state variable appear in the cost functional,
or when the control function is a deterministic function.  We argue that
the control function must contain a deterministic component in order for
the problem to constitute a control problem, versus a stochastic inverse
problem, hence the non-intrusivity property of the stochastic collocation
method will be lost for one-shot formulations.  Applying a stochastic
Galerkin method does not impose additional difficulties compared to
solving stochastic PDEs, hence, for the method presented in this work the
stochastic Galerkin method is preferred over the collocation method for
control problems.  In the case of inverse problems, where the function
to be found is wholly stochastic, it was shown that it is possible in
some cases to preserve the non-intrusivity property of the stochastic
collocation method.

The formulated methods are supported by extensive numerical experiments
that address both optimal control and inverse problems.  The computed
results illustrate the impact of various options in the formulation and
the difference between the considered cost functionals. In particular,
examples show the impact, for the considered model problem, of the
different ways in which statistical data can be included in the cost
functionals.
\subsection*{Acknowledgement}

This work was performed while E.~Rosseel was a visiting researcher at
the University of Cambridge with support from the Research Foundation
Flanders (Belgium).
\appendix
\section{Evaluating derivatives of stochastic functions}
\label{app}

Analytical expressions for inner products of the form
\begin{equation}\label{eq:int}
 \int_\Gamma \nabla_y \psi_i(y) \cdot \nabla_y\psi_j(y) \rho \dif y
\end{equation}
are derived here based on the orthogonality properties of the generalised
polynomial chaos basis functions~$\psi_i$ in \eqref{eq:psi}.
Given~\eqref{eq:psi} and $\delta$ the Kronecker delta, the
integral~\eqref{eq:int} can be rewritten as
\begin{align}\label{eq:int2}
  \int_\Gamma \nabla_y \psi_i(y) \cdot \nabla_y \psi_j(y) \rho \dif y
    = \sum_{k=1}^L \prod_{t = 1, t\neq k}^L \delta_{i_t,j_t}
\int_{\Gamma_k} \frac{d \varphi_{i_k}}{dy_k}\frac{d \varphi_{j_k}}{dy_k}
 \rho \dif y_k.
\end{align}

\paragraph{Hermite polynomials.}
In the case of Hermite polynomials, which are normalised with respect
to a standard normal distribution, the derivative of $\varphi_{i_k}$
is given by
\begin{equation}
 \frac{d \varphi_{i_k}}{dy_k} = \sqrt{i_k} \varphi_{i_k - 1}.
\end{equation}
Expression~\eqref{eq:int2} then simplifies to
\begin{equation}
 \int_\Gamma \nabla_y \psi_i(y)\cdot \nabla_y \psi_j(y) \rho \dif y
  = \delta_{i,j} \sum_{k = 1}^L j_k.
\end{equation}

\paragraph{Legendre polynomials.}
In the case of Legendre polynomials, normalised and scaled to a uniform
distribution on $\sbr{-\sqrt{3},\sqrt{3}}$, it can be shown that the
following relation holds:
\begin{equation}
 \frac{d \varphi_{i_k}}{dy_k}
    = \frac{\sqrt{3 \del{2i_k + 1}}}{3}
  \sum_{t = 0}^{\left\lfloor \frac{j_k-1}{2} \right\rfloor}
        \sqrt{2\del{i_k - 1 - 2t} + 1} \varphi_{i_k - 1 - 2t}.
\end{equation}
After some calculations, expression~\eqref{eq:int2} corresponds to
\begin{multline}
 \int_\Gamma \nabla_y \psi_i(y) \cdot \nabla_y\psi_j(y) \rho \dif y
  = \sum_{k=1}^L \del{1 -\del{\del{i_k+j_k} \bmod 2}}
\frac{\sqrt{2i_k+1}\sqrt{2j_k+1}}{3}
\\
\left\lfloor \frac{\min\del{i_k, j_k} + 1}{2}
 \right\rfloor \del{  1- 2 \left\lfloor \frac{1-\min\del{i_k,j_k}}{2}
 \right\rfloor } \prod_{t=1, t\neq k}^L \delta_{j_t, i_t}.
\end{multline}
\bibliographystyle{abbrvnat}
\bibliography{references.bib}
\end{document}